# THE DEVELOPMENT OF A HYBRID ASYMPTOTIC EXPANSION FOR THE HARDY FUNCTION $Z(t)$, CONSISTING OF JUST $[2\sqrt{2}-2]\sqrt{t/2\pi}$ MAIN SUM TERMS, SOME 17% *LESS* THAN THE CELEBRATED RIEMANN-SIEGEL FORMULA

D. M. Lewis

*Department of Mathematics, University of Liverpool, M&O Building, Peach St, Liverpool L69 7ZL*

(Tel: *+44 151 794 4014*. Email: *D.M.Lewis@liverpool.ac.uk*)

## Abstract

This paper begins with a re-examination of the Riemann-Siegel Integral, which first discovered amongst by Bessel-Hagen in 1926 and expanded upon by C. L. Siegel on his 1932 account of Riemann's unpublished work on the zeta function. By application of standard asymptotic methods for integral estimation, and the use of certain approximations pertaining to special functions, it proves possible to derive a new zeta-sum for the Hardy function $Z(t)$. In itself this new zeta-sum (whose terms made up of elementary functions, but are unlike those that arise from the analytic continuation of the Dirichlet series) proves to be a computationally inefficient method for calculation of $Z(t)$. However, by further, independent analysis, it proves possible to correlate the terms the new zeta-sum with the terms of the Riemann-Siegel formula, thought, since its discovery by Siegel, to be the most efficient means of calculating $Z(t)$. Closer examination of this correlation reveals that is possible to formulate a hybrid asymptotic formula for $Z(t)$, consisting of a sum containing both Riemann-Siegel terms and terms from the new zeta-sum, in such a way as to reduce the overall CPU time required by a factor between $14-15\%$. Alongside the obvious computational benefits of such a result, the very existence of the new zeta-sum itself highlights new theoretical avenues of study in this field.

## 1. **Introduction**

With the publication of Siegel (1932, [33]) it became clear that Riemann's researches on the zeta function and analytic number theory went far beyond those published in his remarkable paper "On the Number of Primes Less than a Given Magnitude" of 1859 [14]. In particular it was clear that Riemann had devised an extremely sophisticated method of calculating the zeta function $\zeta\left(\frac{1}{2}+\mathrm{i}t\right)$ high inside the critical strip ($|\mathrm{Im}\ t| \leq 1/2\,, |\mathrm{Re}\ t| \gg 1$), subsequently known as the Riemann-Siegel formula far in advance of the classical method based on Euler-Maclaurin summation employed by [18] to find the first 10 or so non-trivial zeros lying (as predicted by Riemann's Hypothesis) along the critical line. (See [14, 23] for the historical background to the subject.) The Riemann-Siegel formula (abbreviated subsequently to RS formula) is an asymptotic approximation to the Hardy function $Z(t) = e^{\mathrm{i}\theta(t)}\zeta\left(\frac{1}{2}+\mathrm{i}t\right)$ with $t \in \mathbb{R}$, consisting of a main sum (7) and an $O\left(t^{-1/4}\right)$ correction (A75). In conjunction with subsequent developments ([3], [4], [6], [15], [23], [24], [26], [32] & [34]), it represents the fastest method currently known for checking the presence of zeros along the critical line, has formed the basis of almost all large scale computations of the zeta function in the critical strip ([7], [19], [28], [30] & [31]) since its discovery in 1932.

In addition to this result, [33] also published an alternative representation of the zeta function, first discovered in 1926 by Bessel-Hagen in Riemann's *Nachlass*, in terms of a definite integral. This representation, now termed as the Riemann-Siegel integral (RSI) formula is given by ([14] p166)

$$\frac{2\xi(s)}{s(s-1)} = F(s) + \overline{F(1-\bar{s})}, \tag{1}$$

where the function $F$ is defined by the formula

$$F(s) = \Gamma\left(\frac{s}{2}\right)\pi^{-s/2}\int_{0\searrow 1}\frac{e^{-\mathrm{i}\pi z^2}z^{-s}}{e^{\mathrm{i}\pi z}-e^{-\mathrm{i}\pi z}}\,dz = \Gamma\left(\frac{s}{2}\right)\pi^{-s/2}\mathrm{RSI}\,. \tag{2}$$

Here $\xi(s)$ is the xi-function, defined in terms of the zeta function (originally by Riemann, although see note on p16 of [14]) by

$$\xi(s) = \Gamma\left(\frac{s}{2}+1\right)(s-1)\pi^{-s/2}\zeta(s)\,, \tag{3}$$

and the symbol $0\searrow 1$ denotes a path of integration along a line of slope $-1$ crossing the real axis between 0 and 1 (the integrand has poles at integer values) and directed from upper left to lower right. [33] and subsequently [25], adopting different methods, used (1-2) to prove the functional equation of the zeta function, viz.

$$\Gamma\left(\frac{s}{2}\right)\pi^{-s/2}\zeta(s) = \Gamma\left(\frac{1-s}{2}\right)\pi^{-(1-s)/2}\zeta(1-s) \Leftrightarrow \xi(s) = \xi(1-s)\,. \tag{4}$$

These proofs complement the two separate proofs of (4) presented in Riemann's original paper.

However, it is interesting to speculate that in devising the representation (1-2), Riemann's main interest may have been to devise more efficient methods for *calculating* the zeta function in the critical strip, given that it is known he was trying to locate zeros, and that the utility of even the RS formula would have been limited even for relatively modest $t$ without access to early to mid 20th century computing resources. Whatever the truth of that, superficially at least integral (2) does have certain attractions in that, as [14] states, "$e^{-i\pi z^2}$ approaches zero very rapidly as $|z| \to \infty$ along any line of the form 0 ↘ 1", which might make it amenable to estimation by local methods. Hence the motivation of this paper to try and answer the question as to whether it is possible obtain an asymptotic approximation for (2) which could be used to calculate $\zeta\left(\frac{1}{2} + \mathrm{i}t\right)$ high along the critical line? Unsurprisingly, given the fundamental role occupied by the zeta function in the solution of unresolved problems relating to prime numbers, this idea is not new and was first pursued by Turing shortly after publication of Siegel's researches (see [37]). Turing's method essentially consisted of integrating (2) around a semi infinite parallelogram with infinite sides $0 \searrow 1$ and $N \nwarrow N + 1$ where $N = \lfloor (t/2\pi)^{1/2} \rfloor$. This results in a formula with a similar main sum to the RS formula (7) but with an additional error term arising from the integral along the line $N \nwarrow N + 1$. Turing was able to estimate this error term more precisely for $t$ in the range 50-1000, than the corresponding remainder terms as then specified by Siegel. However, subsequent improvements in the estimates of these remainder terms have rendered Turing's method redundant and often forgotten. Recently though Turing's methodology has been revived and developed by [26], resulting in a more sophisticated estimate for Turing's error term by using incomplete gamma functions.

In this paper a new approach will be adopted in which Turing's oblique method will be replaced by a more direct attempt to estimate (2) directly along the line $0 \searrow 1$. Somewhat surprisingly, this apparently simple idea eventually yields an approximation from the RSI to Hardy's function $Z(t)$ based on an entirely new zeta-sum (see 6 & A64), whose terms are structurally different to the terms that arise from the corresponding Dirichlet series (on which the RS formula is based). *Even more surprisingly,* it subsequently proves possible to relate this new series directly to the terms of the main sum of RS formula, in such a way as to show that the latter is actually a computationally *inefficient* method for calculating $Z(t)$ when $t$ is large. Instead it is possible to formulate a more efficient hybrid formula for $Z(t)$ based on partial sums of both series. Specifically, the simplest variant of the proposed hybrid formula (63) is made up of the first $\left\lfloor \sqrt{t/4\pi} \right\rfloor$ terms of the RS main sum (7), combined with the following asymptotic approximation, devised from the new series, for the sum of the remaining terms, viz.

$$\sum_{N=\lfloor\sqrt{t/4\pi}\rfloor}^{\lfloor\sqrt{t/2\pi}\rfloor} \frac{2cos\left\{\frac{t}{2}log\left(\frac{t}{2\pi}\right)-\frac{t}{2}-\frac{\pi}{8}-tlog(N)\right\}}{\sqrt{N}}$$

$$= \sum_{\alpha=INT_O(a)+2}^{INT_O\left(3\sqrt{t/\pi}\right)} \frac{2\sqrt{2}\cos\left((t/2)\{log(pc)+1/pc\}+t/2+\pi/8\right)}{(\alpha^2-a^2)^{\frac{1}{4}}} + T_{a+\varepsilon} + error \quad (5)$$

where $a=\sqrt{8t/\pi}$, $pc(\alpha)=2(\alpha/a)^2\left\{1+\sqrt{1-(a/\alpha)^2}\right\}-1$ and $\alpha$ represents successive *odd integers*. (Here $INT_O(x)\equiv$ largest odd integer $\le x$.) The transition term $T_{a+\varepsilon}$, which is only significant if $|\varepsilon|=|a-NINT_O(a)|\le t^{-1/6}$, is discussed in sections A3.6-3.7 (see A65), and the maximum $|error|<6.15t^{-1/12}$. Notice that the sum on the right-hand side of (5) is made up of *purely elementary functions* and is *confined purely* to the odd integers *alone*. Consequently it consists of only $\frac{1}{2}\left(3-2\sqrt{2}\right)\sqrt{t/\pi}$ elementary terms, approximately 58.5% less than the $\left(1/\sqrt{2}-1/2\right)\sqrt{t/\pi}$ terms needed to sum the left-hand side. Overall, this means that the number of terms which must be summed to compute $Z(t)$ using the hybrid formula is only $\left[2\sqrt{2}-2\right]\sqrt{t/2\pi}$, an almost 17.16% reduction compared to the $\sqrt{t/2\pi}$ terms needed using the RS formula alone. Sample computations of the errors that arise when utilising (5) are shown in Table I (Section 3). The small scale of these errors is indicative of the approximation's effectiveness. Although the amplitudes and phases making up the new sum are consist of elementary functions, they are slightly more expensive to calculate than those making up the RS formula. This means that the typical saving in terms of CPU time when computing the RHS of (5) *vis-a-vis* the LHS is about 50% (not the 58.5% based purely on the count of the number of terms). Overall this means that employing the hybrid formula to compute $Z(t)$ results in an efficiency saving in terms of CPU time of some $14-15\%$ (see Section 3.2 for details of a very simple algorithm which achieves these savings).

The objective of this paper will be to explain in detail how the hybrid formula, encapsulated by approximation (5), arises. To this end the paper will be split up into two, effectively distinct, parts. Appendices A & C are devoted to the derivation of the new zeta-sum, of which the right-hand side of (5) forms just a fraction, from the Riemann-Siegel integral (2). The analysis in these Appendices relies on two unproven assertions (summarised in Section A3.8) concerning the convergence properties of this new zeta-sum, but very strongly suggests that it provides an asymptotic approximation $Z(t)$. The first part of the main paper (Section 2 and Appendix B) is devoted to a formal proof of this suggestion, encapsulated by the *Main Theorem* (see Section 2.1) which states that the difference between the main terms of the new zeta-sum and the main terms of the RS formula will tend to zero as $t\to\infty$. In the process this circumvents the need for a proof of the interpretations made in deriving the new zeta-sum. In the second part of the main paper (Section 3) the analysis employed to complete the proof of the *Main Theorem,* is then used to formulate the hybrid (5, 63), together with an examination of the CPU time savings which can be attained.

## 2. Relating the new zeta-sum for $Z(t)$ to the terms in the RS Formula.

2.1 *Main Theorem*

In the new series (A64) devised for Hardy's $Z$ function in Appendix A, the first order main zeta-sum that requires calculation is given by

$$MS(t) = \sum_{\alpha=INT_O(a)+2}^{N_\alpha} \frac{2\sqrt{2}\cos\big((t/2)\{log(pc) + 1/pc\} + t/2 + \pi/8\big)}{(\alpha^2 - a^2)^{\frac{1}{4}}}, \qquad (6)$$

where $a = \sqrt{8t/\pi}$, $\alpha$ are odd integers, $pc(\alpha) = 2(\alpha/a)^2\left\{1 + \sqrt{1-(a/\alpha)^2}\right\} - 1 > 1$ and $N_\alpha \geq INT_O\ (t/\pi) + 2$ (see section A4). The primary goal of this section (alongside Appendix B) will be to establish a framework leading to the proof (cf. 2.6) of the following theorem.

*Main Theorem*

Let $t > t_0 = 30$. Define $N_t = \lfloor (t/2\pi)^{1/2} \rfloor$, let $\theta(t)$ be the Riemann-Siegel theta function (A34) with expansion $(tlog(t/2\pi) - t)/2 - \pi/8 + O(t^{-1})$, and suppose $MS(t)$ represents the sum defined by (6) above. Then there exists an upper bound $\mathrm{E}(t) > 0$, satisfying the $\lim_{t\to\infty} \mathrm{E}(t) = 0$, such that

$$\left| 2\sum_{N=1}^{N_t} \frac{cos\{\theta(t) - tlog(N)\}}{\sqrt{N}} - MS(t) \right| < \mathrm{E}(t). \qquad (7)$$

The sum on the left hand side of (7) forms the dominant contribution to the full RS formula (A75). In the process of establishing a proof of this theorem, it will become apparent how one can utilise (6) in order to establish the more computational efficient method (5) for the calculation of $Z(t)$. This is discussed in detail in Section 3.

2.2 *Initial Comments on a Proof of the Main Theorem*

The minimum value of the upper limit, $N_\alpha = INT_O\ (t/\pi) + 2$ in (6), is dictated by the constraints of Euler-Maclaurin summation technique discussed in section A4. However, this choice of cut-off requires the calculation of many terms in the Bernoulli sum (A72) to guarantee convergence, and a better procedure is to set $N_\alpha = INT_O(2t/\pi - a)$, approximately double the minimum value. As will be seen this increased number of terms presents no extra difficulty regarding the proof of the main theorem, whilst at this cut-off value series (A72) converges very rapidly after just a few terms. In what follows it will be assumed that $a \geq INT_O(a) + O\big(t^{-1/6}\big)$. In those instances when $a$ happens to equal or lie very close to an odd integer, then series (6) must be augmented by the transitional term $T_{a+\varepsilon}$ calculated either from (A65), or by numerical integration as discussed in Section A3.6. The modifications needed to account for the transitional term and complete the proof of the main theorem are discussed in Section 2.6.

From the residue theorem, series (6) is given by

$$Re\left[\frac{\sqrt{2}}{\pi \mathrm{i}}\int_{C(R)}\frac{(-\pi/2)sin(\pi z/2)e^{\mathrm{i}t/2+\mathrm{i}\pi/8}exp\{(\mathrm{i}t/2)\left(log\left(pc(z)\right)+1/pc(z)\right)\}}{cos(\pi z/2)(z^2-a^2)^{1/4}}dz\right]. \quad (8)$$

Here $z = t/\pi + Re^{\mathrm{i}\theta}$, with $\theta \in [0, 2\pi]$, is a complex variable mapping out the circumference of a circle $C(R)$, centred at $t/\pi$ and radius $R$. Near the simple poles at $z = \alpha$ lying inside this circle, $cos(\pi z/2) \approx -(\pi/2)sin(\pi\alpha/2)(z-\alpha)$, so that the real parts of the residues conform to the terms in (6). The choice of $R$ should be slightly less than $t/\pi - a$. Some very specific considerations concerning the estimation of certain integrals necessary for the evaluation of (8) (see Appendix B), suggest that a suitable definition for the maximum of $R$ should be

$$R_{max} = \frac{t}{\pi} - \frac{a(pc_- + 1)}{2\sqrt{pc_-}}, \quad (9)$$

where

$$pc_- = 1 + \eta\frac{\sqrt{8\pi}}{t^{1/3}} = 1 + \frac{1}{2}\left(\frac{1250}{\pi^3}\right)^{1/6}\frac{\sqrt{8\pi}}{t^{1/3}} \Rightarrow R_{max} \approx \frac{t}{\pi} - a - \frac{\eta^2\sqrt{8\pi}}{t^{1/6}}. \quad (10)$$

The following analysis concentrates on the case when $R \approx_p R_{max}$. Here the sign $\approx_p$ denotes equality $=$, except on those occasions when the radius must be adjusted slightly to ensure the contour path does not cross the real axis through a pole of (8). The methodology, with a few specified modifications, also applies for any choice of $R = R_{max} - R_t$ with $R_t \sim O\left(t^{0-1/2}\right)$, so that any sub-series of (6) with lower limit $\alpha > INT_O(a) + 2$ can also be estimated.

Substituting $\theta$ for $z$ by means of $z = t/\pi + Re^{\mathrm{i}\theta} \Rightarrow dz = R\mathrm{i}e^{\mathrm{i}\theta}d\theta$ into the integral (8) and utilising the results

$$\frac{\pi sin(\pi z/2)}{2cos(\pi z/2)}e^{\mathrm{i}t\left(log\left(pc(z)\right)+1/pc(z)\right)/2} = \frac{\pi}{2\mathrm{i}}e^{\mathrm{i}t\left(log\left(pc(\theta)\right)+1/pc(\theta)\right)/2}\begin{cases}\dfrac{e^{\mathrm{i}\pi z(\theta)}-1}{e^{\mathrm{i}\pi z(\theta)}+1} & \text{for } \theta \in (0,\pi)\\[2ex] \dfrac{1-e^{-\mathrm{i}\pi z(\theta)}}{1+e^{-\mathrm{i}\pi z(\theta)}} & \text{for } \theta \in (\pi,2\pi)\end{cases} \quad (11)$$

and

$$(z^2 - a^2)^{1/4} = \sqrt{\frac{t}{\pi}}\left[\left(1 + \frac{\pi R}{t}e^{\mathrm{i}\theta}\right)^2 - \frac{8\pi}{t}\right]^{1/4}, \quad (12)$$

transforms (8) into

$$-\frac{R\sqrt{\pi}e^{\mathrm{i}t/2+\mathrm{i}\pi/8}}{\mathrm{i}\sqrt{2t}}\left\{\int_0^{\pi}\frac{e^{\mathrm{i}\theta+\mathrm{i}t\left(log\left(pc(\theta)\right)+1/pc(\theta)\right)/2}}{[(1+\pi Re^{\mathrm{i}\theta}/t)^2 - 8\pi/t]^{1/4}}\frac{e^{\mathrm{i}\pi z(\theta)}-1}{e^{\mathrm{i}\pi z(\theta)}+1}d\theta\right.$$

$$\left.+\int_{\pi}^{2\pi}\frac{e^{\mathrm{i}\theta+\mathrm{i}t\left(log\left(pc(\theta)\right)+1/pc(\theta)\right)/2}}{[(1+\pi Re^{\mathrm{i}\theta}/t)^2 - 8\pi/t]^{1/4}}\frac{1-e^{-\mathrm{i}\pi z(\theta)}}{1+e^{-\mathrm{i}\pi z(\theta)}}d\theta\right\}. \quad (13)$$

Estimates for the values of these two integrals in (13) will be established separately.

2.3 *Estimation of the first integral in* (13)

*Lemma 2.3.* Provided the path of the circle $C(R)$ centred at $t/\pi$ and radius $R$ avoids crossing the real axis at a pole of integrand (8), then the first integral in (13) is of $O(t^{-1})$.

*Proof.* Along a contour defined by $z = t/\pi + \left(t/\pi - a - \varepsilon\sqrt{t/\pi}\right)e^{i\theta}$ where $\eta^2\pi\sqrt{8}/t^{2/3} \le \varepsilon \le 1$ and $\theta \in (0, 2\pi)$, the real part of the exponential phase of the integrals in (13) is given by

$$Re\left\{\frac{it}{2}\left(log\big(pc(\theta)\big) + \frac{1}{pc(\theta)}\right)\right\} = -\frac{t}{2}\left[Arg\left(\frac{\left(1+e^{i\theta}\right)^2}{2\pi}\right)\right] + \frac{\sqrt{8\pi t}\left(1+\varepsilon/\sqrt{8}\right)sin(\theta)}{2+2cos(\theta)} + O(1),$$

$$\theta \notin \left(\pi \pm O(1)/\sqrt{t}\right). \quad (14)$$

$$Re\left\{\frac{it}{2}\left(log\big(pc(\theta)\big) + \frac{1}{pc(\theta)}\right)\right\} \approx (\theta-\pi)\left[\frac{t^{\frac{3}{2}}}{\sqrt{\pi}}\left(\frac{1}{\sqrt{2}} - \frac{\sqrt{\varepsilon}}{2^{\frac{3}{4}}} + \frac{\varepsilon}{4}\right) - t\left(2 - 2^{\frac{3}{4}}\sqrt{\varepsilon} + \sqrt{2}\varepsilon\right) + O\left((t\varepsilon)^{\frac{3}{2}}\right)\right],$$

$$\theta \in \left(\pi \pm O(1)/\sqrt{t}\right). \quad (15)$$

Equation (15) is necessary because when $\theta \approx \pi$ the variable $z$ changes from $O(t)$ to $O(\sqrt{t})$. The argument term in (14) will be positive when the imaginary part of $\left(1+e^{i\theta}\right)^2 > 0$, that is when $2sin(\theta)[1+cos(\theta)] > 0$. The term in square brackets is clearly positive, so the positivity or otherwise of the argument term is simply determined by the sign of $sin(\theta)$ which, together with (15), shows that the real phase is negative over the range $\theta \in (0, \pi)$. Consequently the first integrand in (13) is exponentially small over the interior of its range and hence it is possible to estimate its size from its behaviour near its endpoints. It is also clear from a comparison of (14) and (15), that the derivatives at $\theta = 0$ and $\theta = \pi$ are of $O(t)$ and $O(t^{3/2})$ respectively. This means that the modulus of the integrand decays much more rapidly around $\theta = \pi$ than $\theta = 0$, and, as a result, the main contribution to this integral will be concentrated around $\theta = 0$. (The size of the contribution near $\theta = \pi$ is complicated by the one quarter power term in the integrand's denominator being potentially near zero and hence it depends on the specific choice of $\varepsilon$. In general it will be $O(\varepsilon^{-1/4}t^{-5/4})$, but even so, at the minimum of $\varepsilon \sim O(t^{-2/3})$ corresponding to $R_{max}$, this will always be less than the contribution at $\theta = 0$.) Now

$$\frac{it}{2}\left(log\big(pc(\theta)\big) + \frac{1}{pc(\theta)}\right) = \frac{it}{2}log(t) + \frac{it}{2}\left[log\left|\frac{\left(1+e^{i\theta}\right)^2}{2\pi}\right| + iArg\left(\frac{\left(1+e^{i\theta}\right)^2}{2\pi}\right)\right] - i\frac{\sqrt{8\pi t}\left(1+\varepsilon/\sqrt{8}\right)e^{i\theta}}{\left(1+e^{i\theta}\right)} + O(1)$$

(16)

which means that near $\theta = 0$ (and specifically including the $O(1)$ terms in 16)

$$\frac{\mathrm{i}t}{2}\left(log\big(pc(\theta)\big)+\frac{1}{pc(\theta)}\right)$$

$$=\frac{\mathrm{i}t}{2}log\left(\frac{2t}{\pi}\right)-\mathrm{i}\sqrt{2\pi t}\left(1+\frac{\varepsilon}{\sqrt{8}}\right)-\mathrm{i}\frac{\pi}{4}-\mathrm{i}\pi\left(1+\frac{\varepsilon}{\sqrt{8}}\right)^2-\left[\frac{t}{2}-\sqrt{\frac{\pi t}{2}}\left(1+\frac{\varepsilon}{\sqrt{8}}\right)\right]\theta+O(t\theta^2)$$

$$=\mathrm{i}\psi(t)-\left[\frac{t}{2}-\sqrt{\frac{\pi t}{2}}\left(1+\frac{\varepsilon}{\sqrt{8}}\right)\right]\theta+O(t\theta^2). \tag{17}$$

Hence the first integral in (13) can be approximated by

$$\int_0^{\pi}\frac{e^{\mathrm{i}\theta+\mathrm{i}t(log(pc(\theta))+1/pc(\theta))/2}}{[(1+\pi Re^{\mathrm{i}\theta}/t)^2-8\pi/t]^{1/4}}\frac{e^{\mathrm{i}\pi z(\theta)}-1}{e^{\mathrm{i}\pi z(\theta)}+1}d\theta\approx f(0)e^{\mathrm{i}\psi(t)}\int_0^{O(1)}e^{-[t/2-\sqrt{\pi t}(1+\varepsilon/\sqrt{8})/\sqrt{2}]\theta+O(t\theta^2)}d\theta, \tag{18}$$

where

$$f(0)=\frac{\mathrm{i}tan\big(t-\pi a/2-\varepsilon\sqrt{\pi t}/2\big)}{\left[\big(2-(1+\varepsilon/\sqrt{8})\sqrt{8\pi/t}\big)^2-8\pi/t\right]^{1/4}}. \tag{19}$$

This gives an estimate of the first integral in (13) of the form

$$\frac{\mathrm{i}tan\big(t-\pi a/2-\varepsilon\sqrt{\pi t}/2\big)e^{\mathrm{i}\psi(t)}}{\left[\big(2-(1+\varepsilon/\sqrt{8})\sqrt{8\pi/t}\big)^2-8\pi/t\right]^{1/4}\left[t/2-\sqrt{\pi t}\big(1+\varepsilon/\sqrt{8}\big)/\sqrt{2}\right]}+O\big(\varepsilon^{-1/4}t^{-5/4}\big), \tag{20}$$

and provided $\varepsilon$ is chosen so that $t(1+\pi R/t)/2=\big(t-\pi a/2-\varepsilon\sqrt{\pi t}/2\big)\neq(2m+1)\pi/2$, this will be of $O(t^{-1})$ and Lemma 2.3 is proved.

2.4 *Estimation of the second integral in* (13)

This integral is far less straightforward than its companion, because over the range $\theta\in(\pi,2\pi)$ the real phase (14) is positive and the integrand grows rapidly away from its endpoints. Consequently to establish an estimate, some alternative methodology is required. First making the substitution $w=\pi Re^{\mathrm{i}\theta}/t$ transforms the integral to

$$-\frac{\mathrm{i}t}{\pi R}\int_{-\pi R/t}^{\pi R/t}\frac{e^{\mathrm{i}t(log(pc(w))+1/pc(w))/2}}{[(1+w)^2-8\pi/t]^{1/4}}\frac{1-e^{-\mathrm{i}t-\mathrm{i}wt}}{1+e^{-\mathrm{i}t-\mathrm{i}wt}}dw. \tag{21}$$

Now over this integration range $Im(w) < 0$, which means that $\left|e^{-\mathrm{i}t-\mathrm{i}wt}\right| < 1$ and the denominator can be replaced by its geometric series giving

$$-\frac{\mathrm{i}t}{\pi R}\sum_{N=0}^{\infty}(-1)^N\int_{-\pi R/t}^{\pi R/t}\frac{e^{\mathrm{i}t\left(log\left(pc(w)\right)+1/pc(w)\right)/2}\left(1-e^{-\mathrm{i}t(1+w)}\right)e^{-\mathrm{i}t(1+w)N}}{[(1+w)^2-8\pi/t]^{1/4}}dw\,. \tag{22}$$

The path of $w$ is the circumference of a semi-circle, radius just less than unity, centred on the origin in the plane $Im(w) < 0$. Joining the points $-\pi R/t$ and $\pi R/t$ creates the boundary of a semi-circular disc, whose interior is free of any poles of the various integrands in (22). Consequently application of Cauchy's theorem means that the individual integrals in (22) are equivalent to the *same* integrals obtained by replacing the *complex* variable $w$ by the *real* variable $x$ and integrating along the real axis from $-\pi R/t$ to $\pi R/t$. Splitting this real integration range into two equal sections $(-\pi R/t\,,0)$ and $(0,\,\pi R/t)$, and applying the substitution $y = -x$ to the former transforms (22) to

$$-\frac{\mathrm{i}t}{\pi R}\sum_{N=0}^{\infty}(-1)^N\left\{\int_0^{\pi R/t}\frac{e^{\frac{\mathrm{i}t\left(log\left(pc(x)\right)+1/pc(x)\right)}{2}}\left(1-e^{-\mathrm{i}t(1+x)}\right)e^{-\mathrm{i}t(1+x)N}}{[(1+x)^2-8\pi/t]^{\frac{1}{4}}}dx+\int_0^{\pi R/t}\frac{e^{\frac{\mathrm{i}t\left(log\left(pc(y)\right)+1/pc(y)\right)}{2}}\left(1-e^{-\mathrm{i}t(1-y)}\right)e^{-\mathrm{i}t(1-y)N}}{[(1-y)^2-8\pi/t]^{\frac{1}{4}}}dy\right\}. \tag{23}$$

Making the further substitutions $u = 1/x$ and $u = 1/y$ respectively, and combining those integrals that appear twice in the sum, gives the second integral in (13) as an exact sum

$$-\frac{\mathrm{i}t}{\pi R}\left[C(0,R)+D(0,R)+2\sum_{N=1}^{\infty}(-1)^N\{C(N,R)+D(N,R)\}\right], \tag{24}$$

where

$$C(N,R)=\int_{\frac{t}{\pi R}}^{\infty}\frac{exp\left[-\mathrm{i}tN(1+1/u)+\mathrm{i}(t/2)\left(log\left(pc(u)\right)+1/pc(u)\right)\right]}{u^2[(1+1/u)^2-8\pi/t]^{1/4}}du, \tag{25}$$

$$D(N,R)=\int_{\frac{t}{\pi R}}^{\infty}\frac{exp\left[-\mathrm{i}tN(1-1/u)+\mathrm{i}(t/2)\left(log\left(pc(u)\right)+1/pc(u)\right)\right]}{u^2[(1-1/u)^2-8\pi/t]^{1/4}}du, \tag{26}$$

with $pc(u) \equiv pc(z = t(1+1/u)/\pi)$ in (25), and $pc(u) \equiv pc(z = t(1-1/u)/\pi)$ in (26). The integrals $C(N,R)$ and $D(N,R)$ are the subject of intensive analysis presented in Appendix B, the results of which can be summarised as follows.

For $R \approx_p R_{max}$, the sum of (25) & (26) together gives rise to three main contributions (the first three terms of B55) which are listed, together with their associated error/correction terms, by equations (B53-55). If $N$ lies within about $t^{1/6}$ of $\sqrt{t/2\pi}$, two of these main contributions become inaccurate and certain specified upper bound terms must be substituted instead. If $R = R_{max} - R_t$ these latter substitutions are unnecessary, and (B55) remains valid for all $N$ (subject to the restrictions listed in Section B1.7). The objective of the next sub-sections is to analyse the results of the sums of these three main contributions (and the substitute upper bounds, listed in parts b)-d) of *Lemma B1.2*, if necessary) which occur when (B55) is incorporated into (24), to ultimately produce an estimate for series (6) as a whole and prove the *Main Theorem*.

2.4a *Estimation of the second integral in* (13): The *sum of first main contribution of* (B55).

The first main contribution to $C(N,R) + D(N,R)$, which applies for all $N$ ≥0, is given by the first term in equation (B55). Substituting this term and its error into (24) yields the following sum

$$-\frac{t}{\pi R}\left(\frac{2e^{[\mathrm{i}t\{log(pc_+)+1/pc_+\}/2]}}{t[(1+\pi R/t)^2 - 8\pi/t]^{\frac{1}{4}}}\right)\left[\frac{2\sqrt{pc_+}}{a} + 2\sum_{N=1}^{\infty}(-1)^N\left\{\frac{e^{-\mathrm{i}tN(1+\pi R/t)}}{\left(a/2\sqrt{pc_+} - 2N\right)} + O\left(\frac{1}{t\left(a/2\sqrt{pc_+} - 2N\right)^2}\right)\right\}\right], \tag{27}$$

where $pc_+ = pc(z = t(1+\pi R/t)/\pi) \approx a^2/4$, so that $a/2\sqrt{pc_+} = 1 + O(t^{-1/2})$. The sum in (27) can be calculated using standard results involving Fourier series and the LerchPhi function. From [16] (1.445.7 & 1.445.8) one has

$$\sum_{N=1}^{\infty}\frac{(-1)^N N sin(Nx)}{(N^2-q^2)} = \frac{1}{2}\left[\sum_{N=1}^{\infty}\frac{(-1)^N sin(Nx)}{N-q} + \sum_{N=1}^{\infty}\frac{(-1)^N sin(Nx)}{N+q}\right] = \frac{\pi}{2}\frac{sin[q(2m\pi - x)]}{sin(q\pi)}, \tag{28}$$

$$\sum_{N=1}^{\infty}\frac{(-1)^N cos(Nx)}{(N^2-q^2)} = \frac{1}{2q}\left[\sum_{N=1}^{\infty}\frac{(-1)^N cos(Nx)}{N-q} - \sum_{N=1}^{\infty}\frac{(-1)^N cos(Nx)}{N+q}\right] = \frac{1}{2q^2} - \frac{\pi}{2}\frac{cos[q(2m\pi - x)]}{q sin(q\pi)}, \tag{29}$$

where $q > 0$ (but not an integer) and $(2m-1)\pi < x < (2m+1)\pi$. Combining $q \times (29) - \mathrm{i}\times(28)$ gives

$$\frac{1}{2}\left[\sum_{N=1}^{\infty}\frac{(-1)^N e^{-\mathrm{i}Nx}}{N-q} - \sum_{N=1}^{\infty}\frac{(-1)^N e^{\mathrm{i}Nx}}{N+q}\right] = \frac{1}{2q} - \frac{\pi}{2}\frac{e^{\mathrm{i}q(2m\pi - x)}}{sin(q\pi)}$$

$$\Rightarrow \sum_{N=1}^{\infty}\frac{(-1)^N e^{-\mathrm{i}Nx}}{N-q} = \frac{1}{q} - \frac{\pi e^{\mathrm{i}q(2m\pi-x)}}{sin(q\pi)} + \sum_{N=0}^{\infty}\frac{(-1)^N e^{\mathrm{i}Nx}}{N+q} - \frac{1}{q} = -\frac{\pi e^{\mathrm{i}q(2m\pi-x)}}{sin(q\pi)} + \Phi(-e^{\mathrm{i}x}, 1, q), \tag{30}$$

where $\Phi(-e^{\mathrm{i}x}, 1, q)$ is the LerchPhi function ([16], 9.55) (NB. not to be confused with the confluent hypergeometric function $\Phi(a,c;z)$.) Hence (27) becomes

$$-\frac{t}{\pi R}\left(\frac{2e^{[\mathrm{i}t\{log(pc_+)+1/pc_+\}/2]}}{t[(1+\pi R/t)^2 - 8\pi/t]^{\frac{1}{4}}}\right)\left[\frac{2\sqrt{pc_+}}{a} + \frac{\pi e^{\mathrm{i}a(2m\pi - t(1+\pi R/t))/4\sqrt{pc_+}}}{sin\left(\pi a/4\sqrt{pc_+}\right)} - \Phi\left(-e^{\mathrm{i}t(1+\pi R/t)}, 1, \frac{a}{4\sqrt{pc_+}}\right)\right]. \tag{31}$$

Provided $e^{it(1+\pi R/t)} \neq -1 \Rightarrow t(1+\pi R/t) \neq (2m+1)\pi$, (already excluded to avoid the poles of integrand 8, cf. equation 20) the LerchPhi function can be represented as an integral ([16], 9.556)

$$\Phi\left(-e^{it(1+\pi R/t)}, 1, \frac{a}{4\sqrt{pc_+}}\right) = \int_0^\infty \frac{e^{-ay/4\sqrt{pc_+}}}{1+e^{it(1+\pi R/t)}e^{-y}} dy \approx e^{-\mathrm{i}(t(1+\pi R/t)/2+\pi/2)} log\left[\mathrm{i}cot\left(\frac{t(1+\pi R/t)}{4}+\frac{\pi}{4}\right)\right]. \tag{32}$$

The last approximation is valid because $a/4\sqrt{pc_+} \approx 1/2$. Hence (31) reduces to

$$\approx -\frac{t}{\pi R}\left(\frac{2e^{[\mathrm{i}t\{log(pc_+)+1/pc_+\}/2]}}{t[(1+\pi R/t)^2-8\pi/t]^{\frac{1}{4}}}\right)\left[1+\pi e^{\mathrm{i}(2m\pi - t(1+\pi R/t))/2} - e^{-\mathrm{i}\left(\frac{t(1+\pi R/t)}{2}+\pi/2\right)} log\left[\mathrm{i}cot\left(\frac{t(1+\pi R/t)}{4}+\frac{\pi}{4}\right)\right]\right]. \tag{33}$$

The upshot of all this is that for both $R \approx_p R_{max}$ and $R = R_{max} - R_t$, the sum of the first main contribution of (B55) substituted into (24) yields a term of $O(t^{-1})$, given by (33). The sum of the associated error terms in (27) cannot amount to anything to alter this conclusion.

2.4b *Estimation of the second integral in* (13): *The sum of second main contribution of* (B55).

The second main contribution to $C(N,R)+D(N,R)$, which applies for $N \notin (N_{2\eta}, N_t]$ (see *Lemma B1.2* for definitions of $N_{2\eta}$ and $N_\mu$, with $\mu = 5/\sqrt{8\pi}$) when $R \approx_p R_{max}$, is given by the second term in equation (B55). Substituting this term into (24) gives the following sums

$$\frac{t}{\pi R}\left(\frac{2e^{[\mathrm{i}t\{log(pc_-)+1/pc_-\}/2]}}{t[(1-\pi R/t)^2-8\pi/t]^{\frac{1}{4}}}\right)\left[\frac{2\sqrt{pc_-}}{a} - 2\left\{\sum_{N=1}^{N_1}\frac{(-1)^N e^{-\mathrm{i}tN(1-\pi R/t)}}{2N-a/2\sqrt{pc_-}} + \sum_{N=N_2}^{\infty}\frac{(-1)^N e^{-\mathrm{i}tN(1-\pi R/t)}}{2N-a/2\sqrt{pc_-}}\right\}\right], \tag{34}$$

where $N_1 = N_{2\eta}$ and $N_2 = \left\lceil\sqrt{t/2\pi}\right\rceil$. If $N_2 > q = a/4\sqrt{pc_-}$, then the series

$$\sum_{N=N_2}^{\infty}\frac{(-1)^N e^{-\mathrm{i}Nx}}{N-q} = e^{-\mathrm{i}N_2 x}(-1)^{N_2}\sum_{p=0}^{\infty}\frac{(-1)^p e^{-\mathrm{i}px}}{(p-(q-N_2))}$$

$$= e^{-\mathrm{i}N_2 x}(-1)^{N_2}\left[\Phi\left(-e^{-\mathrm{i}x}, 1, (N_2-q)\right)\right] = e^{-\mathrm{i}N_2 x}(-1)^{N_2}\int_0^\infty \frac{e^{-(N_2-q)y}}{1+e^{-\mathrm{i}x}e^{-y}} dy$$

$$\approx e^{-\mathrm{i}N_2 x}(-1)^{N_2}\left[\frac{1}{(N_2-q)(1+e^{-\mathrm{i}x})} + O\left(\frac{1}{[(N_2-q)(1+e^{-\mathrm{i}x})]^2}\right)\right], \tag{35}$$

provided $e^{\mathrm{i}x} \neq -1$. Similarly if $N_1 < q = a/4\sqrt{pc_-}$,

$$\sum_{N=1}^{N_1}\frac{(-1)^N e^{-iNx}}{N-q}=\sum_{N=-\infty}^{N_1}\frac{(-1)^N e^{-iNx}}{(N-q)}-\sum_{N=-\infty}^{0}\frac{(-1)^N e^{-iNx}}{(N-q)}$$
$$=-e^{-iN_1x}(-1)^{N_1}\sum_{p=0}^{\infty}\frac{(-1)^{-p}e^{ipx}}{\big(p+(q-N_1)\big)}+\sum_{p=0}^{\infty}\frac{(-1)^{-p}e^{ipx}}{(p+q)}$$
$$=-e^{-iN_1x}(-1)^{N_1}\{\Phi\big(-e^{ix},1,(q-N_1)\big)\}+\Phi(-e^{ix},1,q)$$
$$\approx\frac{1}{(1+e^{ix})}\left\{-\frac{e^{-iN_1x}(-1)^{N_1}}{(q-N_1)}+\frac{1}{q}\right\}+O\left(\frac{1}{[(q-N_1)(1+e^{-ix})]^2}\right). \qquad (36)$$

Now if $R\approx R_{max}$, then $a/4\sqrt{pc_-}\approx\sqrt{t/2\pi}-\eta t^{1/6}$, so the differences $|N_1-a/4\sqrt{pc_-}|$ and $|N_2-a/4\sqrt{pc_-}|$ are both of $O(\eta t^{1/6})$. Hence provided $e^{ix}\equiv e^{it(1+\pi R/t)}\neq -1$, the error terms in (35-36) will be $O(t^{-1/6})$ smaller than the main terms. Consequently (34) can be approximated by the first order terms of (35-36), giving

$$\frac{t}{\pi R}\left(\frac{2e^{[it\{log(pc_-)+1/pc_-\}/2]}}{t[(1-\pi R/t)^2-8\pi/t]^{\frac{1}{4}}}\right)\left[\frac{2\sqrt{pc_-}}{a}\right.$$
$$\left.-\frac{1}{(1+e^{it(1-\pi R/t)})}\left\{-\frac{e^{-iN_1t(1-\pi R/t)}(-1)^{N_1}}{\big(a/4\sqrt{pc_-}-N_1\big)}+\frac{4\sqrt{pc_-}}{a}+\frac{e^{-iN_2t(1-\pi R/t)}(-1)^{N_2}}{\big(N_2-a/4\sqrt{pc_-}\big)}\right\}+O\big(t^{-1/3}\big)\right]$$
$$(37)$$

If $R\approx_p R_{max}$ then $[(1-\pi R/t)^2-8\pi/t]^{\frac{1}{4}}\approx 2\sqrt{\pi\eta}t^{-5/12}$, which means that (37) is a term of $O(t^{-3/4})$, considerably larger than the sum of the first main contribution (33). If $R=R_{max}-R_t$ then $N_1$ and $N_2$ can be allowed to converge, as the relative error term (B32) is now a negligible $O(t^{-1/2})$ for all $N$. So in this case (34) becomes

$$\frac{t}{\pi R}\left(\frac{2e^{[it\{log(pc_-)+1/pc_-\}/2]}}{t[(1-\pi R/t)^2-8\pi/t]^{\frac{1}{4}}}\right)\left[\frac{2\sqrt{pc_-}}{a}-2\left\{\sum_{N=1}^{\infty}\frac{(-1)^N e^{-itN(1-\pi R/t)}}{2N-a/2\sqrt{pc_-}}\right\}\right], \qquad (38)$$

provided $a/4\sqrt{pc_-}\neq$ integer (in which case the term $N=a/4\sqrt{pc_-}$ vanishes and the corresponding $D(N,R)$ integral contains only a saddle point term of the form B18). Applying (30) to (38) one obtains

$$\frac{t}{\pi R}\left(\frac{2e^{[it\{log(pc_-)+1/pc_-\}/2]}}{t[(1-\pi R/t)^2-8\pi/t]^{\frac{1}{4}}}\right)\left[\frac{2\sqrt{pc_-}}{a}+\left\{\frac{\pi e^{ia(2m\pi-t(1-\pi R/t))/4\sqrt{pc_-}}}{sin\big(a\pi/4\sqrt{pc_-}\big)}-\Phi\left(-e^{it(1-\pi R/t)},1,\frac{a}{4\sqrt{pc_-}}\right)\right\}\right].$$
$$(39)$$

When $R=R_{max}-R_t$ the term $[(1-\pi R/t)^2-8\pi/t]^{\frac{1}{4}}\sim O(t^{-1/4})$, which means that (39), like (37), is also a term of $O(t^{-3/4})$. Since the choices of $\eta$ and $\mu$ set out in the proof of *Lemma B1.2* guarantees that the magnitude of the associated error term never exceeds 2% of the second main contribution, the sum of such error terms incorporated into (24) cannot amount to anything greater than (37) or (39).

2.4c *Estimation of the second integral in* (13): *Sum of third main contribution and associated error of* (B55).

First define $N_t^- = min\{N_\mu, \lfloor a/4\sqrt{pc_-}\rfloor\}$ for both $R \approx_p R_{max}$ and $R = R_{max} - R_t$. Then substituting the third main contribution to $C(N,R) + D(N,R)$ and its associated error into (24) gives

$$\frac{2\text{i}t}{\pi R}\left\{\frac{\sqrt{2\pi}e^{-\text{i}\left(t-\frac{t}{2}log\left(\frac{t}{2\pi}\right)+\pi+\frac{2\pi^2}{t}\dots\right)-\text{i}\pi/4}}{\sqrt{t}(1-8\pi/t)^{\frac{1}{4}}}+\frac{1.1691}{(1-8\pi/t)^{\frac{1}{4}}t^{\frac{2}{3}}}\right\}$$

$$-\frac{2\text{i}t}{\pi R}\sum_{N=2}^{N_t^-}(-1)^N\left\{\sqrt{\frac{2\pi}{Nt}}e^{-\text{i}\pi/4}exp\left(\text{i}\frac{t}{2}log\left(\frac{t}{2\pi N^2}\right)-\text{i}t-\text{i}\pi N^2\right)+O_1\left(\frac{2^{3/2}5\sqrt{\pi}(1/N+2\pi N/t)^2}{16(tN)^{3/2}(1/N-2\pi N/t)^3}\right)\right\}. \quad (40)$$

The question here is not so much the value of the main sum in (40) as the potential size of the sum of the $O_1$ error terms. One can simplify things a little by noting that

$$\left|\frac{(-1)^N 2^{3/2}5\sqrt{\pi}(1/N+2\pi N/t)^2}{16(tN)^{\frac{3}{2}}(1/N-2\pi N/t)^3}\right| < \frac{1.566\times 4}{t^{\frac{3}{2}}N^{\frac{7}{2}}(1/N-2\pi N/t)^3}, \quad (41)$$

for all $N$ in (40). Now define $M = \lfloor -log\left(1 - N_t^-\sqrt{2\pi}/\sqrt{t}\right)/log(2)\rfloor$ and divide up the range of $N$ into $(M+1)$ subsections given by $\left[\frac{2}{1}, \frac{\sqrt{t}}{2\sqrt{2\pi}}\right], \dots, \left[\frac{(2^{p-1}-1)\sqrt{t}}{2^{p-1}\sqrt{2\pi}}, \frac{(2^p-1)\sqrt{t}}{2^p\sqrt{2\pi}}\right], \dots, \left[\frac{(2^M-1)\sqrt{t}}{2^M\sqrt{2\pi}}, N_t^-\right]$ respectively. Now in the first subsection, (41) has a maximum value of $\approx 1.566 \times 2\sqrt{2}t^{-3/2}$, when $N = 2$. (In fact equation 41 is minimized when $N \approx \sqrt{t}/\sqrt{26\pi}$.) Within the other subsections, (41) is maximised at the upper endpoint. At this endpoint, $N = (2^p - 1)\sqrt{t}/2^p\sqrt{2\pi}$, and one has

$$\frac{1.566\times 4}{t^{\frac{3}{2}}N^{\frac{7}{2}}(1/N-2\pi N/t)^3} = \frac{1.566(2\pi)^{\frac{1}{4}}2^{3p-1}}{t^{\frac{7}{4}}(1-2^{-p})^{\frac{1}{2}}(1-2^{-p-1})^3}, \quad (42)$$

whilst in each subsection (apart from the last) there are approximately $\sqrt{t}/2^p\sqrt{2\pi}$ integers. If the $O(t^{-1})$ term from sum over the first subsection is excluded, the next $M-1$ subsections give rise to

$$\sum_{N=\sqrt{t}/2\sqrt{2\pi}}^{(2^M-1)\sqrt{t}/2^M\sqrt{2\pi}} O_1\left(\frac{(-1)^N 2^{3/2}5\sqrt{\pi}(1/N+2\pi N/t)^2}{16(tN)^{\frac{3}{2}}(1/N-2\pi N/t)^3}\right) < \frac{1.566}{(2\pi)^{\frac{1}{4}}t^{\frac{5}{4}}}\sum_{p=2}^{M}\frac{2^{2p-1}}{(1-2^{-p})^{\frac{1}{2}}(1-2^{-p-1})^3}$$

$$< \frac{1.566}{(2\pi)^{\frac{1}{4}}t^{\frac{5}{4}}}\times 0.87\times 4^M \approx \frac{0.86(4^M)}{t^{5/4}}, \quad \text{(for all } M \geq 3). \quad (43)$$

Now if $pc_{-} = 1 + \eta\sqrt{8\pi}t^{-1/3}$ as is the case for $R \approx_p R_{max}$, then $M < log(t)/3log(2)$. Hence

$$log\left[4^M/t^{5/4}\right] < log(t)log(4)/3log(2) - 5log(t)/4 = -7log(t)/12$$

$$\Rightarrow \frac{0.86(4^M)}{t^{5/4}} < \frac{0.86}{t^{7/12}}. \quad (44)$$

Of course if $R = R_{max} - R_t$ then $pc_{-} = 1 + O(1)$ and (43) is only $O\left(t^{-5/4}\right)$. For the last subsection the maximum value of (41) occurs when $N = N_t^{-} = N_\mu$, where

$$\frac{1.566 \times 4}{t^{\frac{3}{2}}N^{\frac{7}{2}}(1/N - 2\pi N/t)^3} \le \frac{1.566 \times 16(\pi/2^7)^{\frac{1}{4}}}{125t^{\frac{3}{4}}} + O\left(t^{-13/12}\right). \quad (45)$$

The number of integers in the last subsection equates to

$$N_t^{-} - (2^M - 1)\sqrt{t}/2^M\sqrt{2\pi} = \sqrt{t}/2^M\sqrt{2\pi} - 5t^{1/6}/\sqrt{8\pi} \le 5t^{1/6}/\sqrt{8\pi}, \quad (46)$$

giving rise to an upper bound of the sum over this section of the form

$$\frac{5t^{1/6}}{\sqrt{8\pi}} \times \frac{1.566 \times 16(\pi/2^7)^{\frac{1}{4}}}{125t^{3/4}} < \frac{0.08}{t^{7/12}}. \quad (47)$$

Combining (44) and (47) the sum over the error terms in (40) cannot exceed a value of $0.94t^{-7/12}$, although in practice it is likely to be considerably smaller than this.

2.4d *Estimation of the second integral in* (13): *Sums of other possible contributions potentially arising from* (B55).

If $R \approx_p R_{max}$ then both the second and third main terms change character somewhat as one approaches close to $N \approx N_t$. For $N \in (N_\mu, N_\eta]$ the third term derived from the integral through the saddle is bounded above by (B31). The maximum magnitude of (B31) can be taken to be the value at $r_N = \mu \approx 1$, which is about $0.64 \times 1.728\sqrt{2\pi/Nt}$ (as discussed after B31, the factor 0.64 is a good numerical estimate of the reduction brought about by the neglected exponential decay across the saddle) dying away as $r_N^{3/8}$ as $r_N \longrightarrow 0$. Hence an upper bound on the contribution of a term of this size summed across this range of $N$ is

$$\left(N_\eta - N_\mu\right) \times 0.64 \times 1.728\sqrt{\frac{2\pi}{Nt}} \approx (\mu - \eta)t^{1/6} \times 1.106\left(\frac{2\pi}{t}\right)^{3/4} \approx \frac{0.313}{t^{7/12}}. \quad (48)$$

For $N \in (N_{2\eta}, N_t]$, the second main term should be replaced by the upper bound given by (B37) for various choice of integration paths specified in the discussion surrounding (B35). The maximum modulus for (B37) is $3.55t^{-3/4}$ but this is reduced by at least a factor of 0.38 due to the neglected exponential decay. So an upper bound on the contribution of this term summed across this range of $N$ is

$$2\eta t^{\frac{1}{6}} \times 0.38 \times 3.55t^{-3/4} \approx 2.50t^{-7/12}. \tag{49}$$

In practice both (48) and (49) again represent considerable overestimates. If $R = R_{max} - R_t$, then neither (48) nor (49) are relevant, and only results (33), (39) and (44 & 47) apply.

2.5 *Estimation of the second integral in* (13): *Summary*

In summary the results of sections 2.4a-d mean that the second integral in (13) can, via (24) and (B53-55), be approximated for any $R \approx_p R_{max}$ by the following:

$$\frac{2\mathrm{i}t}{\pi R}\left\{\frac{\sqrt{2\pi}e^{-\mathrm{i}\left(t-\frac{t}{2}log\left(\frac{t}{2\pi}\right)+\pi+\frac{2\pi^2}{t}\right)-\mathrm{i}\pi/4}}{\sqrt{t}(1-8\pi/t)^{\frac{1}{4}}}\right.$$
$$\left.-\sum_{N=2}^{N_t^-}(-1)^N\left\{\sqrt{\frac{2\pi}{Nt}}e^{-\mathrm{i}\pi/4}exp\left(\mathrm{i}\frac{t}{2}log\left(\frac{t}{2\pi N^2}\right)-\mathrm{i}t-\mathrm{i}\pi N^2\right)\right\}\right\}+\epsilon,$$

(50)

where $N_t^- = min\{N_\mu, \lfloor a/4\sqrt{pc_-}\rfloor\}$ and the error term satisfies

$$|\epsilon| < \frac{1.1691}{t^{\frac{2}{3}}} + O\left(\frac{1}{t}\right)_{(33)} + O\left(\frac{1}{t^{3/4}}\right)_{\substack{(37)\text{ or}\\(39)}} + \left(\frac{0.94}{t^{7/12}}\right)_{\substack{(44)+\\(47)}} + \left(\frac{0.313}{t^{7/12}}\right)_{(48)} + \left(\frac{2.50}{t^{7/12}}\right)_{(49)}.$$

(51)

In (51) the subscripts indicate the equation number of the specific term(s) from which this order of magnitude estimate is derived. For $R = R_{max} - R_t$ rather than $R \approx_p R_{max}$, the last two terms in (51) do not apply and the third term is derived from (39) not (37).

2.6 *Proof of Main Theorem*

For either $R \approx_p R_{max}$ or $R = R_{max} - R_t$ estimate (20) of the first integral in (13) is a small term of $O(t^{-1})$ which has to be added (but will not significantly alter) the error term $\epsilon$ given by (51). Substituting the main terms given by (50), augmented by this adapted error term (51), into equation (13) gives

$$-\frac{R\sqrt{\pi}e^{\mathrm{i}t/2+\mathrm{i}\pi/8}}{\mathrm{i}\sqrt{2t}}\left[\frac{2\mathrm{i}t}{\pi R}\left\{\frac{\sqrt{2\pi}e^{-\mathrm{i}\left(t-\frac{t}{2}log\left(\frac{t}{2\pi}\right)+\pi+\frac{2\pi^2}{t}\dots\right)-\mathrm{i}\pi/4}}{\sqrt{t}(1-8\pi/t)^{\frac{1}{4}}}\right.\right.$$
$$\left.\left.-\sum_{N=2}^{N_t^-}(-1)^N\left\{\sqrt{\frac{2\pi}{Nt}}e^{-\mathrm{i}\pi/4}exp\left(\mathrm{i}\frac{t}{2}log\left(\frac{t}{2\pi N^2}\right)-\mathrm{i}t-\mathrm{i}\pi N^2\right)\right\}+\epsilon\right\}\right]$$

$$=2e^{\mathrm{i}\frac{t}{2}log\left(\frac{t}{2\pi}\right)-\mathrm{i}t/2-\mathrm{i}\pi/8}\left[\frac{e^{-\mathrm{i}2\pi^2/t\dots}}{(1-8\pi/t)^{\frac{1}{4}}}+\sum_{N=2}^{N_t^-}\frac{exp(-\mathrm{i}tlog(N))}{\sqrt{N}}\right]-\sqrt{\frac{2t}{\pi}}e^{\mathrm{i}t/2+\mathrm{i}\pi/8}\epsilon. \qquad (52)$$

Hence taking the real parts of (52) one deduces from (6 & 8) that

$$\sum_{\alpha=INT_O(a)+2}^{INT_O(2t/\pi-a)}\frac{2\sqrt{2}\cos\left((t/2)\{log(pc)+1/pc\}+t/2+\pi/8\right)}{(\alpha^2-a^2)^{\frac{1}{4}}}$$
$$=2\left[\frac{cos\left\{\frac{t}{2}log\left(\frac{t}{2\pi}\right)-\frac{t}{2}-\frac{\pi}{8}-\frac{2\pi^2}{t}\dots\right\}}{(1-8\pi/t)^{\frac{1}{4}}}+\sum_{N=2}^{N_t^-}\frac{cos\left\{\frac{t}{2}log\left(\frac{t}{2\pi}\right)-\frac{t}{2}-\frac{\pi}{8}-tlog(N)\right\}}{\sqrt{N}}\right]-\sqrt{\frac{2t}{\pi}}\epsilon. \qquad (53)$$

From (53) one clearly sees the connection between the main sum (6) of the new series (A64) with the main sum of the RS formula (7) used to calculate $Z(t)$ along the critical line. Note the distinctive first term and the fact that the cosine phases are $\theta(t)-tlog(N)$ minus the $O(t^{-1})=1/48t+\cdots$ of the full $\theta(t)$ expansion. Incorporating these $O(t^{-1})$ terms, so that the phase on the right of (53) is exactly $\theta(t)$ as stated in the main theorem (7), corresponds to adding an $O(t^{-3/4})$ correction to the error term appearing in (53). This is insignificant because the currently the maximum size of the error term (to first order in $t$) is strictly less than

$$|error|=\sqrt{\frac{2t}{\pi}}|\epsilon|<\sqrt{\frac{2t}{\pi}}\left\{\left(\frac{0.94}{t^{7/12}}\right)+\left(\frac{0.313}{t^{7/12}}\right)+\left(\frac{2.50}{t^{7/12}}\right)\right\}\approx\frac{2.99}{t^{1/12}}, \qquad (54)$$

utilising all the results summarised by equation (51). However, the calculations required to establish this estimate limited the radius of the circle $C(R)$ in (8) to be no greater than $R\approx_p R_{max}$, as defined by (9). If $R\approx_p R_{max}$ then

$$N_t^-=N_\mu=\left\lfloor\sqrt{\frac{t}{2\pi}}-\frac{5t^{1/6}}{\sqrt{8\pi}}\right\rfloor. \qquad (55)$$

Consequently the last $5t^{1/6}/\sqrt{8\pi}$ terms up to $N_t$ of the main RS sum are not included in (53), and replacing $N_t^-$ by $N_t$ (to fully complete the link to the RS sum) has the potential to introduce an extra source of error into (54). One would expect that the sum of these last $(N_t-N_t^-)$ terms of the RS sum will reach a significant in size in the region $|a-INT_O(a)|<$ $t^{-1/6}$ , because that is when the transition term $T_{a+\varepsilon}$ (A65), which was excluded from the analysis leading to (53), is also relatively large. Therefore one would expect

$$T_{a+\varepsilon} \approx \sum_{N=N_t^-+1}^{N_t} \frac{2cos\left\{\frac{t}{2}log\left(\frac{t}{2\pi}\right)-\frac{t}{2}-\frac{\pi}{8}-tlog(N)\right\}}{\sqrt{N}}, \qquad (56)$$

and to become relatively large only when $a$ is within $\pm t^{-1/6}$ of an odd integer. Numerical evidence strongly supports this idea. When $a$ is exactly equally to an odd integer (so that $\varepsilon = g = 0$ in equation A65) the magnitude of the sum on the RHS of (56) is consistently between 13-14% smaller than the fixed value of $T_a = 1.4089t^{-1/12}$. (N.B. when $a$ is an odd integer, the phase in A65 $\cos(t+\pi/24) = \pm\sqrt{3}/2$. The $\pm$ sign alternates between successive pairs of odd integers.) From this point, the magnitude of the RHS of (56) increases slowly, until it reaches a maximum value of $\approx 2.30t^{-1/12}$ whenever $a + 0.93t^{-1/6} = odd\ integer$, (this magnitude $2.30t^{-1/12}$ is consistently 7% or so larger than the magnitude of the first standard term in equation 6, which is the limit for $T_{a+\varepsilon}$ as $\varepsilon$ increases beyond $t^{-1/6}$). Subsequently the sum (56) rapidly falls away to something $o(t^{-1/12})$ as $a$ moves outside the transition zone. This cycle then repeats itself as $a$ moves to within close proximity of the next odd integer, continuing in this way for odd integers up to at least $t = 10^{30}$ and in all probability indefinitely. But more detailed analysis would be necessary to establish a proof of these numerical correlations for all $t$. What is surprising is the relatively large size of this numerical maximum of $2.30t^{-1/12}$ when compared to the absolute maximum of no more than $3.16t^{-1/12}$, which arises from (56) when all the cosine terms are unity. Clearly for such a relatively large maximum to come about the phases of the cosines must correlate, modulo $2\pi$, to a remarkable degree when $a + 0.93t^{-1/6} = odd\ integer$. In the absence of a proof of these numerical results for arbitrarily large $t$, extending the RS sum in (53) from $N_t^-$ to $N_t$ has the potential to increase the upper bound on the $|error|$ in (54) to

$$\mathrm{E}(t) = |error| < \frac{2.99}{t^{1/12}} + \frac{3.16}{t^{1/12}} = \frac{6.15}{t^{1/12}} \xrightarrow[t\to\infty]{} 0, \qquad (57)$$

although the introduction of the transition term into (53) when $a$ is within $\pm t^{-1/6}$ of an odd integer will, in practice, ensure that this last addition is a considerable overestimate. This completes the proof of the main theorem (7). As a corollary, proof of the main theorem demonstrates that series (6, A64) does indeed yield an asymptotic approximation for $Z(t)$, *independent* of the interpretation (A21) made for the RSI, which led to its discovery.

## 3. Increasing the computational efficiency of the RS Formula

3.1 *Introduction*

Improvements to the computational efficiency of calculating $Z(t)$ over and above what can be achieved using the RS formula, can be discerned by considering how the main terms in that formula come to be linked to the main terms in the new series representation (6, 53). The connection between the odd integer regime $\alpha$ of the new series and the integers $N$ of the RS main sum can best be seen by viewing the variable $pc$ as sort of bridge. The RS terms on the right-hand side of (53) arise from the evaluation of the integral $D(N,R)$, when integrating across the saddle points of its exponential phase, which are situated at $N = a/4\sqrt{pc_{sad}}$ (see B12). From (5, A47) one can see that in "$\alpha$ space" these saddles correspond to points given by

$$\alpha^2 = \frac{2t(pc_{sad}+1)^2}{\pi pc_{sad}} = \frac{2t(t/2\pi N^2+1)^2}{\pi\, t/2\pi N^2} = 4N^2\left(\frac{t}{2\pi N^2}+1\right)^2 = \frac{t^2}{N^2\pi^2} + \frac{4t}{\pi} + 4N^2. \qquad (58)$$

Now clearly when $N$ is $o(t^{1/3})$, (58) $\Rightarrow \alpha \approx t/N\pi$ and the situation is as illustrated in Fig. 1a. This shows a circle of radius $R$, centred on $t/\pi$ (as specified in equation 8), enclosing the odd integers $\alpha$ (shown as black diamonds on the figure). Schematically these represent the terms of the new series (6). The integers $N = 1, 2, \ldots$, representing the terms of the RS main sum, are shown as black dots on the figure. If $R < t/2\pi$ then all the "diamonds" within the circle simply sum up to give the first term of the RS main sum *alone*. If $t/2\pi < R < t/3\pi$ then all the "diamonds" within the circle combine to give the sum of the *first and second* terms of the RS main sum and this pattern continues *for a while*. (This comes about because the gradient of the phase in (6) approximates to $2\pi N + O(N^3t^{-1})$, so that the terms combine together when $\alpha \approx t/N\pi$ for $N < t^{1/3}$.) Two observations concerning this process are worth making. First, when summing the new series for increasing $\alpha$ starting from $INT_O(a)+2$, one is actually computing the sum of the terms of the RS main sum *backwards*, ending with the $N = 1$ term. Second, it is now no surprise that when computing $Z(t)$ from the new zeta-sum using Euler-Maclaurin summation (see section A4), imposition of a minimum upper limit to the main sum of $INT_O(t/\pi)$ was necessary for the Bernoulli sum (A72) to converge. This minimum is necessary condition because attempting to try and analytically continue the series at any smaller limit must fall foul of the $N = 1$ term of the RS main sum, which will not have been accounted for.

However, the simple correlation of $\alpha = t/N\pi$ with the $N$th term of the RS main sum cannot continue indefinitely. To see this, consider how many terms of the RS main sum one might expect to find in the interval between $\alpha = a$ (when $pc = 1$ and $N = \sqrt{t/2\pi}$ in equation 58) and the first odd integer above $a$, that is $\alpha = INT_O(a) + 2 = a + \delta$. Now

$$pc(INT_O(a)+2 = a+\delta) = 1 + \left(\frac{8\pi}{t}\right)^{1/4}\delta^{1/2} + O\left(\frac{\delta}{\sqrt{t}}\right), \qquad (59)$$

so in $N$ space the interval $\delta$ is representative of an interval

$$\Delta N = \sqrt{\frac{t}{2\pi}} - \sqrt{\frac{t}{2\pi pc(a+\delta)}} \approx \sqrt{\frac{t}{2\pi}}\left\{1 - \frac{1}{\left(1+(8\pi/t)^{1/4}\sqrt{\delta}\right)^{1/2}}\right\} \approx \left(\frac{t}{2\pi}\right)^{1/4}\sqrt{\frac{\delta}{2}}. \qquad (60)$$

In other words between $\alpha = a$ and $\alpha = a + \delta$ there are potentially as many as $\Delta N \approx (t/2\pi)^{1/4}$ of the final terms (if $\delta$ takes its maximum value of two) of the RS main sum squeezed in the interval $\delta$. Of course for the most part $\delta$ will be much less than two, but the correlation established through equation (58) continues irrespective of the value of $\delta$. Consider the next interval, situated between the first two odd integers above $a$, that is $\alpha = INT_O(a) + 2 = a + \delta$ and $\alpha = INT_O(a) + 4 = a + \delta + 2$. This is representative of

$$\Delta N = \sqrt{\frac{t}{2\pi}}\left(\frac{1}{\sqrt{pc(a+\delta)}} - \frac{1}{\sqrt{pc(a+\delta+2)}}\right) \approx \left(\frac{t}{2\pi}\right)^{1/4}\left(\frac{\sqrt{\delta+2}-\sqrt{\delta}}{\sqrt{2}}\right) \qquad (61)$$

terms near the end of the RS main sum. Hence the *minimum* number of terms between these two odd integers is approximately $(\sqrt{2}-1)(t/2\pi)^{1/4}$ (when $\delta = 2$).

The situation is illustrated schematically in Fig. 1b, which shows that many black dots, representative of terms in the RS main sum, are now distributed between the intervals of odd integers (black diamonds) situated just above $\alpha = a$. (Note this is the inverse of the distributional equivalence near $\alpha = INT_O(t/\pi)$, shown in Fig. 1a.) Also shown in Fig. 1b are small sections of two circumferences of circles centred at $\alpha = t/\pi$ crossing the $\alpha$-axis just above and below $\alpha = INT_O(a) + 2$. Increasing $R$ so that the larger circle $C(R)$ now includes the pole of (8) situated at $\alpha = INT_O(a) + 2$, means that the corresponding *single* term in series in (6) is equivalent to *sum of the* $\Delta N \sim t^{1/4}$ *terms* (61) *of the Riemann-Siegel series*. More specifically, let $\alpha_1 = INT_O(a) + 2$, $\alpha_2 = \alpha_1 + 2$, and denote the corresponding nearest integer values satisfying (58) by $N_1$ and $N_2$ respectively. Then

$$\frac{2\sqrt{2}\cos\left(\frac{t}{2}\left\{log\big(pc(\alpha_1)\big) + \frac{1}{pc(\alpha_1)}\right\} + \frac{t}{2} + \frac{\pi}{8}\right)}{(\alpha_1^2 - a^2)^{\frac{1}{4}}} \approx \sum_{N=(N_1+N_2)/2}^{N_t} \frac{2cos\{\theta(t) - tlog(N)\}}{\sqrt{N}}. \qquad (62)$$

For example, if $t = 10^{24}$, then $\alpha_1 = 1595769121607$, $(N_1 + N_2)/2 = 398941625041$ and $N_t = 398942280401$. The sum of the $655{,}360$ terms on the right of (62) amounts to $4.727..\times 10^{-4}$, compared to the value $4.744..\times 10^{-4}$ on the left. The *relative* error of $0.0034 \equiv 2.8\Lambda_\alpha^{-1}$, in the range predicted by (A64). Consequently by substituting the early terms of new series (6) situated near $\alpha = a$, one effectively computes the sum of $\Delta N \sim t^{1/4}$ of the latter terms of the Riemann-Siegel series *in one go*, with obvious increases the computational efficiency. Hence it is potentially faster to calculate values for the main sum of terms making up $Z(t)$ by means of a hybrid formula of the form

$$\begin{aligned}\sum_{N=1}^{N_t} & \frac{2cos\left\{\frac{t}{2}log\left(\frac{t}{2\pi}\right) - \frac{t}{2} - \frac{\pi}{8} - tlog(N)\right\}}{\sqrt{N}} \\ &= \sum_{N=1}^{N_{CO}} \frac{2cos\left\{\frac{t}{2}log\left(\frac{t}{2\pi}\right) - \frac{t}{2} - \frac{\pi}{8} - tlog(N)\right\}}{\sqrt{N}} \\ &+ \sum_{\alpha=INT_O(a)+2}^{L_{CO}} \frac{2\sqrt{2}\cos\big((t/2)\{log(pc) + 1/pc\} + t/2 + \pi/8\big)}{(\alpha^2 - a^2)^{\frac{1}{4}}} + error, \qquad (63)\end{aligned}$$

where $N_{CO}$ and $L_{CO}$ are cut-off integers chosen to reduce the total number of terms in the two sums in (63) to a minimum. Of course in practice the utility of such an observation depends upon two points; the computational speed of calculating terms in the new series and the scale of the resulting error term that will accrue. These points are discussed below.

3.2 *Maximising the Computational efficiency of the hybrid formula.*

Let $S$ be the total number of terms in the proposed hybrid formula (63), given by

$$S = N_{CO} + \frac{(L_{CO} - INT_O(a) - 2)}{2}, \tag{64}$$

for a fixed value of $t$. Now if $S$ was a continuous function of the variable $pc(\alpha)$, then there would be a cut-off value $pc = pc_{CO}$ which can be used to define both

$$L_{CO} = \sqrt{\frac{2t}{\pi}\frac{(pc_{CO}+1)}{\sqrt{pc_{CO}}}} \qquad \text{and} \qquad N_{CO} = \sqrt{\frac{t}{2\pi pc_{CO}}} \tag{65}$$

(from equations A47 and B12 respectively) as continuous variables. Now to maximise the computational efficiency of the hybrid for calculations of $Z(t)$, one needs to minimise the total CPU time needed to calculate (63). This depends upon both the number of terms in each sum and the relative speed of their computation. The second factor is potentially important because although the amplitudes and phases in the new zeta-sum depend purely on elementary functions, they are somewhat more complex than the familiar amplitudes and phases of the RS formula. Let $\tau_{RS}$ be the average CPU time required to calculate a single term of the RS formula and suppose the corresponding average for a term of the new series is $\Omega\tau_{RS}$, where $\Omega > 1$. Using (65), the average CPU time required to calculate the $S$ terms of (63) then becomes

$$\text{CPU time of } S(pc_{CO}) = \left\{\sqrt{\frac{t}{2\pi pc_{CO}}} + \left[\sqrt{\frac{t}{2\pi}\frac{(pc_{CO}+1)}{\sqrt{pc_{CO}}}} - \frac{(INT_O(a)+2)}{2}\right]\Omega\right\}\tau_{RS}, \tag{66}$$

which will be minimised when

$$\frac{dS}{dpc_{CO}} = \sqrt{\frac{t}{2\pi}}\left\{\Omega\frac{\left[pc_{CO}^{1/2} - pc_{CO}^{-1/2}(pc_{CO}+1)/2\right]}{pc_{CO}} - \frac{1}{2pc_{CO}^{3/2}}\right\} = \sqrt{\frac{t}{8\pi pc_{CO}^3}}\{\Omega pc_{CO} - 1 - \Omega\} = 0. \tag{67}$$

Hence the minimal choice of cut-off occurs at $pc_{CO} = 1 + \Omega^{-1}$, in which case

$$L_{CO} = \sqrt{\frac{2t}{\pi(1+\Omega^{-1})}(2+\Omega^{-1})} \qquad \text{and} \qquad N_{CO} = \sqrt{\frac{t}{2\pi(1+\Omega^{-1})}}. \tag{68}$$

Substituting these choices of $L_{CO}$ and $N_{CO}$ into (66) gives rise to the following optimum CPU time for the hybrid as a function of Ω

$$\text{CPU}_{\text{MIN}} \text{ time of } S(pc_{CO}) = \tau_{RS}\sqrt{\frac{t}{2\pi}}\frac{2\left[1+\Omega\left(1-\sqrt{1+\Omega^{-1}}\right)\right]}{\sqrt{1+\Omega^{-1}}} = \tau_{RS}\sqrt{\frac{t}{2\pi}}X(\Omega). \qquad (69)$$

Suppose initially that $\Omega = 1 \Rightarrow pc_{CO} = 2$. Then adopting the hybrid formula (63) with $N_{CO} = \left\lfloor\sqrt{t/4\pi}\right\rfloor$ and $L_{CO} = INT_O\left(3\sqrt{t/\pi}\right)$ means that the computation of the main sum of the RS formula (7, A75) would reduce to summing

$$\sum_{N=1}^{INT\left(\sqrt{t/4\pi}\right)} \frac{2cos\left\{\frac{t}{2}log\left(\frac{t}{2\pi}\right)-\frac{t}{2}-\frac{\pi}{8}-tlog(N)\right\}}{\sqrt{N}} + \sum_{\alpha=INT_O(a)+2}^{INT_O\left(3\sqrt{t/\pi}\right)} \frac{2\sqrt{2}\cos\left((t/2)\{log(pc)+1/pc\}+t/2+\pi/8\right)}{(\alpha^2-a^2)^{\frac{1}{4}}} + T_{a+\varepsilon}(if\ needed) + error, \qquad (70)$$

with the error still bounded by (57). This is equivalent to summing just

$$S = \left\lfloor\sqrt{t/4\pi}\right\rfloor + \frac{INT_O\left(3\sqrt{t/\pi}\right)}{2} - \frac{(INT_O(a)+2)}{2} \approx \sqrt{\frac{t}{2\pi}}\left[2\sqrt{2}-2\right], \qquad (71)$$

terms, compared to the $N_t = \left\lfloor\sqrt{t/2\pi}\right\rfloor$ previously. The factor $2\sqrt{2}-2 = 0.8284\ldots$ implies a reduction of $\left[3-2\sqrt{2}\right]\sqrt{t/2\pi}$ terms, or $\approx 17.16\%$. The error bound is unchanged because an estimate of the partial sum in (6, A64) can be found by finding the difference between (8) integrated around the two circles $C(R_{max})$ and $C(R_{max}-R_t)$. The latter will yield an error term an order of magnitude smaller (cf. note after equation 44) than (57). This encapsulates the simplest manifestation of the hybrid formula, quoted as (5) in the Introduction.

The value of $\Omega = 1$ is obviously an idealisation. Fig.2 shows a plot of the function $100\left(1-X(\Omega)\right)$ as defined by (69), which is equivalent to the percentage CPU time saving produced by the hybrid for values of $\Omega \geq 1$. This shows that substantial CPU savings are still attainable when Ω is significantly larger than unity. So what sort of Ω values are practically possible? A simple illustration of what can easily be achieved employs the following identity for the phase of (6)

$$\frac{t}{2}\left\{log(pc)+\frac{1}{pc}\right\}+\frac{t}{2}+\frac{\pi}{8} = t\left[(\rho^2-1)^{1/2}\left\{(\rho^2-1)^{1/2}-\rho\right\}+log\left\{(\rho^2-1)^{1/2}+\rho\right\}\right]+t+\frac{\pi}{8}, \qquad (72)$$

where $\rho = \alpha/a$. (This identity can easily be established by writing $pc$ in terms of $\rho$ and expanding the phase as a Taylor series about $\rho = 1$.) Utilising (72) one can develop an efficient algorithm for the calculation of the sum of the new series in terms of standard intrinsic functions.

**Algorithm**

Initialise $a := sqrt(8t/\pi)$, $c := sqrt\big(sqrt(8\pi/t)\big)$, $y := (t + \pi/8)_{mod\ 2\pi}$, $s := 0$;

If $|\varepsilon| = |a - NINT_O(a)| < t^{-1/6}$ set $s := T_{a+\varepsilon}$ as given by (A65) or numerical estimate;

For $\alpha = INT_O(a) + 2$† up to $L_{CO}$ step 2, repeat:

$\rho := \alpha/a\,;\ b := sqrt(\rho^2 - 1);\ s := s + cos(t * \{b * (b - \rho) + log(b + \rho)\} + y)/sqrt(b)\,;$

Output: $c * s$

† If $|\varepsilon| < t^{-1/6}$ and $\varepsilon < 0$ then start the sum at $NINT_O(a) + 2$

The advantage of this algorithm is that in calculating the phase of each term by means of (72), one goes some way towards calculating the amplitude, $(\rho^2 - 1)^{-1/4}$, too. Employing this algorithm in a Fortran programme run on a 2.93 GHz Intel Processor PC to compute the sum (5) of the new series from $\alpha = INT_O(a) + 2$ to $L_{CO} = INT_O\left(3\sqrt{t/\pi}\right)$ for $t$ values between $10^6 - 10^{11}$ (see next section and Table I for exact details) and calculating the ratio of the CPU time needed to sum the left hand side of (5), one finds, on average, that $\Omega = 1.30$. Incorporated into the computer algebra package Maple, the above algorithm yields, for similar calculations, an average value of $\Omega = 1.17$. From Fig. 2 these values of Ω imply savings in CPU time utilising the hybrid (63, 68), for corresponding calculations of $Z(t)$, of 14.2% and 15.3% respectively. Alternatively, if one has no *a priori* knowledge of Ω, but simply assumes the optimum cut off occurs at $pc_{CO} = 2$ and utilises (5), values of $\Omega = 1.30\ \&\ 1.17$ would yield savings of 13.5% and 15.1% respectively. Obviously these values will vary depending upon the platform the algorithm is run on and the exact calculations being performed, but they are indicative of the CPU efficiency savings that can be achieved. As it stands the algorithm relies purely on the intrinsic functions *sqrt, log* etc. that are a feature of any programming language. It may be possible to design a "purpose-built" algorithm for computing terms in the new series which reduces Ω closer to unity and pushing up the savings closer to the idealised maximum of 17.16%.

*3.3 Sample computations and practical error bounds*

Computations of $Z(t)$ using the RS formula are (in the main) restricted to Gram points $g_n$ defined by

$$\theta(t = g_n) = n\pi. \tag{73}$$

The observation that, on average, a single zero of $Z(t)$ is usually to be found situated between successive Gram points [20, 23] allows one to verify the Riemann hypothesis for the first $n$ zeros by checking for the expected $n$ sign changes at the various Gram points. (See [28] and [36] for more details on this issue.) Hence it is appropriate to look at the performance of the asymptotic approximation (5) which forms the basis of the idealised hybrid formula (70) at these points. Table I shows a comparison of the errors obtained when calculating the sum of the last $N \in \left\{\left\lfloor\sqrt{t/4\pi}\right\rfloor + 1, N_t\right\}$ terms of the RS formula (with constant contribution to the phase set equal to $\theta(t) - 1/48t$ ... as in equation 70), and the first $\alpha \in \left[INT_O(a) + 2, INT_O\left(3\sqrt{t/\pi}\right)\right]$ odd integer terms of the new series. The computations are each carried out using the algorithm given above across a million Gram points lying above seven successive $t = 10^{6,7,8,9,10,11\ \&\ 12}$ multiples, to give an indication of the approximation's performance at different orders of magnitude. On the rare occasions when $t = g_n$ is such that $a$ happens to lie close to an odd integer, then the $T_{a+\varepsilon}$ transition term at $\alpha = INT_O(a)$ is also included in the new series, either by using (A65) or by direct numerical calculation.

The results shown in Table I clearly illustrate the effectiveness of utilising the first part of the new series to estimate the final part the RS sum. The average errors are $O\left(t^{-1/4}\right)$ and even the largest errors do not stray much beyond this size. These results can be improved upon. The largest errors are found at Gram points at which either $N_{CO}$ or $L_{CO}$ as given by (68), happens to lie someway from their respective rounded up/down integer values. This suggests it would be better to define both $N_{CO}$ and $L_{CO}$ as the *nearest* suitable integers given by (68). The corresponding results using these definitions are also shown in Table I in **bold**. As can be seen this modification has little effect on the average errors but reduces the largest errors by up to 20%. As expected *all* the errors are *very much* smaller values than the absolute upper bound on the error of $6.15t^{-1/12}$ given by (57). The latter could be reduced in two main ways.

The largest contributions to this upper bound arise from three sources. The sum over $N$ of the integral $D(N,R)$ for $N \in \left[N_{2\eta},\ N_t\right]$ (see Appendix B and section 2.4d), the sum of the first order corrections to the estimates $D(N,R)$ arising from the integration through the saddle points (discussed in sections 2.4c-d) and the nature of the transition term (see section 2.6). The relatively crude estimation procedures employed in those sections take no account of the phases multiplying these terms, and by adapting the methodology of exponent pairs ([5],[10] & [11], [17], [21] & [22]) one should be able to increase the magnitude of the power of $t$ which appears in the upper bound somewhat. However, this would require more sophisticated analysis than adopted here.

Alternatively, adopting a heuristic approach, the reduction in the largest errors shown in Table I when $N_{CO}$ and $L_{CO}$ are defined by the *nearest* suitable integers, suggests that a practical error bound can be postulated from the size of the terms in the vicinity of the cut off point. At the idealised cut off $pc_{CO} = 2$, analysis of equation (58) shows that the terms of the RS sum and the new series are distributed in "$pc$ space" on an alternate one-to-one basis, a "halfway house" between distributional regimes shown in Figs.1a & 1b. Consequently, if the cut-off is made at or near $pc_{CO} = 2$, then the error in the hybrid is *extremely unlikely* to exceed the magnitude of the larger of those two terms from both series which lie closest in proximity to the cut-off point. At the idealised cut-off $pc_{CO} = 2$, the magnitudes of both series are the same because

| Range of $t$ $t_{start} - t_{end}$ | The first $10^6$ Gram points above $t_{start}$ | Average $\|error\|$ for the $10^6$ Gram points | Average exponent $s$, so that $\|error\| = t^s$ | Max$\|error\|$ and associated Gram point | $\left(\frac{64\pi}{t}\right)^{1/4}$ |
|---|---|---|---|---|---|
| $\times 10^6$ 1.0 – 1.5150284 | $g_{1747145} - g_{2747145}$ | $3.72 \times 10^{-2}$ **$3.54 \times 10^{-2}$** | $-0.263$ **$-0.267$** | $1.318 \times 10^{-1}, g_{1757119}$ **$1.044 \times 10^{-1}, g_{1834822}$** | $1.191 \times 10^{-1}$ |
| $\times 10^7$ 1.0 – 1.04393264 | $g_{21136125} - g_{22136125}$ | $2.20 \times 10^{-2}$ **$2.09 \times 10^{-2}$** | $-0.261$ **$-0.266$** | $7.347 \times 10^{-2}, g_{21585534}$ **$5.873 \times 10^{-2}, g_{22099951}$** | $6.696 \times 10^{-2}$ |
| $\times 10^8$ 1.0 – 1.003788547 | $g_{248008024} - g_{249008024}$ | $1.25 \times 10^{-2}$ **$1.18 \times 10^{-2}$** | $-0.260$ **$-0.263$** | $4.158 \times 10^{-2}, g_{248610706}$ **$3.244 \times 10^{-2}, g_{248485478}$** | $3.766 \times 10^{-2}$ |
| $\times 10^9$ 1.0 – 1.0003326980 | $g_{2846548032} - g_{2847548032}$ | $7.01 \times 10^{-3}$ **$6.66 \times 10^{-3}$** | $-0.259$ **$-0.261$** | $2.313 \times 10^{-2}, g_{2846691368}$ **$1.703 \times 10^{-2}, g_{2847140443}$** | $2.118 \times 10^{-2}$ |
| $\times 10^{10}$ 1.0 – 1.00002965447 | $g_{32130158314} - g_{32131158314}$ | $3.76 \times 10^{-3}$ **$3.52 \times 10^{-3}$** | $-0.257$ **$-0.260$** | $1.216 \times 10^{-2}, g_{32130991688}$ **$9.997 \times 10^{-3}, g_{32131011467}$** | $1.191 \times 10^{-2}$ |
| $\times 10^{11}$ 1.0 – 1.000002674767 | $g_{357948363082} - g_{357949363082}$ | $2.26 \times 10^{-3}$ **$2.13 \times 10^{-3}$** | $-0.254$ **$-0.256$** | $6.687 \times 10^{-3}, g_{357949027704}$ **$6.709 \times 10^{-3}, g_{357949027219}$** | $6.696 \times 10^{-3}$ |
| $\times 10^{12}$ 1.0 – 1.0000002435992 | $g_{3945951430271} - g_{3945952430271}$ | $1.62 \times 10^{-3}$ **$1.21 \times 10^{-3}$** | $-0.243$ **$-0.259$** | $3.341 \times 10^{-3}, g_{39459523302550}$ **$2.377 \times 10^{-3}, g_{3945951434899}$** | $3.766 \times 10^{-3}$ |

Table I. Comparison of the errors obtained when using the sum of the first $\alpha \in \left[INT_O(a) + 2, \boldsymbol{N}/INT_O\left(3\sqrt{t/\pi}\right)\right]$ odd integer terms of the new series (6, A64-65), to estimate the sum of the last $N \in \left[\boldsymbol{N}/INT\left(\sqrt{t/4\pi}\right) + 1, N_t\right]$ of the RS main sum as formulated in (5, 70). The computations were carried out at a million Gram-points, starting at $t = 10^{6,7,8,9,10,11\ \&\ 12}$ respectively. Values of the average and maximum errors over these Gram-point ranges are shown.

$$\frac{2\sqrt{2}}{(L_{CO}^2 - a^2)^{\frac{1}{4}}} = \frac{2^{5/4}\pi^{1/4}(pc_{CO})^{1/4}}{t^{1/4}\sqrt{pc_{CO}-1}} = \left(\frac{64\pi}{t}\right)^{1/4} = \frac{2}{\sqrt{a/4\sqrt{pc_{CO}}}} = \frac{2}{\sqrt{N_{CO}}}. \qquad (74)$$

Hence one would expect that $(64\pi/t)^{1/4}$ should form a practical upper bound on the error in this instance. The calculations presented in Table I support this conclusion. If $N_{CO}$ and $L_{CO}$ are defined to be *nearest* suitable integers from expressions (68), just two of the seven million computations shown give rise to an error exceeding $(64\pi/t)^{1/4}$ (at Gram points $g_{3579490272188\&9}$ and then only by less than 0.2%). If $N_{CO}$ and $L_{CO}$ are defined by rounding up and down (68) then errors greater than $(64\pi/t)^{1/4}$ are observed on less than one hundred occasions. It should be emphasized that these results are based purely on numerical observation and are in no way a proof that (74) forms the basis of a *de facto* upper bound on the error. But they are indicative that exceptions are likely to be very rare. It is also indicative that the $O(t^{-1/12})$ error bound, established in the proof of (7), is only this large because the transition term $T_{a+\varepsilon}$ is not included in (6). For the general case, when the cut-off is given by $pc_{CO} = 1 + \Omega^{-1}$, equation (74) suggests a practical upper bound on the error equivalent to $(32\Omega(1+\Omega)\pi/t)^{1/4}$.

## 4. Conclusions

Utilising the RS integral as a starting point, this paper shows that it is possible derive a new *form* (as articulated below) of zeta-sum (6, A64) for calculating Hardy's $Z$ function along the critical line. Utilising just the most significant contributions to the new zeta sum it is possible to prove that the error arising from this approximation will never exceed a term of $O(t^{-1/12})$. In its entirety this new zeta-sum is a computationally inefficient means for calculating $Z(t)$. However, having in the course of the error analysis established a distributional link between the terms of the new zeta-sum and the terms of the RS formula, it proves possible to postulate a means of increasing the computational efficiency of calculating $Z(t)$. This can be done by means of a hybrid formula combining partial sums of the RS terms with terms from the new series, as summarised by (63 & 68). The source of the increase in computational efficiency lies in the fact that the sum of the final $(1/\sqrt{2} - 1/2)\sqrt{t/\pi}$ terms of the RS formula can be calculated using series (5), which contains just $(3/2 - \sqrt{2})\sqrt{t/\pi}$ terms. Potentially this means that $Z(t)$ can be calculated by summing some ≈ 17.16% less terms than by employing the RS formula alone. The slightly more complicated make up of the new series terms means that this reduction translates into a saving of ≈ 15% in CPU time. Given that since its discovery in 1932 the RS formula was thought to be the most efficient method possible for calculating $\zeta\left(\frac{1}{2} + it\right)$ and given the scrutiny it has undergone, these results represent an important advance.

Naturally after such an advance there are a number of strands of future research effort. The error term associated with (5) has been conservatively estimated never to exceed $6.15/t^{1/12}$. From the point of view of employing the hybrid formula to calculate $Z(t)$ at Gram points in order to verify the Riemann hypothesis over larger ranges of $t$, with calculations currently (at the time of writing) confirming its veracity beyond the first 10 trillion zeros ($t{\sim}2.5 \times 10^{13}$), this error bound is already sufficiently small to prove helpful in future calculations. (In most instances the error only needs to be small enough to confirm the sign of $Z(t)$ at a Gram point to extend the range of verification.) However, as the results of Table I indicate, this error bound represents a considerable overestimate, and in practice the error only exceeds $(64\pi/t)^{1/4}$ on very rare occasions. Further investigation could be directed towards refining this observation, to obtain a more precise estimate of the error term associated with (5) for any value of $t$. Another interesting observation is the behaviour of the new series when very close to the point $\alpha = a$, where it effectively acts a smooth asymptotic cut off of the RS formula as $N \to N_t$, similar in substance, if not in form, to the asymptotic cut-off devised by [4]. But this new formulation seems a much more aesthetically appealing way of prescribing such a cut-off, predicting as it does that the sum of the last $N = (N_\mu, N_t]$ RS terms will rarely amount to anything significant. Exceptions to this rule occur periodically whenever $a$ lies very close to an odd integer, when these final terms suddenly combine to a remarkable degree to give a sum of $O\left(t^{-1/12}\right)$ (q.v. the discussion at the end of Section 2.6).

An important focus of research would be to establish the possibility of calculating the new series for multiple values $t = k\tau$ with $k \in [K_1, K_2]$ in a total of just $O\left(\left(2\sqrt{2}-2\right)N_t Log\{K_2 - K_1\}\right)$ operations, utilising the method of non-linear FFT [13] as outlined in [12]. This technique, originally developed by [31], underlies the modern method of 'fast' calculations of $Z(t)$ at many different Gram points employed for Riemann hypothesis verification. If so, the hybrid could be utilised to speed up such calculations. The non-linear FFT technique may not be so easy to apply here though because the terms of the new series, although elementary, are functions of *scaled* odd integers $\rho = \alpha/a$ (see the algorithm in 3.2), not *pure* integers as in the RS formula.

Beyond the computational aspects of the work there are some important theoretical considerations. First, from the point of view of completeness (although this is not necessary for the validation of the new zeta-sum and the hybrid, which relies purely on the proof of the *Main Theorem* devised in Section 2 & Appendix B), it would good to have a proof of the exact cancellation of all the terms of all the terms of $O\left(e^{\alpha\sqrt{\pi t}/\sqrt{2}}\right)$ in the calculation of the imaginary part of the RSI (as asserted in sections A3.8 and Appendix C). Second, it would be interesting to know how the new zeta-sum (6, A64) can be formulated for $z = \sigma + \mathrm{i}t$ off the critical line with $1/2 < \sigma < 1$. This may just prove to be as simple substituting $t \equiv t + \mathrm{i}(1/2 - \sigma)$ through most of the above analysis. Finally there is the idea that (6, A64) is a zeta-sum of a new form. To substantiate this idea, consider the *symmetries* associated with (6, A64). To begin there is the fact that one can replace $\alpha$ by $-\alpha$, compute the sum over the negative odd integers and still obtain the same answer for $Z(t)$. In one sense this is an obvious mathematical artefact of the formulation of the RSI as the sum (A21), in which $\alpha$ can be positive or negative. However, no corresponding natural symmetry over the summands is apparent from any other formulation of $Z(t)$, certainly not one based on the Dirichlet series of the zeta function. Another interesting symmetry becomes apparent when one appreciates that the terms of the new zeta-sum can also arise when definition (A47) for $pc(\alpha\ )$, is replaced by the corresponding one using the negative instead of the positive square root. This gives values for $pc(\alpha) \in (0, 1)$, because one is effectively replacing $pc(\alpha)$ by its reciprocal. But the terms of the new zeta-sum, $Re\{(pc)^{1/4} exp[\mathrm{i}\{t/2\,(log(pc) + 1/pc) + t/2 + \pi/8\}]/\sqrt{pc - 1}\}$, remain the same, irrespective of which definition of $pc(\alpha)$ is employed. So in '$pc$ space' one finds the terms of (6, A64) repeated twice as poles along the real axis in the intervals $(0, 1)$ and $(1, \infty)$. This maybe another artefact of the $\alpha \leftrightarrow -\alpha$ symmetry as noted above, although it is not necessarily exact in the sense that the terms of (6, A64) arise from the approximation of integral (A44). The geometry of the respective $\alpha \leftrightarrow pc \leftrightarrow N$ transformations in this regard is somewhat curious and worthy of further exploration. It also would be interesting to discover if corresponding zeta-sums of this new type exist for generalised Dirichlet $L$-series, and if so what properties they possess.

## References

[1] Abramowitz, M. & Stegun, I. A. 1972 *Handbook of Mathematical Functions*, New York, Dover Publications.

[2] Bender, C. M. & Orszag, S. A. 1999 *Advanced Mathematical Methods for Scientists and Engineers*, New York, Springer-Verlag.

[3] Berry, M. V. 1995 The Riemann-Siegel expansion for the zeta function: high orders and remainders. *Proc. R. Soc. Lond. A* **450,** 439-462.

[4] Berry, M. V. & Keating, J. P. 1992 A new asymptotic representation for $\zeta\left(\frac{1}{2}+\mathrm{i}t\right)$ and quantum spectral determinants. *Proc. R. Soc. Lond.* A **437**, 151-173.

[5] Bombieri, E. & Iwaniec, H. 1986 Some mean value theorems for exponential sums. *Ann. Scuola Norm. Sup. Pisa Cl. Sci.* (4) **13**, 473-486.

[6] Borwein, J. M., Bradley, D. M. & Crandell, R. E. 2000 Computational strategies for the Riemann zeta function. *J. Comp. Appl. Math.* **121**, 247-296.

[7] Brent, R. P. 1979 On the zeros of the Riemann zeta function in the critical strip. *Math. Comp.* **33**, 1361-1372.

[8] de Bruijn, N. G. 1981 *Asymptotic Methods in Analysis,* New York, Dover Publications.

[9] Cohen, H., Rodriguez-Villegas, F. & Zagier, D. 2000 Convergence acceleration of alternating series. Exp. Math. **9**, 3-12.

[10] Van der Corput, J. G. 1921 Zahlentheoretische Abschätzungen, *Math. Ann.* **84**, 53-79.

[11] Van der Corput, J. G. 1922 Verschärfung der Abschätzung beim Teilerproblem, *Math. Ann.* **87**, 39-65.

[12] Crandall, R. & Pomerance, C. 2005 *Prime Numbers-A Computational Perspective*, Springer Science & Business Media, Inc. (USA).

[13] Dutt, A. & Rokhlin, V. 1993 Fast Fourier Transforms for Nonequispaced Data. *SIAM J. Sci. Comput.* **14**, 1368-1393.

[14] Edwards, H. M. 1974 *Riemann's zeta function*, New York, Dover Publications.

[15] Gabcke, W. 1979 Neue Herleitung und Explizite Restabschätzung der Riemann-Siegel-Formel. Ph.D. thesis, Göttingen.

[16] Gradshteyn, I. S. & Ryzhik, I. M. 2007 *Table of Integrals, Series and Products*. ($7^{th}$ edition) Jeffrey, A. & Zwillinger, D. (eds.) Elsevier Academic Press.

[17] Graham, S. W. & Kolesnik, G. 1991 *Van der Corput's Method of Exponential Sums*, LMS Lecture Notes 126, Cambridge University Press.

[18] Gram, J. P. 1903 Sur les Zéros de la Fonction $\zeta(s)$ de Riemann. *Acta Math*. **27**, 289-304.

[19] Haselgrove, C. B. & Miller J. C. P. 1960 Tables of the Riemann zeta Function. *Roy. Soc. Math. Tables*, **6**, Cambridge, Cambridge University Press.

[20] Hutchinson, J. I. 1927 On the roots of the Riemann zeta-function, V. *Trans. Amer. Math. Soc.* **27**, 27-49.

[21] Huxley M. N. 1993 Exponential sums and the Riemann zeta function IV. *Proc. Lond. Math. Soc.* **66**, 1-40.

[22] Huxley M. N. 2005 Exponential sums and the Riemann zeta function V. *Proc. Lond. Math. Soc.* **90**, 1-41.

[23] Ivić A. 2013 *The Theory of Hardy's Z-Function*. Cambridge Tracts in Mathematics, Cambridge University Press.

[24] Keating, J. P. 1992 Periodic orbit resummation and the quantization of chaos. *Proc. R. Soc. Lond.* A **436**, 99-108.

[25] Kuzmin, R. 1934 On the roots of the Dirichlet series. *Izv. Akad. Nauk SSSR Ser. Math. Nat. Sci.* **7**, 1471-1491.

[26] Kuznetsov A. 2007 On the Riemann-Siegel formula. *Proc. R. Soc. Lond.* A **463**, 2557-2568.

[27] Luke, Y. L. 1969 *The special functions and their approximations (Vol. 1),* New York & London, Academic Press.

[28] van de Lune, J., te Riele, H. J. J. & Winter, D. T. 1986 On the zeros of the zeta function in the critical strip. IV, *Math. Comp.* **46**, 667-681.

[29] Nardin, M., Perger, W. F. & Bhalla, A. 1992 Numerical evaluation of the confluent hypergeometric function for complex arguments of large magnitudes. *J. Comp. Appl. Math.* **39**, 193-200.

[30] Odlyzko, A. M. 2001 The $10^{22}$ zero of the Riemann zeta function. *Dynamical, Spectral and Arithmetic Zeta Functions*, M. van Frankenhuysen and M. L. Lapidus (eds.) Amer. Math. Soc., Contemporary Math. **290**, 139-144.

[31] Odlyzko, A. M. & Schönhage A. 1988 Fast algorithms for multiple evaluations of the Riemann zeta function. *Trans. Amer. Math. Soc.* **309**, 797-809.

[32] Paris, R. B. 1994 An asymptotic representation for the Riemann zeta function on the critical line. *Proc. R. Soc. Lond. A* **446**, 565-587.

[33] Siegel, C. L. 1932 Über Riemanns Nachlaß zur analytischen Zahlentheorie. *Quellen Studien zur Geschichte der Math. Astron. Und Phys. Abt. B: Studien* **2**, 45-80.

[34] Titchmarsh, E. C. 1986 *The theory of the Riemann zeta-function,* 2nd edn. Oxford, Clarendon Press.

[35] Tricomi, F. G. 1954 *Funzioni ipergeometriche confluenti* Éd. Cremonese, Rome.

[36] Trudgian, T. S. 2011 On the success and failure of Gram's law and the Rosser rule. **143**, 225-256.

[37] Turing, A. M. 1943 A method for the calculation of the zeta-function. *Proc. London Math. Soc. Ser.2* **48**, 180-197.

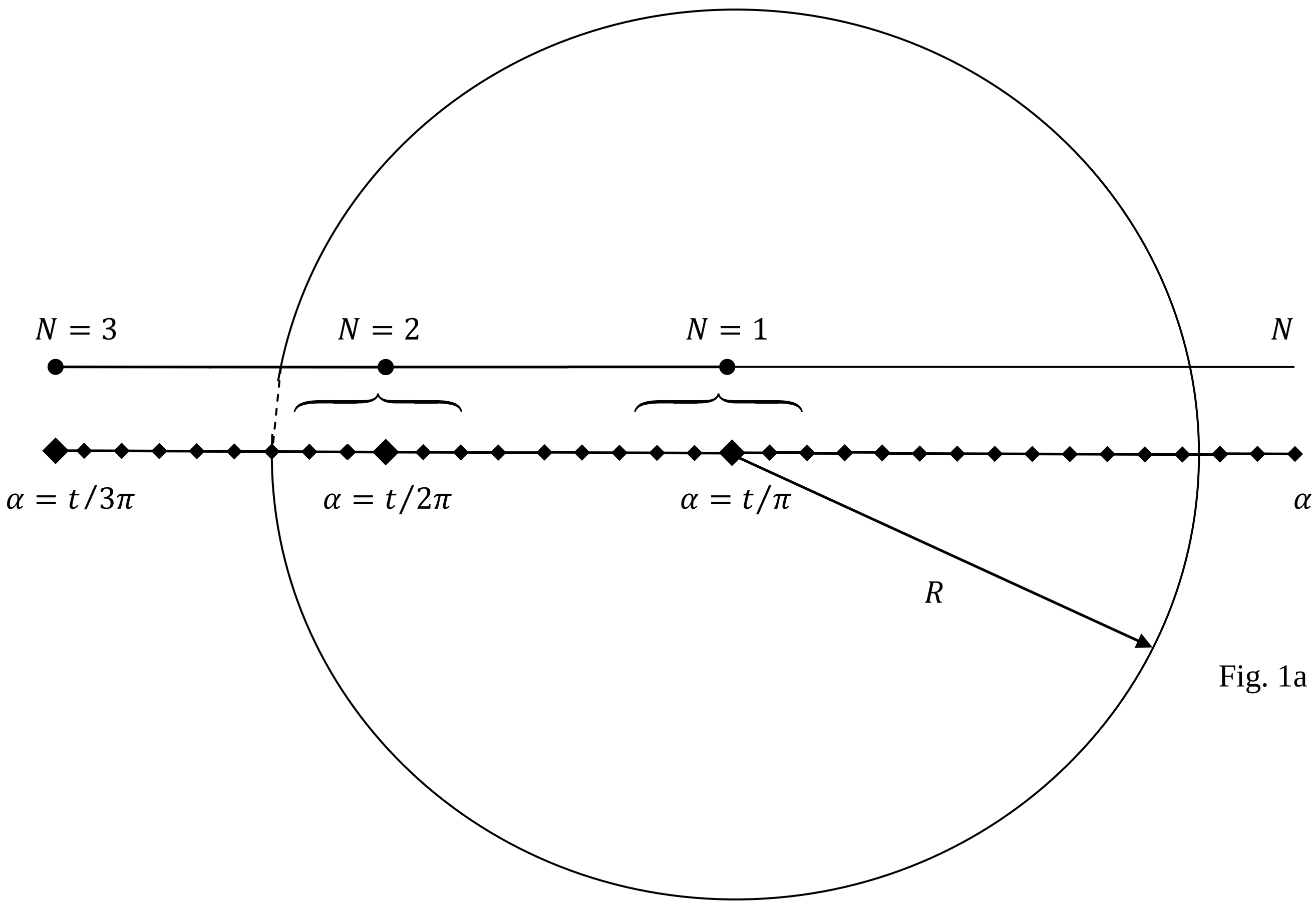

Fig. 1a

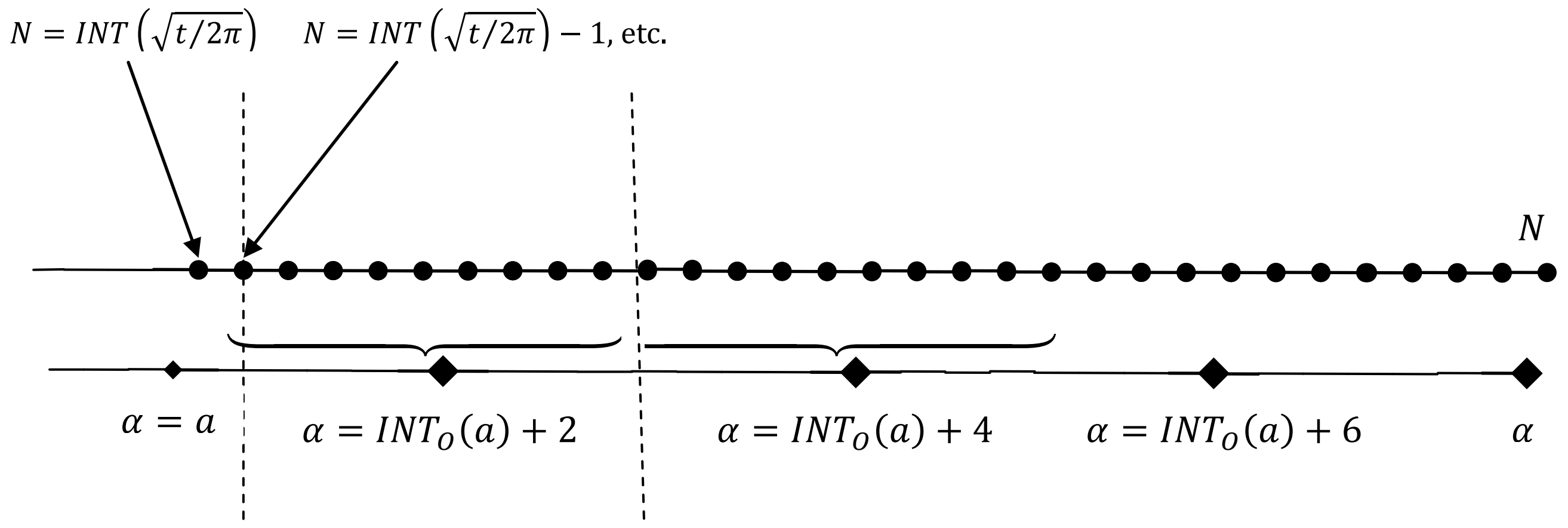

Sections of Circumferences of Circles

of radii $R$ centred at $\alpha = t/\pi$.

Fig. 1b

Figure 1. Schematic of the distribution of the saddle points $N = 1,2 \dots INT\left(\sqrt{t/2\pi}\right)$ of the integral $D(N,R)$ (26), which give rise to the terms of the RS formula and their relation to the terms of the new zeta-sum (6, A64) in '$\alpha$ space'. In Fig. 1a when $\alpha = O(t/\pi)$, the low numbered saddle points are distributed far apart at approximate intervals of $\alpha = t/N\pi$. This means that the terms of the new series within about $O\left(t^{1/3}\right)$ of $\alpha = t/N\pi$ sum together to form the respective RS terms. In Fig.1b when $\alpha \approx a$ this distributional pattern is reversed, with now $O\left(t^{1/4}\right)$ saddle points squeezed between each successive odd integer $\alpha$. This means that in this region each term of the new series is representative of the sum of $O\left(t^{1/4}\right)$ terms of the RS formula.

Fig. 2

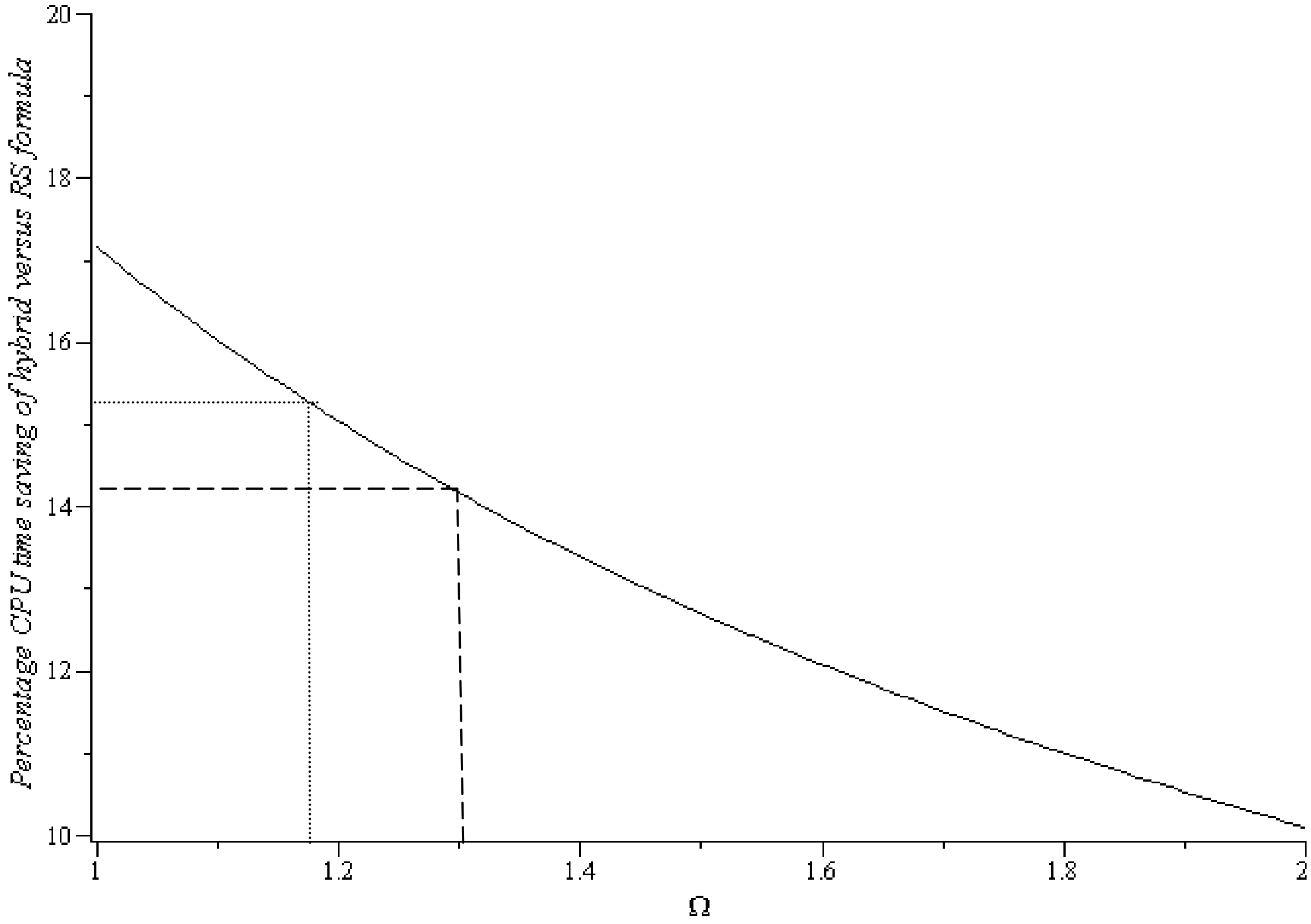

Figure 2. Plot of the percentage CPU time savings possible by substituting the hybrid formula (63 & 68) in place of the RS formula against Ω, the average CPU time ratio for the computation of single terms from each formulae. Employing the algorithm listed in section 3.2, typical values of $\Omega \approx 1.17 - 1.30$, imply CPU time savings of about 14-15%.

# Appendix A. Derivation of the new Zeta-Sum from the RS Integral

## A1. Expansions of some relevant integrals

### A1.1 *General remarks and estimation strategy*

One seeks an asymptotic expansion of the integral in (2), which for $s = \frac{1}{2} + \mathrm{i}t$ (here $t$ will assumed to be real and positive, although this initial analysis can be adjusted so that $s$ lies off the critical line) can be rewritten as

$$\mathrm{RSI}{=}\mathrm{RSI}_{-} + \mathrm{RSI}_{+} = \frac{1}{2\mathrm{i}}\int_{-\infty+\mathrm{i}\infty}^{0+\mathrm{i}/2} \frac{e^{-\mathrm{i}\pi z^2} z^{-1/2-\mathrm{i}t}}{\sin(\pi z)} dz + \frac{1}{2\mathrm{i}}\int_{0+\mathrm{i}/2}^{\infty-\mathrm{i}\infty} \frac{e^{-\mathrm{i}\pi z^2} z^{-1/2-\mathrm{i}t}}{\sin(\pi z)} dz, \qquad \text{(A1)}$$

where the path integration is a line of slope $-1$ definitely chosen to pass through the point $z = 1/2$, midway between the poles at 0 and 1, i.e. $z = 1/2 + qe^{-\mathrm{i}\pi/4}, q \in (-\infty,\ \infty)$. The second term in the numerator is defined in the usual way, viz.

$$\begin{aligned} z^{-1/2-\mathrm{i}t} &= exp[(-1/2 - \mathrm{i}t)\{log|z| + \mathrm{i}Arg(z)\}] \\ &= exp[(-1/2 - \mathrm{i}t)\{log|z| + \mathrm{i}arctan(Im(z)/Re(z))\}], \qquad Re(z) > 0 \\ &= exp\big[(-1/2 - \mathrm{i}t)\big\{log|z| + \mathrm{i}\{\pi + arctan(Im(z)/Re(z))\}\big\}\big]. \quad Re(z) < 0 \end{aligned} \qquad \text{(A2)}$$

This adjustment in the definition of (A2) to ensure analytic continuation as the integration path crosses the imaginary axis at $q = -1/\sqrt{2}$, will, eventually, prove to be of some significance (see Appendix C). A brief overview of the characteristics of the integrand (A1) along the designated integration path reveals it reaches a maximum amplitude when $q \approx -t^{1/3}/\left(2^{\frac{5}{6}}\pi^{\frac{1}{3}}\right)$ (of about $e^{3\pi t/4+(\pi t^2/32)^{1/3}}$). A naïve initial strategy might be to consider local expansions of the integrand around this maximum to estimate the integral. However, this fails because of the extremely rapid sinusoidal oscillations (whose period decreases as $t$ increases) that arise from (A2), which have the effect of cancelling out the influence of the maximum amplitude on the final result of the overall integral. Rather the influence of these oscillations is more akin to that of a delta function, averaging out the contribution to the integral as one travels down the integration path and leaving its value heavily dependent on the nature of the integrand any designated cut-off point. As the point $q = -1/\sqrt{2}$ provides a natural break in the path of integration (because the *arctan* function is not defined at $Re(z) = 0$) it is logical to investigate the behaviour of the integral in two parts, $\mathrm{RSI}_{-}$ and $\mathrm{RSI}_{+}$, as expressed in equation (A1). For $\mathrm{RSI}_{-}$ the integration covers the range $q \in \left(-\infty, -1/\sqrt{2}\right)$, whilst for $\mathrm{RSI}_{+}$ the integration range covers $q \in \left(-1/\sqrt{2}\,, \infty\right)$. The strategy will then be to see how these results can be matched up in a sensible way. This matching process gives rise to a modified interpretation of (A1) in terms of an infinite sum of a certain class of somewhat simpler integrals, which in turn provides the key for establishing the asymptotic expansion of the RSI. To see how this comes about one first needs to establish a result concerning two relatively simple integrals which provides the foundation for the analysis that follows.

A1.2 *An identity for a relatively simple class of integrals*

Consider the generic integral, denoted by $I_n$, defined by

$$I_n = \int_{-\infty+\mathrm{i}/2\sqrt{2}}^{w+\mathrm{i}/2\sqrt{2}} z^n e^{-\lambda z^2} e^{\mathrm{i}\beta Log(z^2)} dz = \int_{-\infty+\mathrm{i}/2\sqrt{2}}^{w+\mathrm{i}/2\sqrt{2}} \frac{e^{\mathrm{i}\beta Log(z^2)}}{z} z^{n+1} e^{-\lambda z^2} dz. \quad \text{(A3)}$$

Here $n$ is a non-negative integer (although this work will be particularly concerned with the two ‘foundation’ cases $n = 0$ or 1); $w \in \mathbb{R}, \beta \in \mathbb{C}$ and $\lambda > 0$ . The appearance of the $e^{-\lambda z^2}$ term in the integrand guarantees the convergence of $I_n$ irrespective of the choice of limit $w$. Integrating by parts gives

$$I_n = \left[\frac{z^{n+1} e^{-\lambda z^2} e^{\mathrm{i}\beta Log(z^2)}}{2\mathrm{i}\beta}\right]_{-\infty+\mathrm{i}/2\sqrt{2}}^{w+\mathrm{i}/2\sqrt{2}} - \frac{(n+1)}{2\mathrm{i}\beta} I_n + \frac{\lambda}{\mathrm{i}\beta} \int_{-\infty+\mathrm{i}/2\sqrt{2}}^{w+\mathrm{i}/2\sqrt{2}} z^{n+2} e^{-\lambda z^2} e^{\mathrm{i}\beta Log(z^2)}\, dz. \quad \text{(A4)}$$

Assuming differentiation under the integral sign is valid, this gives

$$I_n = \frac{\left(w+\mathrm{i}/2\sqrt{2}\right)^{n+1}}{2\mathrm{i}\beta} e^{-\lambda\left(w+\mathrm{i}/2\sqrt{2}\right)^2} e^{\mathrm{i}\beta Log\left(\left(w+\mathrm{i}/2\sqrt{2}\right)^2\right)} - \frac{(n+1)}{2\mathrm{i}\beta} I_n - \frac{\lambda}{\mathrm{i}\beta}\frac{dI_n}{d\lambda},$$

$$\Rightarrow \quad \frac{dI_n}{d\lambda} + \left[\frac{n+1}{2\lambda} + \frac{\mathrm{i}\beta}{\lambda}\right] I_n = \frac{\left(w+\mathrm{i}/2\sqrt{2}\right)^{n+1}}{2\lambda} e^{-\lambda\left(w+\mathrm{i}/2\sqrt{2}\right)^2} e^{\mathrm{i}\beta Log\left(\left(w+\mathrm{i}/2\sqrt{2}\right)^2\right)}. \quad \text{(A5)}$$

This is obviously a standard linear first order ODE, which one can solve by multiplying through by the integrating factor $\lambda^{(n+1)/2+\mathrm{i}\beta}$ giving

$$\lambda^{(n+1)/2+\mathrm{i}\beta} I_n = \frac{\left(w+\mathrm{i}/2\sqrt{2}\right)^{n+1}}{2} e^{\mathrm{i}\beta Log\left(\left(w+\mathrm{i}/2\sqrt{2}\right)^2\right)} \int \lambda^{(n-1)/2+\mathrm{i}\beta} e^{-\lambda\left(w+\mathrm{i}/2\sqrt{2}\right)^2} d\lambda + C_n,$$

$$\Rightarrow I_n = \frac{\left(w+\mathrm{i}/2\sqrt{2}\right)^{n+1}}{2\lambda^{(n+1)/2+\mathrm{i}\beta}} e^{\mathrm{i}\beta Log\left(\left(w+\mathrm{i}/2\sqrt{2}\right)^2\right)} \int \lambda^{(n-1)/2+\mathrm{i}\beta} e^{-\lambda\left(w+\mathrm{i}/2\sqrt{2}\right)^2} d\lambda + C_n \lambda^{-(n+1)/2-\mathrm{i}\beta}.$$

(A6)

Here $C_n$ is a constant of integration that is independent of both the variables $\lambda$ and $w$. If $|\beta| \gg \lambda$, it is possible to develop an asymptotic approximation‡ for $I_n$ in reciprocal powers of $\beta$. This can be done by writing the exponential term in the integrand of (A6) as a series, integrating, rewriting $((n-m)/2 + \mathrm{i}\beta)^{-1}$ as geometric series and then summing the terms of $O(\beta^{-1}, \beta^{-2}, \ldots)$. However, beyond the fact that to first order $I_n \sim \lambda^{n+1} e^{-\lambda w^2 + \mathrm{i}\beta Log(w^2)}/2\mathrm{i}\beta$, the details of this turn out to be of little relevance for the RSI, and identity (A6) can be utilised as it stands.

‡These details are available from the author on request.

A1.3 *Exponential Sums of the Integrals* $I_n$

Consider the somewhat more complicated integral denoted by $Q$, which can be written as an exponential sum over the $I_n$ as follows:

$$Q=\int_{-\infty+\mathrm{i}/2\sqrt{2}}^{w+\mathrm{i}/2\sqrt{2}} e^{\sigma z-\lambda z^2} e^{\mathrm{i}\beta Log(z^2)} dz = \sum_{n=0}^{\infty}\frac{\sigma^n}{n!}\int_{-\infty+\mathrm{i}/2\sqrt{2}}^{w+\mathrm{i}/2\sqrt{2}} z^n\, e^{-\lambda z^2} e^{\mathrm{i}\beta Log(z^2)} dz = \sum_{n=0}^{\infty}\frac{\sigma^n}{n!} I_n \qquad \text{(A7)}$$

Here $\sigma \in \mathbb{C}$, and because the exponential sum is absolutely convergent for all $|\sigma|$ and as the integrals $I_n$ increase only algebraically as $\lambda^{n+1}$, the interchange of the order of summation and integration will be valid. In what follows it will prove more useful to write $e^{\sigma z} = sinh(\sigma z) + cosh(\sigma z)$, giving

$$Q = \sum_{n=0}^{\infty}\frac{\sigma^{2n}}{(2n)!} I_{2n} + \sum_{n=0}^{\infty}\frac{\sigma^{2n+1}}{(2n+1)!} I_{2n+1}. \qquad \text{(A8)}$$

Now the various integrals $I_n$ can be derived from the foundation integrals $I_0$ and $I_1$ simply by means of successive differentiation. Depending upon whether $n$ is odd or even, one has

$$I_n = \left(-\frac{d}{d\lambda}\right)^{n/2} I_0, \ \text{if } n \text{ is even}, \qquad I_n = \left(-\frac{d}{d\lambda}\right)^{(n-1)/2} I_1, \ \text{if } n \text{ is odd}, \qquad \text{(A9)}$$

with the now standard caveat that the interchanging of the order of the differentiation and integration remains a valid operation. Hence (A8) becomes

$$Q = \sum_{n=0}^{\infty}\frac{\sigma^{2n}}{(2n)!}\left(-\frac{d}{d\lambda}\right)^{n} I_0 + \sum_{n=0}^{\infty}\frac{\sigma^{2n+1}}{(2n+1)!}\left(-\frac{d}{d\lambda}\right)^{n} I_1. \qquad \text{(A10)}$$

Finally on substituting identity (A6) for the two foundation integrals into (A10), gives the following expression for $Q$

$$Q = \sum_{n=0}^{\infty}\frac{\sigma^{2n}}{(2n)!}\left(-\frac{d}{d\lambda}\right)^{n}\left[\frac{\left(w+\mathrm{i}/2\sqrt{2}\right)}{2\lambda^{1/2+\mathrm{i}\beta}} e^{\mathrm{i}\beta Log\left(\left(w+\mathrm{i}/2\sqrt{2}\right)^2\right)}\int \lambda^{-1/2+\mathrm{i}\beta} e^{-\lambda\left(w+\mathrm{i}/2\sqrt{2}\right)^2} d\lambda + C_0\lambda^{-1/2-\mathrm{i}\beta}\right]$$
$$+\sum_{n=0}^{\infty}\frac{\sigma^{2n+1}}{(2n+1)!}\left(-\frac{d}{d\lambda}\right)^{n}\left[\frac{\left(w+\mathrm{i}/2\sqrt{2}\right)^2}{2\lambda^{1+\mathrm{i}\beta}} e^{\mathrm{i}\beta Log\left(\left(w+\mathrm{i}/2\sqrt{2}\right)^2\right)}\int \lambda^{\mathrm{i}\beta} e^{-\lambda\left(w+\mathrm{i}/2\sqrt{2}\right)^2} d\lambda + C_1\lambda^{-1-\mathrm{i}\beta}\right].$$

(A11)

Substituting the asymptotic approximations developed from identity (A6) for the foundation integrals in place of the terms in square brackets in (A11), leads to an asymptotic approximation for integral $Q$ in terms of $w$ and reciprocal powers of the (large) parameter $\beta$. But although useful for validation tests of the assumptions leading to (A11), these details again prove to be of no relevance for the RSI. The next stage of the process is to relate the integral $Q$ to the RSI as defined by (A1). This is will be done in some detail for the first part of (A1), namely $\mathrm{RSI}_-$. The corresponding development for $\mathrm{RSI}_+$ is very similar and most of the detail will be omitted.

## A2. The expansion of the Riemann-Siegel Integral as a Sum

A2.1. *The $RSI_-$ along that part of the contour lying in the region with $Re(z) < 0$.*

Making the substitution $z = 1/2 + qe^{-i\pi/4}, q \in \left(-\infty,\ -1/\sqrt{2}\right)$ in (A1) and utilising (A2) yields, after some algebra

$$e^{-i(1/2+it)\pi}e^{-i\pi/4}\int_{-\infty}^{-1/\sqrt{2}} \frac{e^{-\pi q^2-\pi(1+i)q/\sqrt{2}}e^{(-1/4-it/2)log\left(1/8+(q+1/2\sqrt{2})^2\right)}e^{-i(1/2+it)arctan\left(-q/(q+1/\sqrt{2})\right)}}{2isin[\pi(1/2+qe^{-i\pi/4})]}e^{-i\pi/4}dq.$$

(A12)

Now because $q < 0$ in (A12), one can express the denominator as an (absolutely) convergent sum as follows:

$$\frac{1}{sin[\pi(1/2+qe^{-i\pi/4})]} = \frac{1}{cos\left[\pi q\,(1-i)/\sqrt{2}\right]} = \frac{2}{e^{i\pi q(1-i)/\sqrt{2}} + e^{-i\pi q(1-i)/\sqrt{2}}}$$
$$= \frac{2e^{\pi q(1+i)/\sqrt{2}}}{\left[1+e^{2\pi q(1+i)/\sqrt{2}}\right]}, = 2\sum_{\substack{\alpha=1\\ odd}}^{\infty}(-1)^{(\alpha-1)/2}e^{\pi q(1+i)\alpha/\sqrt{2}}. \qquad \text{(A13)}$$

Substituting (A13) into (A12) and making the change of variables $x = q + 1/2\sqrt{2}$ gives

$$\frac{1}{i}\sum_{\substack{\alpha=1\\ odd}}^{\infty}(-1)^{(\alpha-1)/2}\, e^{-i(1/2+it)\pi}e^{-i\pi/2-\pi/8}e^{-(\alpha-1)(1+i)\pi/4}$$

$$\times\int_{-\infty}^{-\frac{1}{2\sqrt{2}}} e^{-\pi x^2}e^{\left(\alpha+i(\alpha-1)\right)\pi x/\sqrt{2}}\, e^{(-1/4-it/2)log\left(1/8+x^2\right)}e^{-i(1/2+it)arctan\left(-(x-1/2\sqrt{2})/(x+1/2\sqrt{2})\right)}dx.$$

(A14)

One now requires some identities for the $arctan$ function. Adapting equations 1.627(2) and 1.622(3) from [16] gives, for any $\sigma\epsilon\mathbb{C}$

$$exp\left[\sigma arctan\left\{-\frac{\left(x-1/2\sqrt{2}\right)}{\left(x+1/2\sqrt{2}\right)}\right\}\right] = e^{-\frac{3\pi\sigma}{4}}exp\left[-\sigma arctan\left\{2\sqrt{2}x\right\}\right],$$

$$= e^{-\frac{3\pi\sigma}{4}}\exp\left[-\frac{\sigma}{2i}Log\left\{\frac{1+i2\sqrt{2}x}{1-i2\sqrt{2}x}\right\}\right], \qquad x < -\frac{1}{2\sqrt{2}}. \qquad \text{(A15)}$$

Hence the last two exponentials in the integrand (A14) can be written as

$$e^{-\frac{3\pi}{4}(t-i/2)}exp\left[(-1/4-it/2)log(1/8+x^2) + (1/4+it/2)Log\left\{\frac{1+i2\sqrt{2}x}{1-i2\sqrt{2}x}\right\}\right],$$

$$= e^{-\frac{3\pi}{4}(t-i/2)}exp\left[(-1/4-it/2)Log\left(1/2\sqrt{2}-ix\right)^2\right],$$

$$= e^{-\frac{3\pi}{4}(t-\mathrm{i}/2)} exp\left[(-1/4-\mathrm{i}t/2)Log\left\{(-\mathrm{i})^2\left(x+\mathrm{i}/2\sqrt{2}\right)^2\right\}\right],$$

$$= e^{-\frac{3\pi}{4}(t-\mathrm{i}/2)} exp\left[(-1/4-\mathrm{i}t/2)\left\{Log(-1)+Log\left(x+\mathrm{i}/2\sqrt{2}\right)^2\right\}\right],$$

$$= e^{-\frac{3\pi}{4}(t-\mathrm{i}/2)} e^{-\mathrm{i}\pi/4+\pi t/2} exp\left[(-1/4-\mathrm{i}t/2)\left\{Log\left(x+\mathrm{i}/2\sqrt{2}\right)^2\right\}\right]. \quad \text{(A16)}$$

Substituting (A16) into (A14) and making one further change of variables $v = x + \mathrm{i}/2\sqrt{2}$, yields, with some algebra

$$\mathrm{RSI}_- =$$

$$-\frac{e^{3\pi t/4}e^{\mathrm{i}3\pi/8}}{\mathrm{i}}\sum_{\substack{\alpha=1\\odd}}^{\infty}(-1)^{(\alpha-1)/2}\,e^{-\mathrm{i}\pi\alpha/2}\int_{-\infty+\mathrm{i}/2\sqrt{2}}^{-1/2\sqrt{2}+\mathrm{i}/2\sqrt{2}} e^{-\pi v^2}e^{\pi\alpha(1+\mathrm{i})v/\sqrt{2}}\,e^{(-1/4-\mathrm{i}t/2)Log(v^2)}dv.$$

(A17)

The integrals in (A17) are clearly of the form $Q$ as defined by equation (A7), with $\lambda \to \pi$, $\sigma \equiv \pi\alpha(1+\mathrm{i})/\sqrt{2}$ , $\beta \equiv -t/2 + \mathrm{i}/4$, and $w \equiv -1/2\sqrt{2}$. (One technical point is that the ubiquitous appearance of $\pi$ in all these formulae means that one has to separate it from the use of $\lambda$ as an independent variable in the development of (A17). Hence one chooses to take the limit $\lambda \to \pi$, rather than actually setting it equal to $\pi$.)

A2.2. *The $RSI_+$ along that part of the contour lying in the region with $Re(z) > 0$.*

Making the substitution $z = 1/2 + qe^{-\mathrm{i}\pi/4}$ this time for with $q \in \left(-1/\sqrt{2}, \infty\right)$ in (A1) allows one to develop an analogous expression to (A17) for the integral $\mathrm{RSI}_+$ along that part of the contour with the $Re(z) > 0$. The argument is very similar to that discussed in section 3.1 and only a few minor adjustments are necessary. Specifically the sum (A13) for the denominator changes as the contour crosses the real axis (when $q > 0$) to

$$\frac{1}{sin[\pi(1/2+qe^{-\mathrm{i}\pi/4})]} = 2\sum_{\substack{\alpha=1\\odd}}^{\infty}(-1)^{(\alpha-1)/2}e^{-\pi q(1+\mathrm{i})\alpha/\sqrt{2}}\ . \quad \text{(A18)}$$

In addition the corresponding *arctan* identity analogous to (A15) which is applicable for all $q > -1/\sqrt{2}$ (or $x > -1/2\sqrt{2}$ after the substitution $x = q + 1/2\sqrt{2}$), is given by

$$exp\left[\sigma arctan\left\{-\frac{\left(x-1/2\sqrt{2}\right)}{\left(x+1/2\sqrt{2}\right)}\right\}\right] = e^{\frac{\pi\sigma}{4}}exp\left[-\sigma arctan\left\{2\sqrt{2}x\right\}\right]. \quad \text{(A19)}$$

With these modifications one finds an equivalent sum to (A17) for $\mathrm{RSI}_+$, but in two parts namely

$\mathrm{RSI}_+ =$

$$-\frac{e^{3\pi t/4}e^{\mathrm{i}3\pi/8}}{\mathrm{i}}\sum_{\substack{\alpha=1\\ odd}}^{\infty}(-1)^{(\alpha-1)/2}\,e^{-\mathrm{i}\pi\alpha/2}\int_{-1/2\sqrt{2}+\mathrm{i}/2\sqrt{2}}^{1/2\sqrt{2}+\mathrm{i}/2\sqrt{2}} e^{-\pi v^2}e^{\pi\alpha(1+\mathrm{i})v/\sqrt{2}}\,e^{(-1/4-\mathrm{i}t/2)Log(v^2)}dv$$

$$-\frac{e^{3\pi t/4}e^{\mathrm{i}3\pi/8}}{\mathrm{i}}\sum_{\substack{\alpha=1\\ odd}}^{\infty}(-1)^{(\alpha-1)/2}\,e^{+\mathrm{i}\pi\alpha/2}\int_{1/2\sqrt{2}+\mathrm{i}/2\sqrt{2}}^{\infty+\mathrm{i}/2\sqrt{2}} e^{-\pi v^2}e^{-\pi\alpha(1+\mathrm{i})v/\sqrt{2}}\,e^{(-1/4-\mathrm{i}t/2)Log(v^2)}dv$$

(A20a, b)

Just like (A17), equation (A20) is essentially made up of two sums involving forms of the integral $Q$, with $\lambda \to \pi$, $\sigma \equiv \pm\pi\alpha(1+\mathrm{i})/\sqrt{2}$ , $\beta \equiv -t/2 + \mathrm{i}/4$ and minor modifications to the ranges of the integration.

Now at this stage it is perfectly possible to replace the integrals of type $Q$ in (A17) and (A20) by their asymptotic expansions in terms of reciprocal powers of the parameter $\beta$. This gives rise to geometric series in terms of the variable $r_- = e^{\pi(1+\mathrm{i})(Re(v)-1/2\sqrt{2})/\sqrt{2}}$ for (A17, A20a), and its reciprocal $r_+ = 1/r_-$ for (A20b). Here $Re(v)$ represents the real part of the cut-off of the integration range. These series can be summed to give asymptotic approximations (in powers of $\beta^{-1},\beta^{-2},\ldots$) for both the $\mathrm{RSI}_-$ and $\mathrm{RSI}_+$ which give good agreement with numerical calculations when the cut-offs are defined by $q \in (-\infty,\ -1/\sqrt{2})$ and $q \in (1/\sqrt{2},\infty)$, for $\mathrm{RSI}_-$ and $\mathrm{RSI}_+$ respectively. This is indicative that the methodology developed above to generate these approximations is quite correct over these particular ranges of $q$. However, over the range of $q \in (-1/\sqrt{2},\ 1/\sqrt{2})$ encapsulated by (A20a) agreement between the asymptotic approximations and the numerical calculations seriously breaks down. On extending the approximations for $\mathrm{RSI}_-$ and $\mathrm{RSI}_+$ so they now meet at the point $q = -1/\sqrt{2}$, one finds an exact cancellation between all the various terms of $O(\beta^{-1},\beta^{-2},\ldots)$‡, suggesting that over the entire integration range the $\mathrm{RSI} = 0$. This result is clearly erroneous, in light of both equation (2) and the numerical results shown in Table CII. So why does the methodology suddenly breakdown as one tries to match the approximations for the $\mathrm{RSI}_-$ and $\mathrm{RSI}_+$? The resolution of this anomaly lies in the interpretation of the expansion of the RSI in terms of a sum of integrals of the type $Q$, underlying equations (A17) and (A20).

A2.3. *The RSI written as a sum.*

The expansion (A13) for $\left[sin\left[\pi\left(1/2 + qe^{-\mathrm{i}\pi/4}\right)\right]\right]^{-1}$, used to develop sum (A17) for the $\mathrm{RSI}_-$, is only convergent for $q < 0$. Consequently (A13) has to be replaced by the expansion (A18) in order to develop the analogous sum (A20) for the $\mathrm{RSI}_+$. But having generated the sums (A17) and (A20), it then proves impossible to match the asymptotic expansions of the individual terms to give a suitable approximation for the RSI itself across the entire range of $q$, despite the fact that the methodology works well over the semi-infinite ranges mentioned above. However, there is an alternative way of linking the RSI to these sums, and that is to *start* with say (A13) as a (convergent, $q < 0$) geometric series of exponentials, and *only then* integrate across the entire range of $q$. Specifically, suppose that the $e^{\pi q(1+\mathrm{i})\alpha/\sqrt{2}}$ terms in (A13) are multiplied by $1/2\mathrm{i}$ ×the numerator of the integrand of equation (A12). When integrated over $q \in (-\infty,\ \infty)$, the resulting integrals are convergent for all $\alpha \in \mathbb{N}$ because of the presence of the $e^{-\pi q^2}$factor. Utilising the substitution $x = q + 1/2\sqrt{2}$, equations (A15) and (A19) for the arctan function and the substitution $v = x + \mathrm{i}/2\sqrt{2}$, one can heuristically interpret the RSI as the following sum (cf. A17 & A20)

$$\mathrm{RSI} = -\frac{e^{3\pi t/4} e^{\mathrm{i}3\pi/8}}{\mathrm{i}} \sum_{\substack{\alpha=1 \\ odd}}^{\infty} (-1)^{(\alpha-1)/2}\, e^{-\mathrm{i}\pi\alpha/2} \int_{-\infty+\mathrm{i}/2\sqrt{2}}^{\infty+\mathrm{i}/2\sqrt{2}} e^{-\pi v^2} e^{\pi\alpha(1+\mathrm{i})v/\sqrt{2}}\, e^{(-1/4-\mathrm{i}t/2)Log(v^2)} dv. \tag{A21}$$

Of course at this stage it is not entirely clear that this improvised continuation of the expansion (A13) to create a series of integrals evaluated over $\pm\infty$ is actually the correct interpretation of the RSI. The main reservation concerning (A21) is of course the convergence properties, or otherwise, of such a sum as $\alpha \to \infty$. This issue is clearly of crucial importance, but for the moment it will be put to one side, under the assumption that the terms in (A21) are sufficiently well behaved to ensure any convergence problems can be satisfactorily resolved. (The qualifications surrounding this assumption are discussed in sections C1.3 and A3.8. The proof that the series arising from this assumption does indeed yield an asymptotic expansion for the RS formula the subject of Section 2 and Appendix B). The upside of this interpretation is that the reason for the exact matching of the asymptotic approximations becomes apparent, but in turn this matching no longer leads to the erroneous conclusion that the $\mathrm{RSI} = 0$.

The integrals in (A21) are all of the type $Q$, but with $w \to \infty$ in the upper limit in (A7). Now as shown by (A10), such integrals can be written as sums of derivatives of the two foundation integrals $I_0$ and $I_1$. However, in the limit $w \to \infty$ the factor on the right-hand side of equation (A5) vanishes, which means (A6) reduces to

$$I_0 = C_0 \lambda^{-1/2-\mathrm{i}\beta} \qquad \text{and} \qquad I_1 = C_1 \lambda^{-1-\mathrm{i}\beta}. \tag{A22}$$

The disappearance of the first term on the right-hand side of (9) in this limit explains the matching of the $\mathrm{RSI}_-$ and the $\mathrm{RSI}_+$ expansions. All the terms in the various asymptotic approximations in powers of $\beta^{-1}, \beta^{-2}, \ldots$ were originally derived by assuming a finite value for $w$, so that the right-hand side of (A5) was non-zero. Piecing together the various terms in the asymptotic approximations derived at the different cut-off values is the equivalent of taking the limit $w \to \infty$ initially, and so it should it is of no surprise that the terms match up and cancel out to zero. However, if the asymptotic series in powers of $\beta^{-1}, \beta^{-2}, \ldots$ vanishes, the importance of the two extra terms (A22), which depend upon the unknown (and presumably non-zero) constants $C_0$ and $C_1$, is highlighted. These constants may, or may not, be sufficiently large to produce a significant amendment to the appropriate asymptotic approximation at any given finite cut-off value $w$ (as mentioned in section A2.2 numerical calculations suggest they are important when $-1/2\sqrt{2} < w < 3/2\sqrt{2}$), but in the limit $w \to \infty$ when they are the only terms left, their nature assumes a critical importance.

Substituting expressions (A22) for the foundation integrals into equation (A10) gives the result

$$\lim_{w\to\infty} Q = C_0 \sum_{n=0}^{\infty} \frac{\sigma^{2n}}{(2n)!}\left(-\frac{d}{d\lambda}\right)^n \lambda^{-1/2-\mathrm{i}\beta} + C_1 \sum_{n=0}^{\infty} \frac{\sigma^{2n+1}}{(2n+1)!}\left(-\frac{d}{d\lambda}\right)^n \lambda^{-1-\mathrm{i}\beta}$$
$$= C_0\lambda^{-1/2-\mathrm{i}\beta}\Phi\left(\frac{1}{2}+\mathrm{i}\beta, \frac{1}{2}; \frac{\sigma^2}{4\lambda}\right) + C_1\sigma\lambda^{-1-\mathrm{i}\beta}\Phi\left(1+\mathrm{i}\beta, \frac{3}{2}; \frac{\sigma^2}{4\lambda}\right), \qquad \text{(A23)}$$

where $\Phi(\alpha, \gamma; z)$ is the confluent hypergeometric function or Kummer's function ([1] Chapter 13). Taking the limit $\lambda \to \pi$ and making the identifications $\sigma \equiv \pi\alpha(1+\mathrm{i})/\sqrt{2}$ and $\beta \equiv -t/2 + \mathrm{i}/4$ appropriate for casting $Q$ into the form of the integrals in (A21), leads to

$$\text{RSI} = -\frac{e^{\frac{3\pi t}{4}} e^{\mathrm{i}3\pi/8}}{\mathrm{i}} \sum_{\substack{\alpha=1 \\ odd}}^{\infty} (-1)^{\frac{(\alpha-1)}{2}} e^{-\frac{\mathrm{i}\pi\alpha}{2}} \left[ C_0 \pi^{-\frac{1}{4}+\frac{\mathrm{i}t}{2}} \Phi\left(\frac{1}{4}-\frac{\mathrm{i}t}{2}, \frac{1}{2}; \frac{\pi\alpha^2 \mathrm{i}}{4}\right) + C_1 \frac{\pi\alpha(1+\mathrm{i})}{\sqrt{2}} \pi^{-\frac{3}{4}+\frac{\mathrm{i}t}{2}} \Phi\left(\frac{3}{4}-\frac{\mathrm{i}t}{2}, \frac{3}{2}; \frac{\pi\alpha^2 \mathrm{i}}{4}\right) \right].$$

(A24)

So the key to establishing the RSI as a sum now lies in the identification of the constants $C_0$ and $C_1$.

A2.4. *The constant terms* $C_0$ and $C_1$.

Consider the following generic integral denoted by $J_k$, defined by

$$J_k = \int_{-\infty}^{-\frac{1}{2\sqrt{2}}} x^k e^{-\lambda x^2} e^{(-1/4-\mathrm{i}t/2)Log\left((x+\mathrm{i}/2\sqrt{2})^2\right)} dx, \qquad \text{(A25)}$$

where $\lambda > 0$ as before. With $k = 0$, making the substitution $v = (x + \mathrm{i}/2\sqrt{2})$ gives the result

$$J_0 = e^{\frac{\lambda}{8}} \int_{-\infty+\mathrm{i}/2\sqrt{2}}^{-1/2\sqrt{2}+\mathrm{i}/2\sqrt{2}} e^{-\lambda v^2} e^{\lambda \mathrm{i} v/\sqrt{2}} e^{(-1/4-\mathrm{i}t/2)Log(v^2)} dv. \qquad \text{(A26)}$$

$J_0$ is clearly of the form $Q$, with $\sigma \equiv \mathrm{i}\lambda/\sqrt{2}$ and $\mathrm{i}\beta \equiv -1/4 - \mathrm{i}t/2$, the same factor multiplying the $Log(v^2)$ term in the integrals of (A21). With these identifications and utilising results (A11) and (A23) one can write $J_0$ in the form

$$J_0 = e^{\frac{\lambda}{8}} \left[ \sum_{n=0}^{\infty} \frac{(\mathrm{i}\lambda/\sqrt{2})^{2n}}{(2n)!}\left(-\frac{d}{d\lambda}\right)^n A_0(-1/2\sqrt{2}) + \sum_{n=0}^{\infty} \frac{(\mathrm{i}\lambda/\sqrt{2})^{2n+1}}{(2n+1)!}\left(-\frac{d}{d\lambda}\right)^n A_1(-1/2\sqrt{2}) \right.$$
$$\left. + C_0\lambda^{-1/4+\mathrm{i}t/2}\Phi\left(\frac{1}{4}-\frac{\mathrm{i}t}{2}, \frac{1}{2}; -\frac{\lambda}{8}\right) + C_1 \frac{\mathrm{i}\lambda}{\sqrt{2}} \lambda^{-3/4+\mathrm{i}t/2}\Phi\left(\frac{3}{4}-\frac{\mathrm{i}t}{2}, \frac{3}{2}; -\frac{\lambda}{8}\right) \right]. \qquad \text{(A27)}$$

Essentially (A27) is an expanded form of (A23) evaluated at $w = -1/2\sqrt{2}$, where the terms $A_0(-1/2\sqrt{2})$ and $A_1(-1/2\sqrt{2})$ are the appropriate integral terms derived from (A6), namely

$$A_{n=0,1} = \frac{\left(-1/2\sqrt{2} + \mathrm{i}/2\sqrt{2}\right)^{n+1}}{2\lambda^{(n+1)/2-1/4-\mathrm{i}t/2}} e^{(-1/4-\mathrm{i}t/2)Log\left(\left(-1/2\sqrt{2}+\mathrm{i}/2\sqrt{2}\right)^2\right)} \int \lambda^{(n-1)/2-1/4-\mathrm{i}t/2} e^{-\lambda\left(-1/2\sqrt{2}+\mathrm{i}/2\sqrt{2}\right)^2} d\lambda.$$

(A28)

Here the primitive for the indefinite integral in (A28) has no integration constant (because it has already been accounted for by $C_{0,1}$ in A27).

Alternatively by making the changes of variables $x = -u$, followed by $w = \lambda u^2$ in (A25), one can rewrite $J_0$ in the form

$$J_0 = \frac{1}{2}\lambda^{-1/4+\mathrm{i}t/2} \int_{\frac{\lambda}{8}}^{\infty} e^{-w} w^{-3/4-\mathrm{i}t/2} \left(1 - \frac{\mathrm{i}\sqrt{\lambda}}{2\sqrt{2w}}\right)^{2(-1/4-\mathrm{i}t/2)} dw. \qquad \text{(A29)}$$

Now equating equations (A27) and (A29) and multiplying both sides by the factor $\lambda^{1/4-\mathrm{i}t/2}$ gives

$$\frac{1}{2}\int_{\frac{\lambda}{8}}^{\infty} e^{-w} w^{-3/4-\mathrm{i}t/2} \left(1 - \frac{\mathrm{i}\sqrt{\lambda}}{2\sqrt{2w}}\right)^{2(-1/4-\mathrm{i}t/2)} dw$$
$$= e^{\frac{\lambda}{8}}\Bigg[\lambda^{1/4-\mathrm{i}t/2}\sum_{n=0}^{\infty}\frac{\left(\mathrm{i}\lambda/\sqrt{2}\right)^{2n}}{(2n)!}\left(-\frac{d}{d\lambda}\right)^n A_0(-1/2\sqrt{2}) + \lambda^{1/4-\mathrm{i}t/2}\sum_{n=0}^{\infty}\frac{\left(\mathrm{i}\lambda/\sqrt{2}\right)^{2n+1}}{(2n+1)!}\left(-\frac{d}{d\lambda}\right)^n A_1(-1/2\sqrt{2})$$
$$+ C_0\Phi\left(\frac{1}{4} - \frac{it}{2}, \frac{1}{2}; -\frac{\lambda}{8}\right) + C_1\frac{\mathrm{i}\lambda}{\sqrt{2}}\lambda^{-1/2}\Phi\left(\frac{3}{4} - \frac{it}{2}, \frac{3}{2}; -\frac{\lambda}{8}\right)\Bigg].$$

(A30)

Now consider what happens in the limit $\lambda \to 0$. On the right hand side of (A30) the first and second terms are of $O(\lambda^{1/4})$ and $O(\lambda^{5/4})$ respectively. Hence in this limit they both vanish. (The terms involving the smallest powers of $\lambda$ are found when $n = 0$ in (A30), in which case the integrand in (A28) can expanded as a sum and integrated term by term to show that both $A_0$ and $A_1$are of $O(1)$ for small $\lambda$.) The fourth term also vanishes as $\Phi(\alpha, \gamma; 0) = 1$. Consequently the right-hand side of (A30) reduces to just $C_0 \times 1$. The left-hand side of (A30) reduces to a gamma function integral, and so for equality to hold one must have

$$C_0 = \frac{1}{2}\Gamma\left(\frac{1}{4} - \frac{\mathrm{i}t}{2}\right). \qquad \text{(A31)}$$

A very similar argument, but this time using the integral $J_1$ as a starting point, can be used to deduce $C_1$. Analogous expressions to (A27) and (A29) are first derived and equated. This time though the result analogous to (A27) is derived from the derivative of the expansion for $Q$ (equation A11) with respect to $\sigma$. (Again this assumes that such a differentiation under the integral in equation A7 is a valid operation.) Finally the two equated expressions for $J_1$ are examined in the limit $\lambda \to 0$. The details are omitted, but ultimately one finds that

$$C_1 = -\frac{1}{2}\Gamma\left(\frac{3}{4}-\frac{\mathrm{i}t}{2}\right) = -\frac{1}{2}\frac{\pi\sqrt{2}}{cosh(\pi t/2)+\mathrm{i}sinh(\pi t/2)}\frac{1}{\Gamma(1/4+\mathrm{i}\,t/2)}. \qquad \text{(A32)}$$

The latter relation ([1] eq. 6.1.32) is useful because $\Gamma(1/4+\mathrm{i}\,t/2)$ plays a fundamental rôle in the expansion of the zeta function along the critical line.

When these identities (A31-32) for the constants $C_0$ and $C_1$ are substituted into (A24), then the interpretation of the RSI as an infinite sum of the form (A21), leads one to conclude that

$$\mathrm{RSI} = -\frac{e^{3\pi t/4}e^{\mathrm{i}3\pi/8}}{\mathrm{i}}\sum_{\substack{\alpha=1\\ odd}}^{\infty}(-1)^{(\alpha-1)/2}\,e^{-\mathrm{i}\pi\alpha/2}\,\frac{1}{2}\left[\Gamma\left(\frac{1}{4}-\frac{\mathrm{i}t}{2}\right)\pi^{-\frac{1}{4}+\frac{\mathrm{i}t}{2}}\Phi\left(\frac{1}{4}-\frac{\mathrm{i}t}{2},\frac{1}{2};\frac{\pi\alpha^2\mathrm{i}}{4}\right)\right.$$

$$\left.-\frac{\pi\sqrt{2}}{cosh(\pi t/2)+\mathrm{i}sinh(\pi t/2)}\frac{1}{\Gamma(1/4+\mathrm{i}\,t/2)}\frac{\pi\alpha(1+\mathrm{i})}{\sqrt{2}}\pi^{-\frac{3}{4}+\frac{\mathrm{i}t}{2}}\Phi\left(\frac{3}{4}-\frac{\mathrm{i}t}{2},\frac{3}{2};\frac{\pi\alpha^2\mathrm{i}}{4}\right)\right],$$

(A33)

provided that (A33) satisfies suitable convergence criteria. To make any further progress in this regard one needs to examine the properties of the gamma and confluent hypergeometric functions appearing in (A33).

## A3. The Real and Imaginary Parts of the RS Integral

A3.1 *Asymptotic expansions for the gamma functions* $\Gamma(1/4 \pm \mathrm{i}\,t/2)$.

The gamma function can be written as

$$\Gamma(1/4+\mathrm{i}\,t/2) = exp\left[Re\left(Log\Gamma(1/4+\mathrm{i}\,t/2)\right)+\mathrm{i}Im\left(Log\Gamma(1/4+\mathrm{i}\,t/2)\right)\right]$$
$$= e^{Re\left(Log\Gamma(1/4+\mathrm{i}t/2)\right)}e^{\mathrm{i}\left(\theta(t)+(t/2)log(\pi)\right)} = e^{Re\left(Log\Gamma(1/4+\mathrm{i}t/2)\right)}e^{\mathrm{i}\theta(t)}\pi^{\mathrm{i}t/2}. \qquad \text{(A34)}$$

Here $\theta(t)$ is the Riemann-Siegel theta function which has an asymptotic expansion given by $\theta(t) = t(log(t/2\pi)-1)/2-\pi/8+1/48t+7/5760t^3\ldots$ (see [13] Chapter 6). Using Stirling's asymptotic expansion for the gamma function ([1] eq. 6.1.41) one finds that for the real part

$$Re\left(Log\Gamma(1/4+\mathrm{i}\,t/2)\right) = -\frac{1}{4}log\left(\frac{t}{2}\right)-\frac{\pi t}{4}+\frac{log(2\pi)}{2}+\frac{1}{32t^2}+\frac{5}{256t^4}+\frac{61}{1536t^6}\cdots \qquad \text{(A35)}$$

$$\Rightarrow e^{Re\left(Log\Gamma(1/4+\mathrm{i}t/2)\right)} = \left(\frac{t}{2}\right)^{-1/4}e^{-\pi t/4}\sqrt{2\pi}e^{1/32t^2+\varpi(t)}, \qquad \text{(A36)}$$

where $\varpi(t) = 5/256t^4+61/1536t^6$ ... One further identity that proves useful is

$$\frac{1}{cosh(\pi t/2)+\mathrm{i}sinh(\pi t/2)} = \frac{(1+e^{-\pi t})+\mathrm{i}(e^{-\pi t}-1)}{e^{\pi t/2}(1+e^{-2\pi t})}. \qquad \text{(A37)}$$

Substituting results (A34, 36-37) and the corresponding formulae for $\Gamma(1/4 - \mathrm{i}\,t/2)$, into (A33) gives

$$\mathrm{RSI} = -\frac{e^{\frac{\pi t}{2}+\mathrm{i}3\pi/8-\mathrm{i}\theta(t)}}{2\mathrm{i}(1+e^{-2\pi t})}\sum_{\substack{\alpha=1\\ odd}}^{\infty}(-1)^{(\alpha-1)/2}\,e^{-\mathrm{i}\pi\alpha/2}\left[\left(\frac{t\pi}{2}\right)^{-\frac{1}{4}}\sqrt{2\pi}\,e^{[1/32t^2+\varpi]}(1+e^{-2\pi t})\Phi\left(\frac{1}{4}-\frac{\mathrm{i}t}{2},\frac{1}{2};\frac{\pi\alpha^2\mathrm{i}}{4}\right)\right.$$

$$\left.-\frac{2\alpha}{\sqrt{2\pi}}\left(\frac{t\pi^5}{2}\right)^{1/4}e^{-[1/32t^2+\varpi]}(1+\mathrm{i}e^{-\pi t})\Phi\left(\frac{3}{4}-\frac{\mathrm{i}t}{2},\frac{3}{2};\frac{\pi\alpha^2\mathrm{i}}{4}\right)\right]$$

(A38)

From (A38), one can see that the proposed expansion (A33) for the RSI is made up of four terms, of significantly different orders of magnitude. There are two large terms, $O(t^{-1/4})\times\Phi$ and $O(\alpha t^{1/4})\times\Phi$, which turn out to be of the same order of magnitude (see Appendix C, eq. C6), followed by a term which is smaller by a factor of $\mathrm{i}e^{-\pi t}$, and then a further term smaller by a factor $(e^{-\pi t})^2$. One would expect (A38) to be dominated by the sum of the difference between the two largest terms, which is the subject of Appendix C. It turns out that not only does this sum converge absolutely, but the resulting asymptotic expansion (C11) agrees very well with numerical estimates of the RSI, indicating that the postulated sum (A33) is indeed the correct interpretation of the integral. (There is a qualification to these results concerning the requirement that terms of $O(e^{\alpha\sqrt{\pi t}/\sqrt{2}})$ in the expansion should cancel for all orders of reciprocal powers of $t$, see Sections C1.3 and A3.8) However, this analysis also shows that (C11) provides only an estimate of the $Im[e^{\mathrm{i}\theta(t)}\mathrm{RSI}(1/2+\mathrm{i}t)]$. So despite forming the overwhelmingly dominant part of (A38), the expansion (C11) is unrelated to the zeta function because, on combining results (1-3) and using (A34), one finds that

$$\zeta(1/2+\mathrm{i}t) = e^{-\mathrm{i}\theta(t)}2Re\left[e^{\mathrm{i}\theta(t)}\mathrm{RSI}(1/2+\mathrm{i}t)\right] = e^{-\mathrm{i}\theta(t)}Z(t), \qquad \text{(A39)}$$

($Z(t)$ being the Hardy function, see [5] & [13]). Consequently in order to establish connection between the Hardy function and (A38), it is the next largest term, that involving $\mathrm{i}e^{-\pi t}$, which now becomes the main focus of this investigation.

### A3.2 *Integral representations of the confluent hypergeometric function*

The term of second order in powers of $e^{-\pi t}$ in the expansion for the RSI in (A38) is given by

$$\mathrm{RSI}_{2^{nd}\,order} = \frac{e^{-\frac{\pi t}{2}+\mathrm{i}3\pi/8-\mathrm{i}\theta(t)}}{(1+e^{-2\pi t})\sqrt{2\pi}}\left(\frac{t\pi^5}{2}\right)^{1/4}e^{-[1/32t^2+\varpi]}\sum_{\substack{\alpha=1\\ odd}}^{\infty}(-1)^{(\alpha-1)/2}\,e^{-\mathrm{i}\pi\alpha/2}\,\alpha\Phi\left(\frac{3}{4}-\frac{\mathrm{i}t}{2},\frac{3}{2};\frac{\pi\alpha^2\mathrm{i}}{4}\right). \qquad \text{(A40)}$$

The methodology employed in Appendix C to represent $\Phi$ and estimate the first order sum is of little use for calculating (A40), because the terms of $O(e^{\alpha\sqrt{\pi t}/\sqrt{2}})$ that appear in the resulting expansions cannot be cancelled out. Rather one is left with an extremely complicated alternating series in powers of $\alpha$ and reciprocal powers of $t$, for which there seems little prospect of establishing any sum-able simplification. Nor do the advanced computational methods of [29] for direct evaluation of the confluent hypergeometric function numerically seem to offer much prospect of estimating (A40) for large $t$. However, there are some interesting features one can discern when evaluating (A40) numerically for very modest $t$ (in the range 10 to 100). First there is the *extremely small* contribution of the first few terms when $\alpha$ is small, which is followed by a sudden increase in the size of the terms, before taking on an oscillatory character, with amplitudes decaying at a rate very close to $\alpha^{-1/2}$ characteristic of the Dirichlet series along the critical line. The key to understanding this curious behaviour lies in the evaluation of a particular integral representation of the confluent hypergeometic function.

First one can rewrite $\Phi(3/4-\mathrm{i}\,t/2\,,3/2\,;\pi\alpha^2\mathrm{i}/4)$ in terms of Euler's integral formulation ([27], eq. 4.2(1)), viz.

$$\Phi\left(\frac{3}{4}-\frac{\mathrm{i}t}{2},\frac{3}{2};\frac{\pi\alpha^2\mathrm{i}}{4}\right)=\frac{\Gamma(3/2)}{\Gamma(3/4+\mathrm{i}t/2)\Gamma(3/4-\mathrm{i}t/2)}\int_0^1 e^{\mathrm{i}\pi\alpha^2x/4}\frac{\exp[\mathrm{i}(t/2)log\{(1-x)/x\}]}{[x(1-x)]^{1/4}}dx\,,$$

$$=\frac{\pi^{-1/2}e^{\pi t/2}e^{2/(12t^2+3)+O(t^{-4})}}{4(t^2/4+1/16)^{1/4}e^{1/24t^2+O(t^{-4})}}\int_0^1 e^{\mathrm{i}\pi\alpha^2x/4}\frac{\exp[\mathrm{i}(t/2)log\{(1-x)/x\}]}{[x(1-x)]^{1/4}}dx\,,\qquad\text{(A41)}$$

using Stirling's asymptotic expansion for the gamma functions. Substituting (A41) into (A40) gives

$$\mathrm{RSI}_{2^{nd}\,order}=\frac{(\pi/2)^{1/4}e^{\mathrm{i}3\pi/8-\mathrm{i}\theta(t)}\mathcal{H}(t)}{4t^{1/4}}\sum_{\substack{\alpha=1\\odd}}^{\infty}(-1)^{(\alpha-1)/2}\,e^{-\mathrm{i}\pi\alpha/2}\,\alpha\int_0^1 e^{\mathrm{i}\pi\alpha^2x/4}\frac{\exp[\mathrm{i}(t/2)log\{(1-x)/x\}]}{[x(1-x)]^{1/4}}dx,$$

where $$\mathcal{H}(t)=\frac{e^{-[1/32t^2+O(t^{-4})]}e^{2/(12t^2+3)+O(t^{-4})}}{(1+e^{-2\pi t})(1+1/4t^2)^{1/4}e^{[1/24t^2+O(t^{-4})]}}\approx 1+\frac{1}{32t^2}.\qquad\text{(A42a, b)}$$

With its singularities and undefined *log* term at the endpoints, this integral form does not appear a particularly promising means for estimation. However, by first splitting the range of integration in two, and then applying the substitution $y=1-x$ to the top half one finds

$$\int_0^1 e^{\mathrm{i}\pi\alpha^2x/4}\frac{\exp[\mathrm{i}(t/2)log\{(1-x)/x\}]}{[x(1-x)]^{1/4}}dx=\int_0^{1/2} e^{\mathrm{i}\pi\alpha^2x/4}\frac{\exp[\mathrm{i}(t/2)log\{(1-x)/x\}]}{[x(1-x)]^{1/4}}dx$$

$$+e^{\mathrm{i}\pi\alpha^2/4}\int_0^{1/2} e^{-\mathrm{i}\pi\alpha^2y/4}\frac{\exp[-\mathrm{i}(t/2)log\{(1-y)/y\}]}{[y(1-y)]^{1/4}}dy\,,\ \text{(A43)}$$

where the latter integral is clearly just the complex conjugate of the former. Secondly the substitution $x=1/(w+1)$ into the former integral yields

$$B=\int_0^{1/2} e^{\mathrm{i}\pi\alpha^2x/4}\frac{\exp[\mathrm{i}(t/2)log\{(1-x)/x\}]}{[x(1-x)]^{1/4}}dx=\int_1^{\infty}\frac{\exp[\mathrm{i}\pi\alpha^2/4(w+1)+\mathrm{i}(t/2)log(w)]}{w^{1/4}(w+1)^{3/2}}dw.\qquad\text{(A44)}$$

Integral $B$ looks much more promising, as the integrand is clearly defined at the endpoints and tends to zero as $w\to\infty$. It turns out that utilising a suitable closed contour, which avoids the singularities are at $w = 0$ and $w = -1$, and applying Cauchy's theorem one can devise very accurate asymptotic estimates for $B$ which prove suitable for the calculation of (A42a).

A3.3 *Asymptotic Analysis of the integral B – The exponential cut-off at* $\alpha = (8t/\pi)^{1/2}$

This analysis begins by treating $B$ defined above as contour integral with $w \in \mathbb{C}$. The first important observation is that on the circle $w = Re^{i\phi}$ with $R \gg \alpha^2$, the exponential term in the integrand will be dominated by the second term $\mathrm{i}(t/2)Log(w) = \mathrm{i}(t/2)log(R) - (t\phi/2)$. So provided $\phi > 0$, the modulus of the integrand will be $O\big(e^{-t\phi/2}/R^{7/4}\big)$ and the contribution to the integral along such a circular contour as $R \to \infty$ will tend to zero as $R^{-3/4}$. Notice too that the integrand becomes exponentially small as $|w|$ gets large, provided $Arg(w) > 0$. This means that the integrand along any contour which approaches the circle $w = Re^{i\phi}$ from a direction defined by $Arg(w) > 0$ will be exponentially decreasing along its entire length.

The next observation concerns the behaviour of the integral along the unit circle $w = e^{\mathrm{i}\phi}$. Consider integrating around the unit circle in both the anti and clockwise senses, starting from $w = 1$ so that $\phi \in [0, \pm\varphi]$ where $\varphi \le \pi/2$. This transforms $B$ into the form

$$B = \mathrm{i}e^{\mathrm{i}\pi\alpha^2/8}\int_0^{\pm\varphi} \frac{\exp\big[\alpha^2\pi sin(\phi)/8\big(1+cos(\phi)\big) - t\phi/2\big]e^{\mathrm{i}\phi}}{e^{\mathrm{i}\phi/4}(1+e^{\mathrm{i}\phi})^{3/2}}\,d\phi. \qquad \text{(A45)}$$

Now for small $\phi$ the exponent of the exponential term behaves $\approx [\alpha^2\pi/8 - t]\phi/2$ . That means that if $\alpha < (8t/\pi)^{1/2}$ and $\phi > 0$ (moving anticlockwise) then the integrand decays exponentially in this direction. Similarly the integrand decays exponentially when moving clockwise, $\phi < 0$, with $\alpha > (8t/\pi)^{1/2}$. The parameter $(8t/\pi)^{1/2}$ turns out to be extremely important, and will be denoted subsequently by the letter $a$. Making a further substitution $x = sin(\phi)/\big(1 + cos(\phi)\big)$ so that $x \in [0, \pm p]$ where $p \le 1$, one can transform (A45) into

$$\mathrm{i}e^{\mathrm{i}\pi\alpha^2/8}\int_0^{\pm p} \frac{\exp[(\mathrm{i}\,3/4 - t/2)arccos(\{1-x^2\}/\{1+x^2\})]e^{\alpha^2\pi x/8}}{2^{3/4}e^{\pm\mathrm{i}(3/2)arctan(x)}}\left[\frac{2}{1+x^2}\right]^{1/4} dx$$

$$= \mathrm{i}e^{\mathrm{i}\pi\alpha^2/8}\int_0^{\pm p} \frac{\exp[(-t/2)arccos(\{1-x^2\}/\{1+x^2\})]e^{\alpha^2\pi x/8}}{2^{3/4}}\left[\frac{2}{1+x^2}\right]^{1/4} dx, \qquad \text{(A46)}$$

which holds because $arccos(\{1-x^2\}/\{1+x^2\}) = \pm 2arctan(x)$ for $x \gtreqless 0$. Consequently the integral in (A46) is *real*. This observation means that the contribution of $B$ to the integral (A43), along a contour confined to $w = e^{\mathrm{i}\phi}$, must be zero, simply because $\mathrm{i}e^{\mathrm{i}\pi\alpha^2/8}real + e^{\mathrm{i}\pi\alpha^2/4} \times \big[-\mathrm{i}e^{-\mathrm{i}\pi\alpha^2/8}real\big] = 0$. Consequently $B$ makes no contribution to the sum (A42a) when the integration takes place along the unit circle to any point $\phi \in [-\pi/2, \pi/2]$. Now when $\alpha < a$, one can imagine integrating anticlockwise around the unit circle to some 'suitable point' (for a specific definition see section A3.6) where the integrand is exponentially small (an operation that contributes nothing to A43), and then integrating down a contour towards the circle $w = Re^{i\phi}$ with $\phi > 0$, along which the integrand *starts* and *remains* exponentially small (actually decreases) over its entire length (see Fig. A1a). In such circumstances the contribution to integral (A43) is dominated by the size of the integral $B$ close to the 'suitable point' (the first point on the contour which actually contributes to A43). Consequently the contribution of *all* the terms satisfying $\alpha < a$ to the sum (A42a) will be *exponentially small*. This exponential cut-off at $\alpha = a$ explains why the initial few terms to (A43) were observed to contribute almost nothing when the hypergeometric function in (A40) was calculated numerically. Hence the sum (A42a) will be dominated by those terms for which $\alpha > a$, and their estimation is the subject of the next two sections.

A3.4 *Estimating integral B when* $\alpha > a$.

In principle one can obtain very good estimates for (A44) by means of contour integration, provided the contour passes through certain saddle points where the modulus of the integrand is very large compared to points on the contour in the immediate neighbourhood of the saddle. Fig. A2a illustrates just such a suitable contour for the case $\alpha > a$. First define 'zero-levels' as those points $w \in \mathbb{C}$ such that the $Re[\mathrm{i}\pi\alpha^2/4(w+1) + \mathrm{i}(t/2)Log(w)] = 0$, where the exponential term in the integrand of $B$ has modulus of unity. The guiding principle used to establish the contour shown in Fig. A2a is move through regions in the complex plane at or below 'zero-level', starting from the point $w = 1$ (avoiding the singularities at 0 and $-1$) and terminating with $|w| \longrightarrow \infty$. Then the various contributions of the integral along each part of the path will be heavily concentrated at the 'zero-level' points through which the contour passes. These contributions can then be estimated using standard local asymptotic analysis, such as the saddle point or Laplace method (e.g. [2] or [8]).

It turns out that when $\alpha > a$ the integrand of $B$ has a 'zero-level' saddle point, situated on the positive real axis at the point $w = pc$, defined by

$$w = pc(\alpha) = (\pi\alpha^2/4t) - 1 + (\pi\alpha^2/4t)(1 - a^2/\alpha^2)^{1/2},$$

$$\Rightarrow \left(\frac{\pi\alpha^2}{4}\right) = \frac{t(pc+1)^2}{2pc}, \qquad \text{(A47)}$$

where $c = \alpha^2\big(1 + (1 - 4t^2/\alpha^4)^{1/2}\big)/2t$ and $p = \big[(\pi\alpha^2/2)\big(1 + (1 - 8t/\pi\alpha^2)^{1/2}\big) - 2t\big]/2tc$. It will transpire that the variable $pc(\alpha)$ plays a fundamental rôle in connecting the traditional representation of the zeta function and that to be derived from the RSI. For the moment it is sufficient to note that $pc > 1$ if $\alpha > a$, and it possess an asymptotic expansion $pc \approx \pi\alpha^2/2t - 2$ for large $\alpha$ ($p \to \pi/2 - 2t/\alpha^2$). One can exploit the existence of these saddle points at (A47) to construct a contour path suitable for estimating integral $B$.

Beginning at $w = 1$, move clockwise around the unit circle to the point $w = -\mathrm{i}$. This defines the initial branch of the contour and is entirely made up of below 'zero-level' points. At $w = -\mathrm{i}$, change direction and move along a second branch of the contour defined by $w = -\mathrm{i} + qe^{-\mathrm{i}\pi/4}$, where $q \geq 0$, as shown in Fig. A2a. Along this branch of the contour the value of $Re[\mathrm{i}\pi\alpha^2/4(w+1) + \mathrm{i}(t/2)Log(w)]$ increases, but before reaching zero it intersects at right angles a third branch of the contour, which emanates from $w = pc$ at an angle to the horizontal of $\pi/4$ radians. This branch is defined by $w = pc + ue^{i\pi/4}$, and the intersection occurs at $u = u_1 = -(pc+1)/\sqrt{2}$. Moving away from $w = pc$ along that part of the third branch with $u < 0$, one finds that the $Re[i\pi\alpha^2/4(w+1) + i(t/2)Log(w)]$ initially decreases to a minimum, before increasing again. The position of this minimum is found to be situated approximately at $u \approx -(pc+1)/\sqrt{2} - \sqrt{2}/4pc$, i.e. roughly where the second and third branches of the contour intersect (and coalescing as $\alpha$ and $pc$ get large). From the intersection point one moves along the third branch of the contour by increasing $u$, passing through the saddle point at $u = 0$ on to a fourth branch, and then continuing on by letting $u \to \infty$ to eventually reach the circumference of a circle of arbitrarily large radius. Along the fourth branch of the contour with $u > 0$, one finds that the $Re[i\pi\alpha^2/4(w+1) + i(t/2)Log(w)]$ continuously decreases, as explained at the start of section A3.3.

Integrand (A44) is analytic on all points on and within the region enclosed by the contour and so from Cauchy's theorem $B = I_1 + I_2 + I_3 + I_4 + I_5$ , where $I_n$ defines the integrals along the various stages of the contour as shown on Fig. A2a. Estimates of each will provide an estimate for $B$. However, three can be disposed of easily. Integral $I_5$ clearly tends to zero as the radius of the sector tends to infinity. The integral $I_2$ is exponentially small, even if $pc(\alpha)$ and hence the length of this branch of the contour tends to infinity. (A rough estimate gives $|I_2| < 2\sqrt{2}e^{t[\pi/8-(pc+1)/2pc]}pc/t(pc+1)$.) And whilst integral $I_1$ does contribute to $B$ (actually a term of $O(t^{-1})$), because the integration takes place along the unit circle it contributes nothing to Euler's integral formulation (A43), for the reasons discussed in section A3.3. The remaining two integrals require more detailed analysis near the saddle point. If $1 < pc \le \sqrt{2}$, the contour requires a certain degree of modification, because in this case the third branch intersects the unit circle at $u = u_2 = \sqrt{2}[cos\{arccos(pc/\sqrt{2}) - \pi/4\} - pc]$ before reaching the second branch at $u = u_1$ (see Fig. A2b). However, this is only a minor inconvenience as one can simply integrate around the unit circle to this first intersection (below 'zero level') and then move along the contour through the saddle, eliminating the negligibly small contribution of integral $I_2$ altogether.

A3.5 *The Saddle Points of Integral* $B$ *when* $\alpha > a$.

Let the branches of the contour passing through the point $w = pc$ be defined as the line $w = pc + ue^{iq}$, where $u$ is real. Expanding the $Re[i\pi\alpha^2/4(w+1) + i(t/2)Log(w)]$ as a power series about $u = 0$, yields

$$\frac{(\pi\alpha^2/4)usin(q)}{[(pc+1)^2 + u^2 + 2u(pc+1)cos(q)]} - \frac{t}{2}arctan\left(\frac{usin(q)}{pc + ucos(q)}\right)$$
$$= sin(q)\left[\frac{(\pi\alpha^2/4)}{(pc+1)^2} - \frac{t}{2pc}\right]u - \frac{tsin(q)cos(q)(pc-1)}{2(pc)^2(pc+1)}u^2 + \sum_{i=3}^{\infty} d_i(pc)u^i \,. \text{(A48)}$$

From the definition of $pc$ given by (A47), one can see that the first coefficient is zero and that $w = pc$ is indeed a saddle point. What is more when $\alpha > a, pc > 1$, and the coefficient multiplying $u^2$ will be negative provided $q$ is an angle in the first quadrant. The choice that maximises this coefficient and the line of steepest descent is clearly $q = \pi/4$ and hence the reason for the approach angle shown in Fig. 2a. With this choice, $d_2 = -\frac{t(pc-1)}{4(pc)^2(pc+1)} < 0$, whilst the other coefficients are given by $d_{4i} = 0$, $d_i = \frac{(-1)^{(i+1)/2}t\left(1/(1+pc^{-1})^{i-1}-1/i\right)}{2\sqrt{2}(pc)^i}$ for $i$ odd and $d_{4i+2} = \frac{(-1)^{i+1}t\left(1/(1+pc^{-1})^{4i+1}-1/(4i+2)\right)}{2(pc)^{4i+2}}$. Concentrating on the integral $I_4$ along the branch of the contour with $u > 0$, one finds that

$$I_4 = \frac{(1+\mathrm{i})}{\sqrt{2}}\int_0^\infty e^{d_2u^2}f(u)du = \frac{(1+\mathrm{i})}{\sqrt{2}}\int_0^\infty \frac{e^{d_2u^2}h(u)}{[(pc)^2+u^2+u\sqrt{2}pc]^{1/8}[(pc+1)^2+u^2+u\sqrt{2}(pc+1)]^{3/4}}du, \qquad \text{(A49)}$$

where

$$h(u) = exp\left[\sum_{i=3}^{\infty} d_i(pc)u^i + \mathrm{i}\left\{g(u) - \frac{arctan\left(u/(pc\sqrt{2}+u)\right)}{4} - \frac{3arctan\left(u/(pc\sqrt{2}+\sqrt{2}+u)\right)}{2}\right\}\right],$$

and

$$g(u) = \frac{(\pi\alpha^2/4)(pc+1+u/\sqrt{2})}{[(pc+1)^2+u^2+u\sqrt{2}(pc+1)]} + \frac{tlog((pc)^2+u^2+u\sqrt{2}pc)}{4}. \qquad \text{(A50a, b)}$$

Hence

$$I_4 = \frac{(1+\mathrm{i})\sqrt{\pi}}{\sqrt{2}}\sum_{n=0}^{\infty}\frac{(2n-1)!!\,f^{(2n)}(0)}{(2n)!\,2^{n+1}|d_2|^{n+1/2}} + \frac{(1+\mathrm{i})}{\sqrt{2}}\sum_{n=0}^{\infty}\frac{n!\,f^{(2n+1)}(0)}{(2n+1)!\,2|d_2|^{n+1}}\,. \qquad \text{(A51)}$$

Now the character of the series (A51) depends upon the relative magnitudes of $f^{(2n)}(0)$ to $|d_2|^{n+1/2}$. Now if $(pc-1)\sim O(1)$, then $|f(0)|/|d_2|^{1/2} \sim O\left(t^{-1/2}\right)$. Now in this $pc$ regime, the dominant term of (A50a) is the $O(t)$ function $g(u)$ given by (A50b). Hence to first order one would expect

$$f'(0) \sim \frac{g'(0)h(0)}{pc^{1/4}(pc+1)^{3/2}} = O(t) \quad \text{and} \quad f''(0) \sim \frac{[g''(0)+(g'(0))^2]h(0)}{pc^{1/4}(pc+1)^{3/2}} = O(t^2). \qquad \text{(A52)}$$

However, *both* $g'(0) = 0$ (as a consequence of $pc$ satisfying A47) and $g''(0) = 0$ (by differentiation), which means that both $f'(0)$ and $f''(0)$ are in fact terms of $O(1)$, just like $f(0)$ itself. The derivatives of the function $f$ only reach a size of $O(t)$ at its third derivative, which means that the series (A51) starts with a term $O\left(t^{-1/2}\right)$ and is then followed by terms no larger than $O(t^{-1})$. Likewise in the regime when $\alpha$ is large, so that $pc \sim O(\pi\alpha^2/2t)$, the first term in the series $f(0)/|d_2|^{1/2}$ is of $O\left(t^{1/4}\alpha^{-3/2}\right)$, whilst the next largest terms (actually the second and fourth) are no more than $O\left(t^{-1/4}\alpha^{-3/2}\right)$. Consequently series (A51) for both these regimes is completely dominated by its initial term and the integral $I_4$ can be written as

$$I_4 = \frac{(1+\mathrm{i})\sqrt{\pi}}{\sqrt{2}}\frac{(pc)^{3/4}exp[\mathrm{i}\{\pi\alpha^2/4(pc+1) + tlog(pc)/2\}]}{(pc+1)\sqrt{t(pc-1)}}\left[1 + O\left(\frac{1}{\sqrt{t}}\right)\right]. \qquad \text{(A53)}$$

Obviously (A53) cannot hold as $pc \to 1$ because $|d_2|$ will no longer be of $O(t)$. To see what happens in this limit let $\alpha = a + \varepsilon$, which implies that $pc \approx 1 + \sqrt{2\varepsilon}(2\pi/t)^{1/4}$. This means that $|d_2| \sim O\left(\sqrt{\varepsilon}t^{\frac{3}{4}}\right)$ and the first term in the series will be of $O(\varepsilon^{-1/4}t^{-3/8})$. The next largest term in series (A51) will be the fourth, when $f'''(0) \sim O(t)$ and $f'''(0)/|d_2|^2 \sim O(\varepsilon^{-1}t^{-1/2})$. Letting $\varepsilon = t^{-s}$ one finds that the first term in (A51) ceases to dominate the series when

$$s/4 - 3/8 = s - 1/2\,, \Rightarrow s = 1/6. \quad \text{(A54)}$$

So provided $\alpha > a + t^{-1/6}$, the initial term in the series will continue to dominate the others and (A53) will provide an excellent approximation for the integral emanating from the saddle point. (More detailed calculations show that the initial and fourth terms of equation A51 become equal in magnitude when $\alpha \approx a + 0.086t^{-1/6}$, $\Rightarrow pc \approx 1 + 0.657t^{-1/3}$. When $\alpha = a + t^{-1/6}$ the initial term dominates by a factor approximately equal to $2\pi$.) So equation (A53) should provide an excellent approximation for the integral $I_4$ for all the odd integer values that lie above $a$, *except possibly the first*. And even then, this first odd integer must lie *extremely* close to $a$ before a serious problem arises with (A53). This important transitional region $[a, a + t^{-1/6}]$ will re-occur again in the analysis of Appendix B.

For the integral $I_3$ defined along the branch of the contour emanating from the saddle point with $u < 0$, the analysis is virtually identical. (There is now a cut-off point at either $u = u_{1,2}$ rather than the semi-infinite contour defined above, but unless $pc \approx 1$, the integrand at these points is approximately $e^{-tO(1)}$ and hence any corrections arising from it will be of a similar order of magnitude.) The leading order term for $I_3$ is identical to $I_4$, and consequently to leading order the contribution to integral $B$ from these two sections of the contour is simply twice the result given by (A53).

A3.6 *Estimates of Integral B when $\alpha = a, \alpha \approx a$ and $\alpha < a$.*

The particular case when $\alpha = a$ is important because only part of the above analysis of the previous two sections applies. When $\alpha = a$, then $pc = 1$, and the saddle point (A47) coincides with the initial point of the range of integration for $B$ in (A44). So to obtain an estimate for the latter one can simply integrate along the integration contour emanating from $w = 1$ as shown in Fig.1b, at an angle of $\pi/4$ radians outward towards a circle of infinite radius. Essentially this is nothing more than the integral $I_4$ as above, without the integrals along the other branches, but with the important difference that $d_2$ is zero. This means that $Re[\mathrm{i}\pi\alpha^2/4(w+1) + \mathrm{i}(t/2)Log(w)]$ behaves as $d_3u^3 = -tu^3/24\sqrt{2}$ along $w = 1 + ue^{\mathrm{i}\pi/4}$. The resulting integral is given by

$$I_4 = \frac{(1+\mathrm{i})}{\sqrt{2}}\int_0^\infty e^{d_3u^3} f(u)du = \frac{(1+\mathrm{i})}{\sqrt{2}}\int_0^\infty \frac{e^{d_3u^3}h(u)}{\left[1+u^2+u\sqrt{2}\right]^{1/8}\left[4+u^2+u2\sqrt{2}\right]^{3/4}}du, \quad \text{(A55)}$$

where

$$h(u) = exp\left[\sum_{i=5}^{\infty} d_i(1)u^i + \mathrm{i}\left\{\frac{2t(2+u/\sqrt{2})}{[4+u^2+u2\sqrt{2}]} + \frac{tlog(1+u^2+u\sqrt{2})}{4} - \frac{arctan(u/(\sqrt{2}+u))}{4} - \frac{3arctan(u/(2\sqrt{2}+u))}{2}\right\}\right]. \quad \text{(A56)}$$

Unfortunately the asymptotic expansion of (A55) is far more complicated than the relatively simple result (A53) for the case $\alpha > a$. That's because on expanding $f(u)$ in (A55) about $u = 0$, one finds that all the terms involving powers of $u^{3n}$ give rise (to leading order) contributions of $O(t^{-1/3})$, whilst those involving $u^{3n+1}$ and $u^{3n+2}$ contribute terms of $O(t^{-2/3})$ and $O(t^{-1})$ respectively. Hence it is necessary to calculate *all* the terms arising from the powers of $u^{3n}$ and (if possible) sum the infinite series. Unfortunately there seems no set pattern for these terms making such a summation problematic. Utilising the results (assuming $d_3 < 0$)

$$\int_0^\infty e^{d_3u^3}du = \frac{\Gamma(1/3)}{3|d_3|^{1/3}}, \quad \int_0^\infty ue^{d_3u^3}du = \frac{\Gamma(2/3)}{3|d_3|^{2/3}}, \quad \int_0^\infty u^2e^{d_3u^3}du = \frac{1}{3|d_3|}, \qquad \text{(A57)}$$

with similar expressions for integrals involving higher powers of $u$, one finds that to leading order

$$I_4 = \frac{\Gamma(1/3)\left(24\sqrt{2}\right)^{\frac{1}{3}}e^{\mathrm{i}t+\mathrm{i}\frac{\pi}{4}}}{3\times 4^{3/4}t^{1/3}}\left[\sum_{m=0}^{K}\frac{(-1)^m\Gamma(2m+1/3)}{\Gamma(1/3)(2m)!}+\mathrm{i}\left\{\sum_{m=1}^{K}\frac{(-1)^m\Gamma(2m-2/3)}{\Gamma(1/3)(2m-1)!}\right\}\right]+O\left(t^{-\frac{2}{3}}\right),$$
$$\underset{K\to\infty}{\longrightarrow}\ \frac{\Gamma(1/3)\left(24\sqrt{2}\right)^{\frac{1}{3}}e^{\mathrm{i}t+\mathrm{i}\frac{\pi}{4}}}{6\sqrt{2}t^{1/3}}\left[\frac{e^{-\mathrm{i}\frac{\pi}{12}}}{2^{1/6}}\right]+O\left(t^{-\frac{2}{3}}\right). \qquad \text{(A58)}$$

This result was established by calculating the first $K = 10$ terms of the two respective alternating series in (A58) long hand in order to establish a general pattern. Application of the convergence acceleration methods of [9] to these first 20 terms yields an estimate of $0.86054207 - \mathrm{i}0.2305815$ ($\pm 2.2\times 10^{-8} \pm \mathrm{i}1.3\times 10^{-7}$), for the respective partial sums. Assuming the terms continue indefinitely in the manner prescribed for $K > 10$, one obtains two convergent alternating series. In the limit as $K \to \infty$, one can easily establish an exact value of $2^{-1/6}e^{-\mathrm{i}\pi/12} = 0.86054208\ldots - \mathrm{i}0.23058155\ldots$ for these series, using the Maple 17 algebra package. Comparisons of numerical calculations of integral $B$ with $\alpha = a$ and (A58), for values of $t = 10^2, 10^3, 10^4$ and $10^5$, gives rise to relative errors of 1.3, 0.27, 0.059 and 0.013% (decaying as $t^{-2/3}$ as expected) respectively. Whilst not a proof, these results are indicative that for larger values of $t$ (for which this paper will be primarily concerned) estimate (A58) is likely to be an extremely good approximation for $B$ when $\alpha = a$. The representative nature of (A58) regarding the Riemann-Siegel sum is discussed in section 3.

Should $\alpha = a + \varepsilon$, where $\varepsilon = g t^{-1/6}$ and $g \in [0, 0.25]$, then the integral $B = I_3 + I_4$ can be estimated along the contour shown in Fig. A2b. The starting point of $I_3$ occurs at $u = u_2 \approx -(32\pi)^{1/4}\sqrt{g}/t^{1/3}$, where the real part of the phase of the integrand is made up of two terms $d_2u^2$ and $d_3u^3$ both of $O(1)$. Similar, but more complicated analysis to that used in the case when $\alpha = a$ above, yields the following estimate

$$B \approx \frac{\left(24\sqrt{2}\right)^{\frac{1}{3}}e^{\mathrm{i}t+\mathrm{i}\sqrt{\frac{\pi}{2}}t^{1/3}g+\mathrm{i}\frac{\pi}{4}}}{6\sqrt{2}t^{1/3}exp\left((32\pi^3)^{1/4}g^{3/2}/3\right)}\left[\frac{\Gamma(1/3)e^{-\mathrm{i}\frac{\pi}{12}}}{2^{1/6}}+\frac{\Gamma(2/3)e^{-\mathrm{i}\frac{\pi}{6}}}{2^{1/6}}3^{1/3}\sqrt{\pi}g+O(g^2)\right]+O\left(t^{-\frac{2}{3}}\right), \quad \text{(A59)}$$

with a (maximum) relative error of 24% when $\mathcal{g} = 0.25$. It is very difficult to obtain simple estimates to any greater precision for $\mathcal{g} \in [0.25, 1]$ because the complexity of the integrand (A49-50) gives rise to a power series in $\mathcal{g}$ on the right of (A59) which converges very slowly. It is better to resort to direct numerical computation of (A43) in these instances. Beyond the transitional region $\alpha = [a, a + t^{-1/6}]$, the standard approximation (A53) takes over. It is worth noting that $B$ initially increases until $\mathcal{g} = 0.53$, when it has a value $\approx 1.35$ times larger than (A58), before declining to a value ≈ 1.54 times smaller than (A58) when $\mathcal{g} = 1$. So within the transitional region one moves smoothly from (A58) via (A59) to (A53), with no dramatic changes in magnitude.

For the case when $\alpha = a - \varepsilon$ and $\varepsilon$ itself is extremely small $\leq 0.05/\sqrt{t}$ , substituting $\mathcal{g} \leq -0.05t^{-1/3}$ into (A59) yields a very accurate estimate for $B$. In this instance the contour of integral $I_4$ commences at points lying on the unit circle (Fig. A1). Finally for the case when $\varepsilon = \mathcal{g}t^{-1/6}$, one can estimate $B$ by simply integrating anticlockwise around the unit circle to a point with argument $\phi = arccos(2\alpha^2/a^2 - 1)$ (where the saddle point is shifted) before integrating outward towards infinity along a line at $\pi/4$ radians to the horizontal (Fig. A1a). (If $\alpha \leq a/\sqrt{2}$, one can just integrate around to $\phi = \pi/2$ and then along the positive imaginary axis.) As discussed in section A3.3 although integration around the unit circle contributes to $B$, it can make no contribution to Euler's integral formulation (A43) which appears in the zeta function series. Hence only the integration along the branch starting at $e^{i\phi}$ is relevant here. A general asymptotic analysis of the integral from this point yields a estimate of $O\left(e^{-t(\phi - sin\phi)/2}\right)$, an insignificant contribution unless $|\varepsilon| = \mathcal{g}t^{-1/6}$, when it is of $O\left(exp\left(-(2\pi^3)^{1/4}\mathcal{g}^{3/2}/3\right)\right)$. Again one can resort to direct numerical computation of (A43) in these instances.

Summarizing all these results, one finds that Euler's integral (A43) is exponentially small when $\alpha < a$, increases rapidly to a peak of $O\left(t^{-1/3}\right)$ between $\alpha = a \pm O\left(t^{-1/6}\right)$, and then decays, ultimately as $O\left(t^{1/4}\alpha^{-3/2}\right)$ for $\alpha \gg a$.

A3.7 *Hardy's Z Function.*

Substituting twice (A53) for $I_3 + I_4$ as the main estimate for $B$ into (A43), and utilising the fact that $\pi(\alpha^2 - 1)/4 \equiv 2m\pi$ when $\alpha$ is an odd integer, one finds

$$\int_0^1 \frac{e^{i\pi\alpha^2 x/4 + i(t/2)log\left\{\frac{(1-x)}{x}\right\}}}{[x(1-x)]^{\frac{1}{4}}} dx = \frac{e^{\frac{i\pi}{8}}(pc)^{3/4} 4\sqrt{\pi}\cos(\pi\alpha^2/4(pc+1) + (t/2)log(pc) + \pi/8)}{(pc+1)\sqrt{t(pc-1)}}[1 + O(\Lambda_\alpha^{-1})], \tag{A60}$$

where $\Lambda_\alpha = MIN\left\{t^{1/8}(\alpha - a)^{3/4}, \sqrt{t}\right\}$. In turn substituting (A60) into (A42a) and noting that $(-1)^{(\alpha-1)/2}e^{-i\pi\alpha/2 + i3\pi/8 + i\pi/8} = +1$ for all odd integers $\alpha$, gives the result

$$\mathrm{RSI}_{2^{nd}\,order} \approx \frac{\pi^{3/4}e^{-i\theta(t)}\mathcal{H}(t)}{(2t)^{1/4}} \sum_{\substack{\alpha > a \\ odd}}^{\infty} \left\{\frac{(pc)^{3/4}\alpha\cos(\pi\alpha^2/4(pc+1) + (t/2)log(pc) + \pi/8)}{(pc+1)\sqrt{t(pc-1)}}[1 + O(\Lambda_\alpha^{-1})]\right\}$$

(A61)

Finally, substituting this result into (A39) and noting $\pi\alpha^2/4(pc+1) = t(pc+1)/2pc$ from (A47) one obtains the expression

$$Z(t) \approx \mathcal{H}(t)2^{5/4}\pi^{1/4}\sum_{\substack{\alpha>a \\ odd}}^{\infty}\left\{\frac{(pc)^{1/4}\cos\big((t/2)\{log(pc)+1/pc\}+t/2+\pi/8\big)}{t^{1/4}\sqrt{pc-1}}[1+O(\Lambda_{\alpha}^{-1})]\right\}$$

(A62)

One further consequence of (A47) is that

$$\frac{2^{5/4}\pi^{1/4}(pc)^{1/4}}{t^{1/4}\sqrt{pc-1}} = \frac{2\sqrt{2}}{(\alpha^2-a^2)^{\frac{1}{4}}}, \quad \text{(A63)}$$

which means that (A62) can be written in the slightly more compact form

$$Z(t) \approx \mathcal{H}(t)2\sqrt{2}\sum_{\substack{\alpha>a \\ odd}}^{\infty}\left\{\frac{\cos\big((t/2)\{log(pc)+1/pc\}+t/2+\pi/8\big)}{(\alpha^2-a^2)^{\frac{1}{4}}}[1+O(\Lambda_{\alpha}^{-1})]\right\}. \quad (A64)$$

This new zeta-sum representation for $Z(t)$ is the *main result* of this Appendix and forms the starting point of the main paper (Section 2). The reason for the ≈ rather than an equals sign in (A64) is that whenever $\alpha = a + \varepsilon$ the transitional term arising from (A58-59) reaches significant size and most be accounted for. For $\varepsilon = g t^{-1/6}$ and $g \in [0, 0.25]$, following through the same procedure from (A60) to (A64) as above, gives the approximation

$$T_{a+\varepsilon} = \frac{\mathcal{H}(t)2^{3/4}exp\big(-(32\pi^3)^{1/4}g^{3/2}/3\big)}{3^{2/3}\pi^{1/4}t^{1/12}}\left[\Gamma\left(\frac{1}{3}\right)cos\left(t+\sqrt{\frac{\pi}{2}}t^{1/3}g+\frac{\pi}{24}\right)\right.$$
$$\left.+\,\Gamma\left(\frac{2}{3}\right)3^{1/3}\sqrt{\pi}g\,cos\left(t+\sqrt{\frac{\pi}{2}}t^{1/3}g-\frac{\pi}{24}\right)+O(g^2)\right]. \quad \text{(A65)}$$

which should be added to (A64) every time $a$ happens to lie close to an odd integer. In the interval $g \in (-1, 0) \cup (0.25, 1)$ this approximation of the transitional term is increasingly inaccurate and it is better (and much easier) to evaluate it numerically.

A3.8 *Convergence Properties of the new series for* $Z(t)$.

For $\alpha \gg t,\ pc \approx \pi\alpha^2/2t$ and the terms in the series (A64) approximate to

$$\frac{\cos[t/2\,(log(\pi/2t)+1)+\pi/8+t\,log(\alpha)]}{\sqrt{\alpha}}. \quad \text{(A66)}$$

Hence the series (A64) cannot converge in the limit as $\alpha \to \infty$ and it is necessary to consider (A64) as the analytic continuation of a corresponding analytical function derived from the $Re\left[e^{\mathrm{i}\theta(t)}\mathrm{RSI}(s)\right]$ for any $s \in \mathbb{C}$ with $Re(s) > 1$, in place of $s = 1/2 + \mathrm{i}t$. From a computational point of view this presents no difficulties. All that is required is the termination the series (A64) at some finite odd integer value $N_\alpha$ and then utilisation of the Euler-Macluarin summation formula ([1] or [13]) to evaluate the remainder. However, terminating series (A64) at a finite value has certain implications in relation to sums (A33 & A38), and the convergence properties of the asymptotic expansion (C11) derived from it in order to estimate the $Im\left[e^{\mathrm{i}\theta(t)}\mathrm{RSI}(1/2 + \mathrm{i}t)\right]$ for large $t$. The unproven assertion that the terms of $O(e^{\alpha\sqrt{\pi t}/\sqrt{2}})$ in the expansion should cancel for all orders of reciprocal powers of $t$, would imply that (A38) can be treated as an infinite sum for $Im\left[e^{\mathrm{i}\theta(t)}\mathrm{RSI}(1/2 + \mathrm{i}t)\right]$. However, the necessity of treating (A64) as a finite sum for the $Re\left[e^{\mathrm{i}\theta(t)}\mathrm{RSI}(1/2 + \mathrm{i}t)\right]$, suggests such a termination might be necessary for the imaginary part as well. Speculatively, this might be necessary because at some point the $O(e^{\alpha\sqrt{\pi t}/\sqrt{2}})$ terms cease to cancel and one must invoke ideas of analytic continuation in order to evaluate a (finite) value for $Im\left[e^{\mathrm{i}\theta(t)}\mathrm{RSI}(1/2 + \mathrm{i}t)\right]$. In that case (C11) might be an example of an ultimately divergent asymptotic series, much like Stirling's series for the gamma function. If this was so, logic suggests that the terms of $O(e^{\alpha\sqrt{\pi t}/\sqrt{2}})$ should cease to cancel at some power of $t^{-n}$ when $e^{N_\alpha\sqrt{\pi t/2}}t^{-n}$ is less than the leading order term of (C11), which is $e^{-\sqrt{\pi t/2}}t^{-1/4}$. That is when $n$ reaches a value greater than $(N_\alpha + 1)\sqrt{\pi t/2}/log(t)$. However, because (as will be seen) $N_\alpha \sim O(t)$, this estimated bound diverges as $t \to \infty$ and hence $n \to \infty$ also, the implication is that the cancellation of the $O(e^{\alpha\sqrt{\pi t}/\sqrt{2}})$ terms does indeed continue indefinitely. So in conclusion it must be correct to treat (A38) as a convergent infinite sum for the purposes of evaluating $Im\left[e^{\mathrm{i}\theta(t)}\mathrm{RSI}(1/2 + \mathrm{i}t)\right]$, but for the purposes of evaluating $Re\left[e^{\mathrm{i}\theta(t)}\mathrm{RSI}(1/2 + \mathrm{i}t)\right]$ one must terminate (A38) at $N_\alpha$ and appeal to analytical continuation. This dichotomy underlies the interpretation of the RSI as a sum when formulating (A21). Incidentally there are no problems surrounding the convergence properties of the smallest term in (A38), namely $e^{-2\pi t}\Phi(1/4 - \mathrm{i}t/2, 1/2; \pi\alpha^2\mathrm{i}/4)$. That is because for $\alpha \gg t^2$ the confluent hypergeometic function behaves as a $constant/\alpha^{3/2\pm\mathrm{i}t}$, yielding an absolutely convergent sum over $\alpha$.

## A4. Sample Computations of the new series for $Z(t)$

### A4.1 *Evaluation of the new series by means of Euler-Maclaurin Summation*

One can use (A64) to make computational estimates of Hardy's $Z$ function for large $t$ by means of the standard technique of Euler-Maclaurin summation ([1] or [13]). This provides an estimate for the general sum of the form

$$\sum_{i=L}^{\infty} f(i) = \int_L^{\infty} f(x)dx + \frac{1}{2}[f(L) + f(\infty)] + \sum_{j=1}^{l} \frac{B_{2j}}{(2j)!}\left[f^{(2j-1)}(\infty) - f^{(2j-1)}(L)\right] + E_{2l}, \quad \text{(A67)}$$

where the error term takes the form

$$E_{2l} = \frac{1}{(2l+1)!}\int_L^{\infty} \bar{B}_{2l+1}(x) f^{(2l+1)}(x)dx\,. \tag{A68}$$

Here $B_2 = 1/6,\ B_4 = -1/30,\ B_6 = 1/42$ ... are the Bernoulli numbers and $\bar{B}_{2l+1}(x)$ is the periodic function $B_{2l+1}(x - INT(x))$, where $B_{2l+1}$ denotes $(2l+1)$th Bernoulli polynomial ([1], Tables 23.1&2). The formula (A67) cannot be applied directly to estimate (A64) because at the lower limit $\alpha \approx a$ the Bernoulli sum will diverge. To find the least upper bound on $\alpha$ for which (A67) is applicable, let

$$f\big(pc(\alpha = 2M+1)\big) = cos\{(t/2)\{log(pc) + 1/pc\} + t/2 + \pi/8\}/(\alpha^2 - a^2)^{1/4}. \tag{A69}$$

If $f^{(n)} \equiv d^{(n)}f/dM^{(n)}$, then using the fact that

$$(pc-1)\frac{d(pc)}{d\alpha} = pc\sqrt{\frac{2\pi pc}{t}} \tag{A70}$$

from (A47), plus identity (A63), one finds that to leading order

$$f^{(2n-1)} \sim (-1)^n \left(\frac{\pi pc}{2t}\right)^{1/4} \left(\sqrt{\frac{2\pi t}{pc}}\right)^{2n-1} \frac{sin\{(t/2)\{log(pc) + 1/pc\} + t/2 + \pi/8\}}{\sqrt{pc-1}}. \tag{A71}$$

Consequently

$$\sum_{j=1}^{l} \frac{B_{2j}}{(2j)!} f^{(2j-1)}\big(pc(2M+1)\big) \sim -\frac{(pc)^{3/4} sin\{(t/2)\{log(pc) + 1/pc\} + t/2 + \pi/8\}}{(2t)^{3/4}\pi^{1/4}\sqrt{pc-1}} \sum_{j=1}^{l}\left(\frac{2\pi t}{pc}\right)^j \frac{|B_{2j}|}{(2j)!}, \tag{A72}$$

and given that $|B_{2j}| \sim 2(2j)!/(2\pi)^{2j}$, together with $pc(\alpha) \approx \pi\alpha^2/2t$ , one can ensure that the Bernoulli sum (A72) converges provided $\alpha = N_\alpha = 2M_\alpha + 1 \geq INT_O(t/\pi) + 2 \gg a$. So to utilise (A64) as a means of computation one must explicitly calculate the 'main sum' of the terms for $\alpha \in [INT_O(a) + 2, N_\alpha]$, followed by application of the Euler-Maclaurin formula to the remaining terms $\alpha \geq K = N_\alpha + 2$, to give

$$Z(t) \approx \mathcal{H}(t)2\sqrt{2}\left\{main\ sum + I + \sum_{j=1}^{l}\frac{B_{2j}f^{(2j-1)}\big(pc(K)\big)}{(2j)!} + \frac{f\big(pc(K)\big)}{2} + E_{2l}\right\}, \tag{A73}$$

with the transitional term (A65) when necessary. Here the integral $I$ in (A73) can be estimated as

$$I = \frac{1}{2}\int_K^{\infty} f\big(pc(\alpha)\big)\,d\alpha = -\frac{sin\big\{(t/2)\big\{log\big(pc(K)\big) + 1/pc(K)\big\} + t/2 + \pi/8\big\}}{\pi^{1/4}\sqrt{pc(K)-1}}\left(\frac{pc(K)}{2t}\right)^{3/4}\left\{1 + O\left(\frac{1}{t}\right)\right\}, \tag{A74}$$

using Laplace's method.

Table AI shows some illustrative calculations of the new zeta-sum representation of $Z(t)$ using the Euler Maclaurin Summation formula (A73) starting from $t = 1000$. Here $N_\alpha$ was set to $INT_O(0.35t) + 2$ throughout, except where stated otherwise‡. The two odd integer terms lying either side of $a$ were calculated as follows:

a) If $|a - \alpha| < g t^{-1/2}$ then equation (A65) was utilised;
b) If $\alpha < a - g t^{-1/2}$, or $(\alpha - a) \in [g t^{-1/2}, t^{-1/6}]$, then the term was calculated numerically;
c) If $\alpha > a + t^{-1/6}$, then the term was calculated from (A64).

For comparison the actual values of $Z(t)$ were calculated using the full RS formula

$$Z(t) = 2\sum_{k=1}^{N_t} \frac{cos\{\theta(t) - tlog(k)\}}{\sqrt{k}} + (-1)^{N_t - 1}\left(\frac{t}{2\pi}\right)^{-1/4} \sum_{r=0}^{m} (-1)^r \left(\frac{t}{2\pi}\right)^{-r/2} \Psi_r(p) + R_m(t), \tag{A75}$$

the most efficient known method for computing the zeta function to reasonable precision. In (A75), $N_t = INT[(t/2\pi)^{1/2}]$, $p = (t/2\pi)^{1/2} - N_t$, $\Psi_0(p) = cos\{2\pi(p^2 - p - 1/16)\}/cos\{2\pi p\}$ and $\Psi_r(p)$ are combinations of the derivatives of the function $\Psi_0(p)$. (For bounds on the remainder term $R_m(t)$ , see [4], [5], and [14]. Typically $|R_2(t)| < 0.011t^{-7/4}$ for all $t > 200$.)

The choice of $t$ values in Table I is designed to illustrate the care which must be taken when utilising (A73) to estimate $Z(t)$. So for instance when $t = 1100$, it turns out that $a = 52.92$ which is within $t^{-1/6}$ of $\alpha = 53$ the first odd integer in the main sum. If this first term is calculated directly from (A64) it leads to a relative error of 29% in the estimated value of $Z(t)$. Numerical calculation of this first term reduces this error to 0.45% (see the second and third lines of the table). For the case when $a =$ an odd integer, take for example $t$ =1103.091720 when $a = 53$ exactly. Here the first term in the series at $\alpha = a$ was calculated using (A65) giving rise to an estimate for $Z(t)$ which is in error by 1.2%. This is somewhat larger than the error for the nearby value $t = 1100$ but not unexpected. That is because the error in (A65) is term of only $O(t^{-1/3})$ smaller than (A65) itself (see equation A58), whereas the errors in the terms in (A64) are $O(t^{-1/2})$ smaller in magnitude. Of course the influence of this $O(t^{-1/3})$ error in (A65) will decline as $t$ increases. Line seven of the table shows the case when $t$ =100148.083310 when $a = 505$ exactly. Using (A65) to compute the first term leads to a much reduced relative error 0.0086%.

A couple of interesting cases are the values $t = 17143.803905$ [13] and $t = 388858886.002$ [32] corresponding to two very small local maxima/minima, separating two very closely spaced zeros. The calculation for $t = 17143.803905$ is also complicated by the fact that $a = 208.94$ within $t^{-1/6}$ an odd integer and so neither (A64) nor (A65) give a particularly good estimate for the first term in the main sum. Hence this was calculated numerically. The calculation for $t = 388858886.002$ is more straightforward as $a = 31467.77$ is not close to an odd integer and the first term can be calculated from (A64) directly. Both calculations result in large relative errors (but small absolute errors) of 48.7% and $\approx 30000\%$ respectively.

| $t$ ($a$) | Main sum ($\alpha$ range) | Bernoulli Sum $l = 60$ | $\frac{f(pc(K))}{2}$ | $I$ (equ. 79) | $Z(t) \approx \sqrt{8} \times Total$ | Actual $Z(t)$ |
|---|---|---|---|---|---|---|
| 1000.0 (50.46) | 0.26431† (49-351) | $7.6192 \times 10^{-2}$ | $1.6683 \times 10^{-2}$ | $-7.3434 \times 10^{-3}$ | 0.98950 | 0.99779 |
| 1100.0 (52.92) | −0.49547 (51-387) | $-9.0923 \times 10^{-2}$ | $2.7258 \times 10^{-4}$ | $8.9646 \times 10^{-3}$ | −1.63245 | −1.26328 |
| 1100.0 (52.92) | −0.36698† (51-387) | $-9.0923 \times 10^{-2}$ | $2.7258 \times 10^{-4}$ | $8.9646 \times 10^{-3}$ | −1.26902 | −1.26328 |
| 1103.091720 (53.00) | 0.48238* (53-387) | $4.0577 \times 10^{-2}$ | $2.2932 \times 10^{-2}$ | $-3.8899 \times 10^{-3}$ | 1.54950 | 1.56826 |
| 17143.803905 (208.94) | $5.922 \times 10^{-3}$† (207-8573)‡ | $-3.7111 \times 10^{-3}$ | $-3.7595 \times 10^{-3}$ | $1.9388 \times 10^{-3}$ | $1.104 \times 10^{-3}$ | $2.153 \times 10^{-3}$ |
| 100000.0 (504.63) | 2.0833 (503-35001) | $-8.0095 \times 10^{-3}$ | $1.6124 \times 10^{-3}$ | $7.4604 \times 10^{-4}$ | 5.87656 | 5.87959 |
| 100148.083310 (505.00) | 2.7636* (505-35053) | $-1.0028 \times 10^{-2}$ | $1.0334 \times 10^{-3}$ | $9.3406 \times 10^{-4}$ | 7.79120 | 7.79053 |
| 2000000.0 (2256.76) | −0.80451 (2255-700001) | $1.1829 \times 10^{-3}$ | $-5.0804 \times 10^{-4}$ | $-1.1014 \times 10^{-4}$ | −2.27389 | −2.27469 |
| $10^7$ (5046.26) | 5.07396 (5045-3500001) | $1.8778 \times 10^{-5}$ | $2.6721 \times 10^{-4}$ | $-1.7484 \times 10^{-6}$ | 14.35212 | 14.35255 |
| 388858886.002 (31467.77) | $6.0029 \times 10^{-5}$ (31465-194429445) ‡ | $-1.5892 \times 10^{-6}$ | $-3.5820 \times 10^{-5}$ | $8.2975 \times 10^{-7}$ | $6.6326 \times 10^{-5}$ | $-2.2183 \times 10^{-7}$ |
| $10^9$ (50462.65) | −1.14258 (50461-35000001) | $3.8502 \times 10^{-5}$ | $-2.4686 \times 10^{-5}$ | $-3.5946 \times 10^{-6}$ | −3.23166 | −3.23130 |

Table AI. Some Illustrative Calculations of Hardy's $Z$ Function based on the Euler-Maclaurin Summation of the new zeta-sum given by equation (A64). The actual values of $Z(t)$ in the last column are calculated using the RS formula.

* First term in main sum estimated using (A65).

† Second term in main sum ($\alpha$ closest to $a$) calculated by numerical integration.

‡ $N_\alpha = INT_O(t/2) + 2$ to improve convergence of Bernoulli sum.

Clearly the computation of $Z(t)$ in these instances would require at least the inclusion of correction terms of $O(t^{-1/2})$ in the main sum (A64) to obtain agreement even to one significant figure. However, it would be pointless to undertake such an analysis because the utility of adopting (A64) as a basis for the calculation of $Z(t)$ is *negligible*. The RS formula (A75), with just $\sqrt{t/2\pi}$ terms in its main sum compared to a minimum of $t/2\pi$ in (A64), will always be *vastly* more efficient. The crucial point about Table AI is not the efficiency of the calculations, but rather the computational support it provides that the formulation (A21) made for the RSI in devising (A64) and (A65) is correct, in the sense that it appears to yield an asymptotic approximation for $Z(t)$, with an absolute error tending to zero as $t \to \infty$. The proof of this statement is the subject of the first part of the main paper. But these results raise the question just how, from a computational perspective, does (A64) bring this about? On one level the answer is obvious, because the fact that RSI satisfies the functional equation of the zeta function means that any asymptotic approximation arising from the RSI must converge on $Z(t)$. But on another level this must mean there is some *connection* between the main sum in (A64) and the main sum in the RS formula (A75), which is by no means obvious. (Not only are the lengths of the sums different, but the magnitudes of the individual terms and phases of the cosines seem completely unrelated.) Uncovering such a connection not only explains how, computationally, (A64) yields an approximation for $Z(t)$, but yields the *remarkable* conclusion that the RS formula, thought since its discovery to be the fastest possible method for computing the zeta function for large $t$, is in fact *inefficient* and can be significantly improved upon.

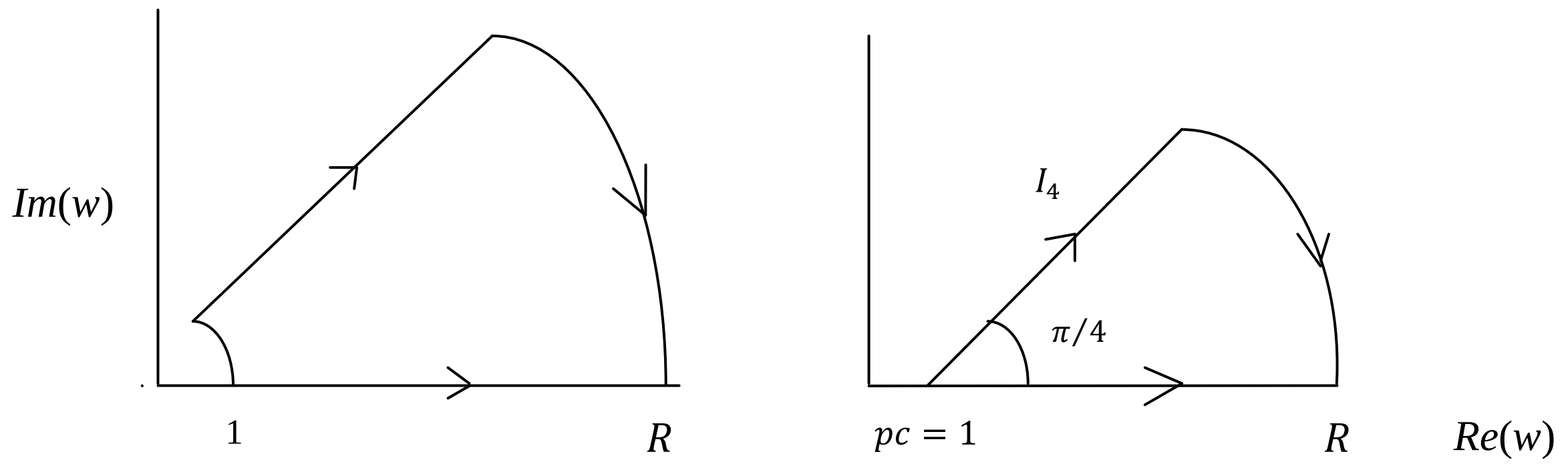

Fig. A1a Fig. A1b

Figure A1. Schematic of the general contour of integration in $w$ space used to estimate integral $B$ equation (A44), for the cases a) when $\alpha < a$ and b) when $\alpha = a \equiv pc = 1$.

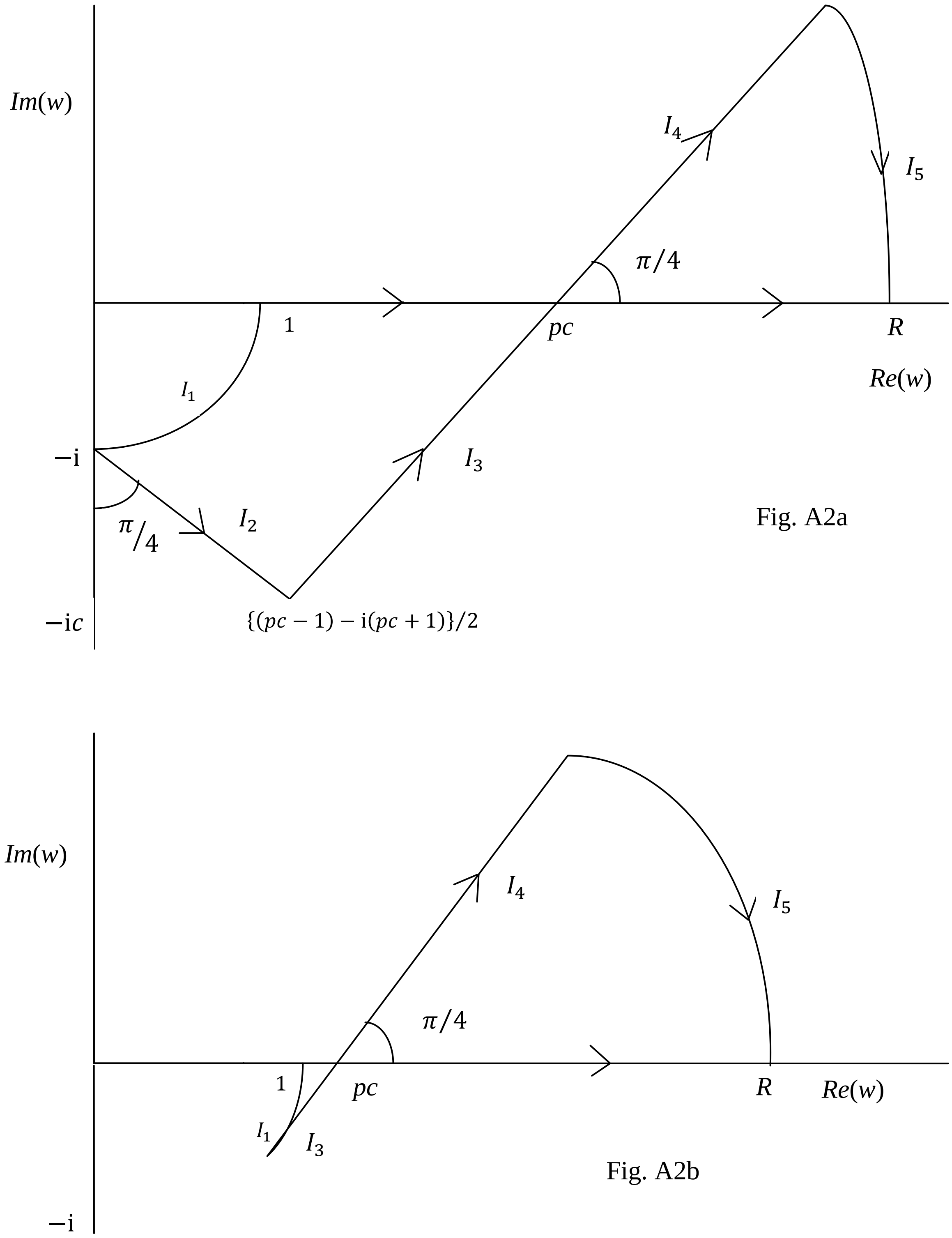

Figure A2. Schematic of the general contour of integration in $w$ space used to estimate integral $B$ equation (A44), when $\alpha > a$ for the cases a) $pc > \sqrt{2}$ and b) $pc \leq \sqrt{2}$.

**Appendix B. Asymptotic estimation of the Integrals $C(N,R)$ and $D(N,R)$ in Section 2.3**

B1.1 *Estimation of the integral* $C(N,R)$.

*Lemma B1.1.* Let $C(N,R)$ be the integral (25) of Section 2.3, defined by

$$C(N,R)=\int_{\frac{t}{\pi R}}^{\infty}\frac{exp\left[-\mathrm{i}tN(1+1/u)+\mathrm{i}(t/2)\left(log\big(pc(u)\big)+1/pc(u)\right)\right]}{u^2[(1+1/u)^2-8\pi/t]^{1/4}}du\,, \qquad \text{(B1)}$$

where $pc(u)\equiv pc(z=t(1+1/u)/\pi)$ as defined by (6, A47), $N=0,1,2,\ldots$, and $R=t/\pi-O\left(t^{1/2}\right)$. Provided $N\neq 1$, then $C(N,R)=O(\{t|1/2-N|\}^{-1})+O(\{t|N-1|\}^{-1})$.

*Proof.* The integrand in (B1) is highly oscillatory and decays only algebraically as $u\to\infty$. However, in the region $Re(u)\geq t/\pi R$ it has no poles and it proves possible to estimate $C(N,R)$ by integrating around two particular contours in the complex plane.

First consider the Taylor expansion of the phase of the numerator along the lines $u=t/\pi R\pm \mathrm{i}\delta$ where $\delta>0$ and small. One obtains (using $dpc/du=-apc^{3/2}/2(pc-1)u^2$)

$$-\mathrm{i}t\left\{N\left(1+\frac{\pi R}{t}\right)-\frac{\left(log(pc_+)+\frac{1}{pc_+}\right)}{2}\right\}\pm\delta\frac{(\pi R)^2}{2t}\left\{\frac{a}{2\sqrt{pc_+}}-2N\right\}+\mathrm{i}\delta^2\frac{(\pi R)^3}{2!\,t^2}\left\{\frac{a}{2\sqrt{pc_+}}-2N+\frac{R}{2(pc_+-1)}\right\}$$

$$\mp\delta^3\frac{(\pi R)^4}{3!\,t^3}\left\{3\left(\frac{a}{2\sqrt{pc_+}}-2N\right)-\frac{3R}{(pc_+-1)}+\frac{\pi R^2pc_+^{\frac{3}{2}}}{\sqrt{2\pi t}(pc_+-1)^3}\right\}\ldots, \qquad \text{(B2)}$$

where $pc_+=pc(t(1+\pi R/t)/\pi)$. The corresponding asymptotic expansion when $\delta\to\infty$ (assuming $t$ is large but with $t/\pi R\approx 1$), is given by

$$-\mathrm{i}\left(tN-\frac{t}{2}log\left(\frac{t}{2\pi}\right)+\pi+\frac{2\pi^2}{t}\ldots\right)\mp\frac{1}{\delta}\left[t(N-1)-2\pi-\frac{8\pi^2}{t}\right]$$
$$-\frac{\mathrm{i}}{\delta^2}\left[t\left\{\left(\frac{t}{\pi R}\right)(N-1)-\frac{1}{2}\right\}-\pi\left\{2\left(\frac{t}{\pi R}\right)+3\right\}-\frac{\pi^2}{t}\left\{8\left(\frac{t}{\pi R}\right)+20\right\}\right]+O\left(\frac{t(N-1)}{\delta^3}\right). \qquad \text{(B3)}$$

Notice in (B3) that the magnitude of the term of $O(1/\delta)$ dominates the remaining terms provided $\delta>t/\pi R\approx 1$. So the combination of (B2) for $\delta<1$ and (B3) for $\delta>1$, together provide a good approximation of the behaviour of the phase of the numerator for all positive $\delta$ values.

Now in (B2) the assumption $t/\pi R\approx 1\Rightarrow a/2\sqrt{pc_+}\approx 1$ also, which means that for $N\geq 2$ moving along the line $u=t/\pi R+\mathrm{i}\delta$ the real part of the phase is initially negative. From (B3), one can see that the real part of the phase *continues to be negative* for all positive values of $\delta$. Consequently the integrand in (B1) along the line $u=t/\pi R+\mathrm{i}\delta$ will be exponentially small along its entire length and hence one can estimate the integral along this line by local methods. A similar conclusion applies for the case $N = 0$ when moving along the line $u = t/\pi R - \mathrm{i}\delta$. The case $N = 1$ is considerably more complicated and is discussed separately in section B1.6. These observations suggest that in order to estimate (B1) one should consider integrating around the two closed contours illustrated in Fig. B1a,b. The contribution of the integrals around the quarter circles $Se^{\mathrm{i}\Phi}$, $\phi \in [-\pi/2\,, 0]$ and $Se^{\mathrm{i}\Phi}$, $\phi \in [0,\ \pi/2]$ must tend to zero as $S \to \infty$, because along these sections the integrand is of $O(1/S^2)$ and the sections themselves are of length $\pi S/2$. (In the limit $S \to \infty$, the phase of the numerator is purely imaginary on these quarter circles because at such points $\delta \sim O(S)$ in equation B3.) So the value of the integral $C(N,R)$ along the real line will be equivalent to the contribution of the respective integrals along the lines $t/\pi R \pm \mathrm{i}\delta$ only. It remains to establish a suitable approximation for the latter. This comes in two parts of similar magnitudes, one part arising from the integrand's behaviour near $\delta = 0$ when the real part of the numerator phase is approximately zero, and one part arising from a contribution as $\delta \to \infty$.

Near $\delta = 0$ the denominator of (B1) can be approximated by

$$\left(\frac{\pi R}{t}\right)^2 \left[\left(1+\frac{\pi R}{t}\right)^2 - \frac{8\pi}{t}\right]^{-1/4} \left\{1 \pm \mathrm{i}\delta \left(\left(\frac{\pi R}{t}\right)^2 \frac{(1+\pi R/t)}{2}\left[\left(1+\frac{\pi R}{t}\right)^2 - \frac{8\pi}{t}\right]^{-1} - \frac{2\pi R}{t}\right) \ldots\right\}, \quad \text{(B4)}$$

giving rise a contribution to the integral near $\delta = 0$ of the form

$$\frac{-\mathrm{i}exp[-\mathrm{i}t\{N(1+\pi R/t) - (log(pc_+) + 1/pc_+)/2\}]}{(t/2)((1+\pi R/t)^2 - 8\pi/t)^{1/4}\left(a/2\sqrt{pc_+} - 2N\right)}\left\{1 + \frac{\mathrm{i}(2\pi R/t - (\pi R)^2(1+\pi R/t)[(1+\pi R/t)^2 - 8\pi/t]^{-1}/2t^2)}{\left(a/2\sqrt{pc_+} - 2N\right)(\pi R)^2/2t}\right\} + O\left(\frac{1}{\left(t\left(a/2\sqrt{pc_+} - 2N\right)\right)^3}\right)$$

(B5)

Now when $\delta \to \infty$ the denominator of (B1) tends to

$$-\frac{1}{\delta^2(1-8\pi/t)^{1/4}} \mp \frac{\mathrm{i}}{\delta^3(1-8\pi/t)^{1/4}}\left(\frac{2t}{\pi R} + \frac{1/2}{(1-8\pi/t)}\right)\ldots, \quad \text{(B6)}$$

which together with (B3) gives rise to a contribution to the integral near the end of each line of the form

$$\pm \mathrm{i}exp\left[-\mathrm{i}\left(tN - \frac{t}{2}log\left(\frac{t}{2\pi}\right) + \pi + \frac{2\pi^2}{t}\right)\right]\int_{\delta_L}^{\infty} e^{\mp\frac{1}{\delta}[t(N-1) - 2\pi - 8\pi/t]}\left\{-\frac{1}{\delta^2(1-8\pi/t)^{1/4}} \mp \frac{\mathrm{i}}{\delta^3(1-8\pi/t)^{1/4}}\left(\frac{2t}{\pi R} + \frac{1/2}{(1-8\pi/t)}\right)\right\}d\delta. \quad \text{(B7)}$$

Assuming the lower limit $\delta_L$ lies in the range $1 \ll \delta_L \ll t$, integration of (B7) means that limit will gives rise to a very small $e^{-t/\delta_L}$ term provided $N \neq 1$. Consequently, integration of (B7) gives rise to the following leading order approximation

$$\frac{\mathrm{i}e^{-\mathrm{i}\left[tN-\frac{t}{2}log\left(\frac{t}{2\pi}\right)+\frac{2\pi^2}{t}\right]}}{(1-8\pi/t)^{1/4}[t(N-1)-2\pi-8\pi/t]}+O\left(\frac{1}{[t(N-1)-2\pi-8\pi/t]^2}\right). \qquad \text{(B8)}$$

So provided $N \neq 1$, an asymptotic approximation for the integral $C(N,R)$ can be found simply by adding results (B5) and (B8). To leading order (B5) is of $O(\{t|1/2-N|\}^{-1})$, whilst (B8) is of $O(\{t|N-1|\}^{-1})$, which proves *Lemma B1.1*.

B1.2 *Estimation of the integral $D(N,R)$, general points.*

The second and more difficult integral $D(N,R)$, (26) of Section 2.3, is defined by

$$D(N,R)=\int_{\frac{t}{\pi R}}^{\infty}\frac{exp\left[-\mathrm{i}tN(1-1/u)+\mathrm{i}(t/2)\left(log\big(pc(u)\big)+1/pc(u)\right)\right]}{u^2[(1-1/u)^2-8\pi/t]^{1/4}}du\,, \qquad \text{(B9)}$$

where $pc(u)\equiv pc(z=t(1-1/u)/\pi)$. Now the corresponding Taylor expansion to equation (B2) of the numerator phase along the lines $u=t/\pi R\pm \mathrm{i}\delta$ where $\delta>0$ and small gives (using $dpc/du=apc^{3/2}/2(pc-1)u^2$)

$$-\mathrm{i}t\left\{N\left(1-\frac{\pi R}{t}\right)-\frac{\left(log(pc_-)+\frac{1}{pc_-}\right)}{2}\right\}\pm\delta\frac{(\pi R)^2}{2t}\left\{2N-\frac{a}{2\sqrt{pc_-}}\right\}+\mathrm{i}\delta^2\frac{(\pi R)^3}{2!\,t^2}\left\{\frac{a}{2\sqrt{pc_-}}-2N+\frac{R}{2(pc_--1)}\right\}$$

$$\mp\delta^3\frac{(\pi R)^4}{3!\,t^3}\left\{3\left(2N-\frac{a}{2\sqrt{pc_-}}\right)-\frac{3R}{(pc_--1)}-\frac{\pi R^2pc_-^{\frac{3}{2}}}{\sqrt{2\pi t}(pc_--1)^3}\right\}\ldots, \qquad \text{(B10)}$$

where $pc_-=pc(t(1-\pi R/t)/\pi)$. The corresponding asymptotic expansion to (B3) of this phase along the lines $u=t/\pi R\pm \mathrm{i}\delta$ with $\delta\to\infty$ and $t$ large, is given by

$$-\mathrm{i}\left(tN-\frac{t}{2}log\left(\frac{t}{2\pi}\right)+\pi+\frac{2\pi^2}{t}\ldots\right)\pm\frac{1}{\delta}\left[t(N-1)-2\pi-\frac{8\pi^2}{t}\right]$$
$$+\frac{\mathrm{i}}{\delta^2}\left[t\left\{\left(\frac{t}{\pi R}\right)(N-1)+\frac{1}{2}\right\}-\pi\left\{2\left(\frac{t}{\pi R}\right)-3\right\}-\frac{\pi^2}{t}\left\{8\left(\frac{t}{\pi R}\right)-20\right\}\right]+O\left(\frac{t(N-1)}{\delta^3}\right).$$

(B11)

Now if the radius $R$ is chosen to be just a bit less than $t/\pi-a$, this will mean that $pc_-$ is very close to unity. So provided $N>a/4\sqrt{pc_-}$ then the real part of the phase (B10 &11) remains negative along the entire length of the line $u=t/\pi R-\mathrm{i}\delta$ and the integral can be estimated by much the same methodology as that used to approximate $C(N,R)$. A similar observation applies for the case when $N = 0$ along the line $u = t/\pi R + \mathrm{i}\delta$. However, for those integers $N \in \left[2, INT\left(a/4\sqrt{pc_-}\right)\right]$ things are more complicated (the special case $N = 1$ is again dealt with in section B1.6). Moving up the line $u = t/\pi R + \mathrm{i}\delta$, initially the real part of the phase is negative (see B10), but as $\delta$ increases it is apparent that at some point a change sign must occur (see B11). Consequently the methodology adopted in the previous section to estimate $C(N, R)$ is not applicable to $D(N, R)$ when $N$ lies in this particular range. Rather one needs to adapt the contour shown in Fig. B1b so that it begins along the line $u = t/\pi R + \mathrm{i}\delta$ for small $\delta$, but ends up along the line $u = t/\pi R - \mathrm{i}\delta$ for large $\delta$, in such a way as to ensure the phase's real part never becomes positive. To achieve this whilst at the same time ensuring that the contour remains confined to the region $Re(u) \geq t/\pi R$, it must cross the real axis at some suitable point. To identify such a point replace $\delta$ by $-\mathrm{i}\delta$ in (B10) to give the behaviour of the phase in the vicinity of some general point on the real axis $u = u_* + \delta$. This gives the result

$$-\mathrm{i}t\left\{N\left(1-\frac{1}{u_*}\right)-\frac{\left(log(pc_*)+\frac{1}{pc_*}\right)}{2}\right\}-\mathrm{i}\delta\frac{t}{2u_*^2}\left\{2N-\frac{a}{2\sqrt{pc_*}}\right\}-\mathrm{i}\delta^2\frac{t}{2!\,u_*^3}\left\{\frac{a}{2\sqrt{pc_*}}-2N+\frac{t/\pi u_*}{2(pc_*-1)}\right\}$$

$$-\mathrm{i}\delta^3\frac{t}{3!\,u_*^4}\left\{3\left(2N-\frac{a}{2\sqrt{pc_*}}\right)-\frac{3\,t/\pi u_*}{(pc_*-1)}-\frac{\pi(t/\pi u_*)^2 pc_*^{\frac{3}{2}}}{\sqrt{2\pi t}(pc_*-1)^3}\right\}\dots, \quad \text{(B12)}$$

where $pc_* = pc(t(1 - 1/u_*)/\pi)$. Examination of (B12) suggests that a suitable point for the contour to cross the real axis is at the saddle point $u_* = u_{sad}$ defined by the equation $2N - a/2\sqrt{pc_{sad}} = 0$. If the contour can be arranged to pass through this saddle in such a way that the real part of the phase will always be non-negative, then the contribution to the integral will be heavily concentrated in the vicinity of $u = u_{sad}$ and can be calculated using standard saddle point methods. This preliminary analysis suggests that in order to estimate the integral $D(N, R)$ when $N \in \left[2,\ Int\left(a/4\sqrt{pc_-}\right)\right]$ a suitable contour would be as illustrated in Fig. B2. These considerations lead to the following lemma.

*Lemma B1.2* Let $D(N, R)$ be the integral defined by (B9) and assume $N \neq 1$. Let $2\eta > \mu > \eta > 0$ be $O(1)$ constants such that $R \approx_p R_{max} = t/\pi - a(pc_- + 1)/2\sqrt{pc_-}$, where $pc_- = 1 + \eta t^{-1/3}\sqrt{8\pi}$. Let $N_\mu = \left\lfloor\sqrt{t/2\pi} - \mu t^{\frac{1}{6}}\right\rfloor$, $N_\eta = \left\lfloor\sqrt{t/2\pi} - \eta t^{\frac{1}{6}}\right\rfloor = INT\left(a/4\sqrt{pc_-}\right)$, $N_t = \left\lfloor\sqrt{t/2\pi}\right\rfloor$, so that $N_{2\eta} < N_\mu < N_\eta < N_t$. Let $H(x)$ be the Heaviside step function and define

$$Y(N,t) = \sqrt{\frac{2\pi}{Nt}}e^{-\mathrm{i}\pi/4}exp\left(\mathrm{i}\frac{t}{2}log\left(\frac{t}{2\pi N^2}\right)-\mathrm{i}t-\mathrm{i}\pi N^2\right)+O_1\left(\frac{2^{3/2}5\sqrt{\pi}(1/N+2\pi N/t)^2}{16(tN)^{3/2}(1/N-2\pi N/t)^3}\right).$$

(Here $f(x) = O_1(g(x)) \Rightarrow |f(x)| \leq 1 \times g(x)$ for $x \geq x_0$ as opposed to the usual $f(x) = O(g(x)) \Longrightarrow |f(x)| \leq C \times g(x)$ for $x \geq x_0$ for some absolute constant $C > 0$.)

Then for appropriate choices of $\eta$ and $\mu$ one can prove the following propositions:

a) If $N = 0$, or $2 \le N \le N_{2\eta}$, or $N > N_t$, then

$$D(N,R) = \left[H(N-1) - H\left(N - N_\mu\right)\right]\{Y(N,t)\} + O(\{t|1-N|\}^{-1}) + max\left[O\left(t^{-3/4}\right)\right].$$

b) If $N \in \left(N_{2\eta}, N_\mu\right]$, then $|D(N,R) - Y(N,t)| < 0.38 \times 3.55t^{-3/4}$.

c) If $N \in (N_\mu, N_\eta]$, then $|D(N,R)| < 0.64 \times 1.728\sqrt{2\pi/Nt} + 0.38 \times 3.55t^{-3/4}$.

d) If $N \in (N_\eta, N_t]$, then $|D(N,R)| < 0.38 \times 3.55t^{-3/4}$.

The *Proof of Lemma B1.2* is discussed in the next three subsections.

B1.3 *Contributions to* $D(N,R)$ *near the saddle point when* $N \in \left[2,\ N_\eta = INT\left(a/4\sqrt{pc_-}\right)\right]$.

From equation (A47), $pc(z)$ is defined by

$$\frac{\pi z^2}{4} = \frac{t(pc+1)^2}{2pc}. \qquad \text{(B13)}$$

Hence when $z = t(1 - 1/u_{sad})/\pi$ one has

$$\left(1 - \frac{1}{u_{sad}}\right) = \sqrt{\frac{2\pi}{t}}\left[\sqrt{pc_{sad}} + \frac{1}{\sqrt{pc_{sad}}}\right]. \qquad \text{(B14)}$$

As the saddle point is defined by $\sqrt{pc_{sad}} = a/4N$, from (B12), one can formulate the following *exact* results:

$$u_{sad} = \frac{1}{1 - 1/N - 2\pi N/t}, \qquad \text{(B15a)}$$

$$-\mathrm{i}t\left\{N\left(1 - \frac{1}{u_{sad}}\right) - \frac{\left(log(pc_{sad}) + \frac{1}{pc_{sad}}\right)}{2}\right\} = \mathrm{i}\frac{t}{2}log\left(\frac{t}{2\pi N^2}\right) - \mathrm{i}t - \mathrm{i}\pi N^2, \qquad \text{(B15b)}$$

$$\frac{-\mathrm{i}\delta^2 t}{2!\,u_{sad}^3}\left\{\frac{a}{2\sqrt{pc_{sad}}} - 2N + \frac{t/\pi u_{sad}}{2(pc_{sad}-1)}\right\} = \frac{-\mathrm{i}\delta^2 tN}{2u_{sad}^4[1/N - 2\pi N/t]}, \qquad \text{(B15c)}$$

$$\frac{1}{u_{sad}^2[(1 - 1/u_{sad})^2 - 8\pi/t]^{1/4}} = \frac{1}{u_{sad}^2[1/N - 2\pi N/t]^{1/2}}. \qquad \text{(B15d)}$$

A couple of remarks concerning equation (B15) are important. Firstly, from (B15a), one can see that when $N = 2$ the saddle point lies just above $u = 2$ and all the remaining saddles lie in the range $u \in (1, 2)$. Secondly, at a point $u_{sad} + \delta$ lying close to a saddle, the phase of the numerator in (B12) is, to $O(\delta^3)$, the sum of (B15b) and (B15c). Consequently if one arranges for the contour to pass through the saddle along a line defined by $u = u_{sad} + qe^{-\mathrm{i}\pi/4}$, $q \in \{-s, s\}$ (see Fig.B2) then the real part of the phase of the numerator will be negative either side of the saddle point (cf. equation B15c with $\delta = qe^{-\mathrm{i}\pi/4}$). Hence integral (B9) through the saddle can now be approximated by

$$\int_{u_{sad}-se^{-\frac{\mathrm{i}\pi}{4}}}^{u_{sad}+se^{-\frac{\mathrm{i}\pi}{4}}} \frac{exp\left[-\mathrm{i}tN\left(1-\frac{1}{u}\right)+\mathrm{i}(t/2)\left(log\big(pc(u)\big)+1/pc(u)\right)\right]}{u^2[(1-1/u)^2-8\pi/t]^{\frac{1}{4}}}du \approx$$

$$e^{-\mathrm{i}\pi/4}\, exp\left(\mathrm{i}\frac{t}{2}log\left(\frac{t}{2\pi N^2}\right)-\mathrm{i}t-\mathrm{i}\pi N^2\right)\int_{-s}^{s}\frac{e^{-q^2tN/2u_{sad}^4[1/N-2\pi N/t]}}{u_{sad}^2[1/N-2\pi N/t]^{1/2}}dq. \qquad \text{(B16)}$$

Unfortunately there is a complication inherent in this procedure. If $N$ happens to lie very close to $N_t$, then the saddle point starts to converge onto the singularity of the integrand at $u = 1/\left(1-2\sqrt{2\pi/t}\right)$ (see Fig. B3). In such a situation the integrand starts fluctuating very rapidly as one passes through the saddle point, not so much because of the exponential term but rather because the denominator is now very close to zero. This problem can be further highlighted by calculation of higher order corrections associated with the procedure (B13-16). If one examines the next term in the Taylor expansion of the denominator about the saddle point beyond (B15d) in more detail (actually the second order term because only even terms in integrand B16 across the saddle will contribute), routine calculation shows (see [2] or equations A49 & A51) that the magnitude of the largest contribution to the first order correction is given by

$$\left(\frac{2^{3/2}5\sqrt{\pi}(1/N+2\pi N/t)^2}{16(tN)^{3/2}(1/N-2\pi N/t)^3}\right), \qquad \text{(B17)}$$

which clearly has the potential to blow up as $N \to \sqrt{t/2\pi}$. In section B1.5 a specific definition for $N_\mu$ will be made that ensures (B17) remains 'small', in a well defined sense, provided $N < N_\mu$. For the moment it is sufficient to remark that it is possible to demonstrate that apart from a 'tiny fraction' of $N \in [N_\mu, N_\eta]$, the correction term (B17) remains extremely small in comparison with the main term obtained by integrating (B16). Provided the integration limit $s$ is chosen to be at least of $O(t^{-7/6})$ (specifically, at least as large as the value to be prescribed by equation B25), this integration gives

$$Y(N,t) = \sqrt{\frac{2\pi}{Nt}}e^{-\mathrm{i}\frac{\pi}{4}}exp\left(\mathrm{i}\frac{t}{2}log\left(\frac{t}{2\pi N^2}\right)-\mathrm{i}t-\mathrm{i}\pi N^2\right)+O_1\left(\frac{2^{3/2}5\sqrt{\pi}(1/N+2\pi N/t)^2}{16(tN)^{3/2}(1/N-2\pi N/t)^3}\right), \quad \text{(B18)}$$

as the contribution to (B9) in the vicinity of the saddle point.

B1.4 *Contributions to $D(N,R)$ away from the saddle point*

Additional contributions to $D(N,R)$ arise first from the integration along the line $u = t/\pi R + \mathrm{i}\delta$ for $N = 0$ and $N \in [2, N_\eta]$ and along the line $u = t/\pi R - \mathrm{i}\delta$ for $N > N_\eta$. In each of these cases $\delta$ is assumed to increase from zero until either it reaches the minimum of the real part of the numerator's phase (between $t^{-1}$ and $t^{-1/2}$ depending on $N$) or joins up with the extremity of the contour passing through the saddle point as in Fig. B2. In any event, the integrand will be at least $O\left(e^{-O(1)\sqrt{t}}\right)$ at this extremity and the integral is then dominated by its behaviour in very close proximity to $\delta = 0$. Using (B10) and the equivalent result to (B4), one obtains a contribution, analogous to (B5), of the form

$$\frac{-\mathrm{i}exp[-\mathrm{i}t\{N(1-\pi R/t) - (log(pc_-) + 1/pc_-)/2\}]}{(t/2)((1-\pi R/t)^2 - 8\pi/t)^{1/4}\left(2N - a/2\sqrt{pc_-}\right)}\left\{1 + \frac{\mathrm{i}(2\pi R/t + (\pi R)^2(1-\pi R/t)[(1-\pi R/t)^2 - 8\pi/t]^{-1}/2t^2)}{\left(2N - a/2\sqrt{pc_-}\right)(\pi R)^2/2t}\right\} + O\left(\frac{[(1-\pi R/t)^2 - 8\pi/t]^{-\frac{9}{4}}}{\left(t\left(2N - a/2\sqrt{pc_-}\right)\right)^3}\right). \quad \text{(B19)}$$

Now if $N \le N_{2\eta}$ or $N \ge N_t$ the first term in (B19) is *at most* of $O\left(\eta^{-5/2}t^{-3/4}/2\sqrt{\pi}\right)$ and dominates the remaining terms. In section B1.5 this 'dominance' is made explicit by a suitable choice of $\eta$, which will define $pc_-$ and $R_{max}$. However, as with the integration through the saddle point, a problem arises when $N$ happens to lie close to $N_t$. In such instances the $\left(2N - a/2\sqrt{pc_-}\right)$ factor gets relatively small, and the second term in (B19) has the potential to dominate the first term, because of the presence of the large reciprocal factor $[(1-\pi R/t)^2 - 8\pi/t]$ (a consequence of the lower limit of the integral $D(N,R)$ lying extremely close to unity and the singularity of the integrand). This situation, which again influences only a (marginally larger) 'tiny fraction' of $N \in [N_{2\eta}, N_t]$, is also addressed in section B1.5.

Finally there is the contribution, analogous to (B8), to the integral when accrues along the lines $u = t/\pi R + \mathrm{i}\delta$ for $N = 0$ and $u = t/\pi R - \mathrm{i}\delta$ for $N \ge 2$ when $\delta \to \infty$. In this region the denominator can be approximated by

$$-\frac{1}{\delta^2(1-8\pi/t)^{1/4}} \mp \frac{\mathrm{i}}{\delta^3(1-8\pi/t)^{1/4}}\left(\frac{2t}{\pi R} - \frac{1/2}{(1-8\pi/t)}\right)\dots, \quad \text{(B20)}$$

which together with (B11) gives rise to a contribution to (B9) of the form

$$\pm\mathrm{i}exp\left[-\mathrm{i}\left(tN - \frac{t}{2}log\left(\frac{t}{2\pi}\right) + \pi + \frac{2\pi^2}{t}\right)\right]\int_{\delta_L}^{\infty} e^{\pm\frac{1}{\delta}[t(N-1)-2\pi-8\pi/t]}\left\{-\frac{1}{\delta^2(1-8\pi/t)^{1/4}} \mp \frac{\mathrm{i}}{\delta^3(1-8\pi/t)^{1/4}}\left(\frac{2t}{\pi R} - \frac{1/2}{(1-8\pi/t)}\right)\right\}d\delta. \quad \text{(B21)}$$

Assuming that $1 \ll \delta_L \ll t$, then to leading order (B21) can be approximated by

$$\frac{-\mathrm{i}e^{-\mathrm{i}\left[tN-\frac{t}{2}log\left(\frac{t}{2\pi}\right)+\frac{2\pi^2}{t}\right]}}{(1-8\pi/t)^{1/4}[t(N-1)-2\pi-8\pi/t]}+O\left(\frac{1}{[t(N-1)-2\pi-8\pi/t]^2}\right). \qquad \text{(B22)}$$

B1.5 *Estimates for* $D(N,R)$ *when* $N \approx N_t$.

For (B18) to represent a good approximation to the integral across the saddle point, one requires that the relative error of the ratio of the magnitudes of the main term in (B18) and the largest part of the first order correction term (B17) should be sufficiently small. That is one requires

$$\frac{2^{3/2}5\sqrt{\pi}(1/N+2\pi N/t)^2\sqrt{Nt}}{16(tN)^{3/2}(1/N-2\pi N/t)^3\sqrt{2\pi}}<\varepsilon_R, \qquad \text{(B23)}$$

for some small value of $\varepsilon_R$. A short calculation reveals that (B23) will be satisfied when

$$N\le\left(\frac{t}{2\pi}\right)^{1/2}\left\{1-\left(\frac{5}{16\varepsilon_R t}\right)^{1/3}\right\}. \qquad \text{(B24)}$$

So for convenience if one prescribes the relative error to be less than 2% , or $\varepsilon_R = 0.02$, then this will be satisfied for all $N \le \sqrt{t/2\pi} - \mu t^{1/6}$ provided one defines $\mu = 5/\sqrt{8\pi} \approx 0.9973$. Hence for all but the 'tiny fraction' of the $N$ lying between $[N_\mu, N_\eta]$, amounting to some $t^{1/6}$ out of $\sqrt{t/2\pi}$ in total, the main term of (B18) will provide approximation to within at least 2% the integral across the saddle point. By contrast if $N \in [N_\mu, N_\eta]$, then (B18) will become increasingly inadequate and a revised estimate of the integral across the saddle is required.

Let $r_N \in (0,1)$ so that $N = \sqrt{t/2\pi} - r_N t^{1/6}$. Investigation along the line $u = u_{sad} + qe^{-\mathrm{i}\pi/4}$ reveals that the denominator of the integral (B16) reaches a minimum at $q_{min} = -2^{-3/2}u_{sad}^2(1/N-2\pi N/t)^2(1/N+2\pi N/t)^{-1}$, (which approximates to about $-11.13r_N^2 t^{-7/6}$) and the reciprocal of the denominator at this point amounts to a relatively small increase of $1.0905077e^{-\mathrm{i}\pi/16}$ times its (real) value at $u = u_{sad}$ itself. Looking at the real part of the exponential phase in integral (B16), one can see that the exponential itself will become small when $q$ lies outside the range

$$\pm f\left\{\frac{2u_{sad}^4[1/N-2\pi N/t]}{tN}\right\}^{1/2}=\pm fX\approx\pm f\frac{(128\pi^3)^{1/4}\sqrt{r_N}}{t^{7/6}}\approx\pm f\frac{7.93\sqrt{r_N}}{t^{7/6}}. \qquad \text{(B25)}$$

Here $f$ is a $O(1)$ factor which should be prescribed so that the decline in the exponential term to the level of $e^{-f^2}$, combined with the rapid increase in the denominator away from $q_{min}$, will mean the integrand is now sufficiently small to kill off any significant contribution to the integral outside the range (B25). A value of $f = 1.5$ is perfectly adequate for this purpose. As $r_N$ decreases from about unity, $q_{min}$ moves increasingly to the centre of the range defined by (B25), and the integrand changes character from a relatively rapidly changing exponential term times an almost constant denominator, to a rapidly changing denominator times a relatively slowly varying exponential term. Now over the small range defined by (B25), the imaginary part of the numerator phase (B10) varies only slightly from (B15b) by a factor $\varphi(q)$ (with a maximum $|\varphi(q)| \le 9/8$ when $r_N = 0.261$ at $q = -1.5X$). So the integral through the saddle can now be written as

$$\int_{u_{sad}-fXe^{-\frac{i\pi}{4}}}^{u_{sad}+fXe^{-\frac{i\pi}{4}}} \frac{exp\left[-itN\left(1-\frac{1}{u}\right)+i(t/2)\left(log\left(pc(u)\right)+1/pc(u)\right)\right]}{u^2[(1-1/u)^2-8\pi/t]^{\frac{1}{4}}}du$$

$$= e^{-i\pi/4}\,exp\left(i\frac{t}{2}log\left(\frac{t}{2\pi N^2}\right)-it-i\pi N^2\right)\int_{-fX}^{fX}\frac{\left(u_{sad}+qe^{-i\pi/4}\right)^{-2}e^{-q^2/X^2+i\varphi(q)}}{[(1-1/\{u_{sad}+qe^{-i\pi/4}\})^2-8\pi/t]^{1/4}}dq,$$

(B26)

from which an upper bound on the modulus can be established by setting $e^{-q^2/X^2+i\varphi(q)} = 1$. Applying the substitution $1-1/\{u_{sad}+qe^{-i\pi/4}\} = \sqrt{8\pi(1+v^4)/t}$, the integral in (B26) can be transformed to

$$e^{i\pi/4}\left(\frac{8\pi}{t}\right)^{\frac{1}{4}}\int\frac{2v^2}{\sqrt{(1+v^4)}}dv \approx e^{i\pi/4}\left(\frac{8\pi}{t}\right)^{\frac{1}{4}}\int 2\,v^2dv = e^{i\pi/4}\left(\frac{8\pi}{t}\right)^{1/4}\frac{2}{3}v^3 \qquad \text{(B27)}$$

The latter approximation can be justified because $|v^4| \ll 1$ across the entire range of (B25). For example a rough estimate (with $f = 1.5$) shows that

$$|v^4| = \left|\left[t\{1-1/\{u_{sad}+qe^{-i\pi/4}\}\}^2/8\pi-1\right]\right| \approx \left|\sqrt{\frac{t}{2\pi}}qe^{-i\pi/4}\right| < 5\sqrt{r_N}t^{-2/3}. \qquad \text{(B28)}$$

Replacing $v$ in (B27) with the original substitution, the integral (B26) is bounded above by

$$e^{i\pi/4}\frac{2}{3}\left(\frac{t}{8\pi}\right)^{\frac{1}{2}}\left\{\left[\left(1-\frac{1}{u_{sad}+fXe^{-i\pi/4}}\right)^2-\frac{8\pi}{t}\right]^{3/4}-\left[\left(1-\frac{1}{u_{sad}-fXe^{-i\pi/4}}\right)^2-\frac{8\pi}{t}\right]^{3/4}\right\}. \text{(B29)}$$

This can be rewritten in a more illuminating form by noting that

$$\left(1-\frac{1}{u_{sad}\pm fXe^{-i\pi/4}}\right)^2 \approx \left(1-\frac{1}{u_{sad}}\right)^2 \pm \frac{2}{u_{sad}^2}\left(1-\frac{1}{u_{sad}}\right)fXe^{-i\pi/4}$$

$$= \left(\frac{1}{N}+\frac{2\pi N}{t}\right)^2 \pm 2\sqrt{2}\left(\frac{1}{N}+\frac{2\pi N}{t}\right)\left\{\frac{[1/N-2\pi N/t]}{tN}\right\}^{\frac{1}{2}}\left(\frac{tN}{tN}\right)^{\frac{1}{6}}fe^{-i\pi/4}$$

$$\approx \frac{8\pi}{t}+\frac{16\pi^2r_N^2}{t^{5/3}} \pm \frac{(2^{13}\pi)^{1/4}r_N^{1/2}}{t^{2/3}}\left(\frac{2\pi}{tN}\right)^{2/3}fe^{-i\pi/4}+O(t^{-2}), \qquad \text{(B30)}$$

if $N = \sqrt{t/2\pi} - r_Nt^{1/6}$. Utilising this result, (B29) approximates to

$$\frac{e^{\mathrm{i}\frac{\pi}{4}}}{3}\sqrt{\frac{t}{2\pi}}\left\{\left[\frac{16\pi^2 r_N^2}{t^{5/3}}+\frac{(2^{13}\pi)^{1/4}r_N^{1/2}}{t^{2/3}}\left(\frac{2\pi}{tN}\right)^{\frac{2}{3}}fe^{-\mathrm{i}\frac{\pi}{4}}\right]^{\frac{3}{4}}-\left[\frac{16\pi^2 r_N^2}{t^{5/3}}-\frac{(2^{13}\pi)^{1/4}r_N^{1/2}}{t^{2/3}}\left(\frac{2\pi}{tN}\right)^{\frac{2}{3}}fe^{-\mathrm{i}\frac{\pi}{4}}\right]^{\frac{3}{4}}\right\}$$

$$=\frac{e^{\mathrm{i}\frac{\pi}{4}}}{3}\left(\frac{2^{31}}{\pi^5}\right)^{\frac{1}{16}}\sqrt{\frac{2\pi}{tN}}r_N^{3/8}f^{\frac{3}{4}}\left\{\left[\left(\frac{N^2\pi^{13/4}}{2^{1/4}t}\right)^{\frac{1}{3}}\frac{r_N^{3/2}}{f}+e^{-\mathrm{i}\frac{\pi}{4}}\right]^{3/4}-\left[\left(\frac{N^2\pi^{13/4}}{2^{1/4}t}\right)^{\frac{1}{3}}\frac{r_N^{3/2}}{f}-e^{-\mathrm{i}\frac{\pi}{4}}\right]^{3/4}\right\}.$$

(B31)

Multiplication of (B31) by the constant exponential terms outside integral (B26) gives an upper bound on the integral through the saddle point for integer values $N=\sqrt{t/2\pi}-r_N t^{1/6}$. (Numerical estimates suggest that the effect of the neglected term $e^{-q^2/X^2+\mathrm{i}\varphi(q)}$ is to multiply (B31) by a factor of 0.64 or so, close to the value of 0.57 for the average value of $e^{-q^2/X^2}$ integrated across the range (B25) with $f=1.5$.) When $r_N=1$ the integral across the saddle point is (as to be expected) of the same order of magnitude as the previous estimate (B18) (at $r_N=\mu\approx 1$, B31 equates to a *maximum* value of $1.728\sqrt{2\pi/Nt}\,e^{-\mathrm{i}0.034}$), but then crucially (B31) dies away as $r_N\to 0$. This analysis confirms that the increase in the error term (B17) as $N$ approaches $\sqrt{t/2\pi}$ is indicative of the changing nature of the integrand rather than an indication the integral itself is actually diverging. Now if $R\approx_p R_{max}$ as defined in *Lemma B1.2*, then for $N\in[2,N_\eta]$ the corresponding saddle points given by $\sqrt{pc_{sad}}=a/4N$ all lie within the integration range of integral $D(N,R)$. Provided $\mu=5/\sqrt{8\pi}>\eta$, the contribution to $D(N,R)$ across each saddle will be given, to within 2%, by the main term of (B18) if $N\le N_\mu$, and if $N\in(N_\mu,N_\eta]$ will be bounded above by the new estimate (B31) which has a maximum value of $0.64\times 1.728\sqrt{2\pi/Nt}$.

In a similar vein, the main term in (B19) will only represent a good approximation of the integral along the lines $u=t/\pi R\pm\mathrm{i}\delta$ provided the first order correction is sufficiently small. Hence one requires a similar condition to (B23) regarding the relative magnitudes of these terms, to determine for which values of $N$ (B19) ceases to be a good approximation. This condition can be formulated as

$$\left|\frac{(1-\pi R/t)}{t\left(2N-a/2\sqrt{pc_-}\right)[(1-\pi R/t)^2-8\pi/t]}\right|<\varepsilon_R. \qquad \text{(B32)}$$

Equation (B32) is a somewhat more complicated constraint than (B23), as it depends not just on $N$, but also on $R$, the size of which determines how close the lower limit of integration lies to the singularity at $u=1/\left(1-2\sqrt{2\pi/t}\right)$. Suppose $R=t/\pi-a-\rho$ with $\rho\ll 1$, which implies that $pc_-\approx 1+\sqrt{\pi a\rho/t}$. Assuming $N=\sqrt{t/2\pi}-r_N t^{1/6}$ as before (although in this instance $r_N$ can potentially be negative if modifications to the integration along $u=t/\pi R-\mathrm{i}\delta$ when $N$ is slightly above $\sqrt{t/2\pi}$ prove necessary), a short calculation transforms (B32) to

$$\frac{1}{2\pi\rho\left|\sqrt{a\rho/2}-2r_N t^{1/6}\right|}<\varepsilon_R. \tag{B33}$$

Now for a single, *independent* $N<\sqrt{t/2\pi}$ value, one can consider the question as to the most appropriate choice of $R$, or equivalently $pc_-$, that guarantees (B32) remains less than $\varepsilon_R$. The saddle point associated with $N$ satisfies $\sqrt{pc_{sad}}=a/4N \Rightarrow pc_{sad}\approx 1+r_N\sqrt{8\pi}t^{-1/3}$. Examination of (B32) shows that if $pc_-\to pc_{sad}$ then first term in the denominator approaches zero, whilst if $pc_-\to 1$ (below $pc_{sad}$) then the second term in the denominator approaches zero. This suggests fixing $pc_-=1+r_N\sqrt{2\pi}t^{-1/3}$, the half way point between the two singular values of (B32) at $pc_-=pc_{sad}$ and $pc_-=1$. This choice of $pc_-$ is equivalent to setting $\rho=r_N^2\sqrt{\pi/2}\,t^{-1/6}$, in which case (B33) becomes

$$\frac{1}{\sqrt{2\pi^3}r_N^3}<\varepsilon_R. \tag{B34}$$

If $\varepsilon_R$ is set to be 2% as above, then $r_N=(1250/\pi^3)^{1/6}=2\eta\approx 1.852\Rightarrow \eta\approx 0.9259<\mu$. Defining $\eta$ and $pc_-=1+\eta\sqrt{8\pi}t^{-1/3}$ this way, ensures that for those $N<N_{2\eta}$ or $N>N_t$ the first order correction to (B19) will amount to no more than 2% of the main term. However, for $N\in[N_{2\eta},\,N_t]$ the factor $t\left(2N-a/2\sqrt{pc_-}\right)$ is insufficiently large to counteract the smallness of the $[(1-\pi R/t)^2-8\pi/t]$ factor arising from this choice of $pc_-$. Consequently for $N$ in this range the main term of (B19) will be an increasingly inadequate approximation to the integral along $u=t/\pi R\pm \mathrm{i}\delta$. Hence in these instances a better estimate (or at least an upper bound) must be substituted. (NB. For those $N<N_\eta$, $pc_{sad}>pc_-$ and $u_{sad}>t/\pi R$, the prescribed contour must pass through the saddle as in Fig. B3a; but for those satisfying $N_\eta<N<N_t$ their respective saddle points now lie between $u=1/\left(1-2\sqrt{2\pi/t}\right)$ and $t/\pi R$, and do not form part of the range of integration of $D(N,R)$, see Fig. B3b.)

To obtain such an upper bound it is necessary to modify the orientation of the contour emanating from $u=t/\pi R$ to a more general direction $u=t/\pi R+e^{\mathrm{i}\sigma}\delta$ (with the choice of $\sigma$ tailored to the specific value of $r_N$) to some appropriate endpoint $\delta_{end}$. Often this will be where it intersects the contour branch passing through the saddle point. It is possible to define four suitable choices of $\sigma$ and $\delta_{end}$ accounting for all the problem range $r_N\in[0,2\eta]$, viz.

$$\delta_{end}=\begin{cases} f^2\left[2t/(\pi R)^2\left(2N-a/2\sqrt{pc_-}\right)\right]\approx f^2(\eta-r_N)^{-1}t^{-7/6}, & 0<r_N\le 0.7,\ \ \sigma=-\pi/2,\\ fX, & 0.7<r_N<\eta,\ \ \sigma=-\pi/4,\\ \sqrt{f^2X^2+(u_{sad}-t/\pi R)^2-fX\sqrt{2}(u_{sad}-t/\pi R)}, & r_N\in[\eta,\ 1.2],\ \ \sigma=\pi-arcsin\left(\dfrac{fX}{\delta_{end}\sqrt{2}}\right),\\ 2.53\left[2t/(\pi R)^2\left(a/2\sqrt{pc_-}-2N\right)\right]\approx 2.53(r_N-\eta)^{-1}t^{-7/6}, & r_N\in(1.2,\ 2\eta),\ \ \sigma=\pi/2. \end{cases}$$

(B35a,b,c,d)

with $f$ and $X$ as given by (B25). The contour branch (B35a) for $r_N \in (0, 0.7]$ is appropriate in those instances when $u_{sad} < t/\pi R$ and the saddle point does not lie in the integration range of $D(N, R)$. The choice of $r_N \leq 0.7$ ensures $\delta_{end}$ does not exceed $fX$. The second branch (B35) is motivated by the fact the non contributing saddles are now almost coincident with $u = t/\pi R$ (see Fig. B3b). Hence defining $\sigma = -\pi/4$ in these instances will ensure the integrand decays very rapidly as $\delta$ increases to $\delta_{end}$, in an analogous fashion to (B26). The third branch (B35c) is defined so that it joins up with the branch through the saddle at the point $u_{sad} - fXe^{-i\pi/4}$ (see Fig. B3a). Crucially the choice of $\delta_{end}$ in (B35c) can never exceed $fX$ and so like the other options is of order $t^{-7/6}$. Finally the choice of branch (B35d) is motivated by the fact that the contributing saddles lie far enough away from $u = t/\pi R$ to integrate at right angles to the real axis, without truncating the branch through the saddle significantly. (The cut off at $r_N \approx 1.2$ in equations B35c,d is where the respective values of $\delta_{end}$ become equal.) These four choices of $\delta_{end}$ and $\sigma$ are sufficient to ensure the contour branch reaches a region where the real part of the phase in (B10) is at least as small as $e^{-f^2} = 0.105$ in magnitude. Hence any further integration beyond $\delta_{end}$, which may be necessary to form a complete contour path from $u = t/\pi R$ to $u = \infty$, will only add an exponentially small contribution.

An upper bound on the integration along the contour branches prescribed by (B35) can be found by noting that when $N \in (N_{2\eta}, N_t]$ the most rapid changes to the integrand arise from the behaviour of the denominator in the vicinity of its root, and that any fluctuations in the exponential phase are, by comparison, much slower. So to establish an upper bound the latter will be assumed to take its maximum value along the entire contour branch, which is specified at $u = t/\pi R$ . With this assumption an upper bound to $D(N, R)$ along the contours prescribed by (B35), can be derived following an analogous procedure to that specified by equations (B26-31) for branches passing through saddle points. The main result is (cf. equation B29)

$$\left\{\left[\left(1 - \frac{1}{t/\pi R + \delta_{end} e^{i\sigma}}\right)^2 - \frac{8\pi}{t}\right]^{\frac{3}{4}} - \left[\left(1 - \frac{1}{t/\pi R}\right)^2 - \frac{8\pi}{t}\right]^{\frac{3}{4}}\right\} \frac{2\sqrt{t} e^{\left[-it\left\{N\left(1-\frac{\pi R}{t}\right) - \left(log(pc_-) + \frac{1}{pc_-}\right)/2\right\}\right]}}{3\sqrt{8\pi}}$$

(B36)

Using the analogous series to (B30) and some algebra (B36) can be approximated to

$$\approx \frac{2\sqrt{8}\pi\eta^{3/2}}{3t^{3/4}} \left\{\left[1 + \frac{\left(t^{7/6}\delta_{end}\right)e^{i\sigma}}{\eta^2 \pi^{3/2}\sqrt{8}}\right]^{3/4} - 1\right\} e^{\left[-it\left\{N\left(1-\frac{\pi R}{t}\right) - \left(log(pc_-) + \frac{1}{pc_-}\right)/2\right\}\right]}. \qquad \text{(B37)}$$

Equation (B37) provides an upper bound (because the exponential decay has been neglected) for the contour integrals along (B35) which should be substituted for (B19) when ever $N \in [N_{2\eta}, N_t]$.The modulus of (B37) at $r_N = 2\eta$ is about $0.8/t^{3/4}$, compared to the value of $\approx 1/\pi t^{3/4}$ for (B19) at this value of $r_N$. The discrepancy between these two values, $1/0.8\pi \approx 0.4$, approximately corresponds to $\left[1 - exp\left(-t^{7/6}\delta_{end}\right)\right]/t^{7/6}\delta_{end} \approx 0.34$ the average value of $e^{-q}$ over the range $q \in \left[0, t^{7/6}\delta_{end} = 2.53/\eta\right]$, which is as one might expect. So in general, the factor $\left[1 - exp\left(-t^{7/6}\delta_{end}\right)\right]/t^{7/6}\delta_{end}$ should give a rough idea of the effect of the neglected exponential decay. This factor has a maximum value of about 0.38 over the ranges prescribed in (B35) when $r_N = 0$ and $t^{7/6}\,\delta_{end} = f^2/\eta$. The maximum modulus of (B37) is $3.55/t^{3/4}$ when $r_N = \eta$ and $\sigma \approx 3\pi/4$, but this corresponds to integrating over a range of $[0, fX]$, four times larger than the $\left[0, f^2 t^{-7/6}/\eta\ \right]$ prescribed for $r_N = 0$. So when $r_N = \eta$ the neglected exponential decay would result in a much greater reduction than the maximum 0.38 factor estimated above. Consequently it is safe to conclude that none of the integrals along the contours defined by (B35) can ever exceed a value given by the product of these two respective maxima, namely $0.38 \times 3.55/t^{3/4}$. This provides an upper bound on the contribution to $D(N, R)$ when integrating along any of the various branches defined by (B35), which can be substituted instead of (B19) when $N \in (N_{2\eta},\ N_t]$.

Summarising all the results derived in the last three sections, one has the following:

Define $O(1)$ constants $\eta = (1250/\pi^3)^{1/6}/2$, and $\mu = 5/\sqrt{8\pi}$ which satisfy the conditions set out in *Lemma B1.2*.

i) If $N = 0$, or $2 \le N \le N_{2\eta}$, or $N > N_t$, then $D(N, R)$ is given by the sum of terms (B18), (B19) and (B22), plus an error that is guaranteed to be less than 2% of the main terms of (B18) and (B19). Hence part a) of the lemma is true.

ii) If, $N \in (N_{2\eta}, N_\mu]$ then (B18) and (B22) are still valid, but (B19) must be replaced by the upper bound $0.38 \times 3.55/t^{3/4}$. Hence part b) of the lemma is true.

iii) If $N \in (N_\mu, N_\eta]$ then neither (B18) nor (B19) are valid, but they can be replaced by their upper bounds, $0.64 \times 1.728\sqrt{2\pi/Nt}$ and $0.38 \times 3.55/t^{3/4}$ respectively. Hence part c) of the lemma is true.

iv) If $N \in (N_\eta, N_t]$ then the integration range of $D(N, R)$ does not include a saddle point contribution. Equation (B19) is not valid for $N$ in this range, but the upper bound $0.38 \times 3.55/t^{3/4}$ can be substituted. Hence part d) of the lemma is true.

Hence all four conditions a)-d) specified in *Lemma B1.2* regarding the integral $D(N, R)$ are true and the lemma is proved.

B1.6 *Evaluation of $C(N, R)$ and $D(N, R)$ for the case when $N = 1$.*

*Lemma B1.6* Let $C(1, R)$ and $D(1, R)$ be the integrals defined by (B1) and (B9) respectively and $R \approx_p R_{max}$ as defined in *Lemma B1.2* above. Then

$$C(1, R) + D(1, R) = \frac{\sqrt{2\pi}e^{-\mathrm{i}\left(t - \frac{t}{2}log\left(\frac{t}{2\pi}\right) + \pi + \frac{2\pi^2}{t}\right) - \mathrm{i}\pi/4}}{\sqrt{t}(1 - 8\pi/t)^{\frac{1}{4}}} + O(t^{-1}) + O\left(t^{-13/12}\right) + O\left(t^{-\frac{2}{3}}\right). \qquad \text{(B38)}$$

*Proof.* The specific problem with the estimation of $C(N,R)$ and $D(N,R)$ when $N=1$ concerns the behaviour of the real part of the numerator's phase along the line $u = t/\pi R + \mathrm{i}\delta$ as $\delta \to \infty$. As one can see from both (B3) and (B11), if $N=1$ it is no longer clear as to the exact behaviour of this term because the previously dominant $t(N-1)$ terms vanish. The real parts of the corresponding asymptotic expansions (B3) and (B11) with $N=1$ (with $t/\pi R = 1$ rather than $\approx 1$ to simplify the algebra somewhat) reduce to

$$\frac{2\pi + 8\pi^2/t}{\delta} - \frac{(4t/3 + 12\pi + 88\pi^2/t)}{\delta^3} + \frac{(26t/5 + 64\pi + 720\pi^2/t)}{\delta^5} + O\left(\frac{t}{\delta^7}\right), \quad \text{(B39)}$$

$$-\frac{(2\pi + 8\pi^2/t)}{\delta} - \frac{(2t/3 - 8\pi^2/t)}{\delta^3} + \frac{(4t/5)}{\delta^5} + O\left(\frac{t}{\delta^7}\right), \quad \text{(B40)}$$

respectively. From (B39), one can see that when $\delta = O(1)\sqrt{t}$, the terms $2\pi/\delta$ and $-4t/3\delta^3$ are comparable in size. More detailed investigation reveals that at $\delta \approx \sqrt{2t/3\pi}$ the real part of the phase changes sign from negative to positive, then increases to a maximum of approximately $2\sqrt{2}\pi^{3/2}/3\sqrt{t} \approx 5.25/\sqrt{t}$ at $\delta \approx \sqrt{2t/\pi}$, before slowly declining to zero as $\delta \longrightarrow \infty$. This is illustrated by the dashed line in Fig. B4, which also shows the behaviour for small $\delta$ derived from equation (B2). The real part of the phase initially falls away very rapidly to a negative minimum of $O(t)$ at $\delta \approx 1/\sqrt{2}$, before slowly climbing back to its maximum of $5.25/\sqrt{t}$ as described. By contrast for (B40), the terms $-2\pi/\delta$ and $-2t/3\delta^3$ are both negative and consequently the real part of the phase is *always* negative in this instance. From this analysis one concludes that when $\delta$ is small there will be contributions to $C(1,R)$ and $D(1,R)$ identical in form to equations (B5) and (B19) for the $N \neq 1$ cases. However, the contribution along the line $u = t/\pi R + \mathrm{i}\delta$ for large $\delta$ will be very different in character to terms (B8) and (B22), because over the ranges $\delta \in \left[(2 \text{ or } 4t/3)^{1/3}, \infty\right)$ the magnitude of the real part of the phase remains very close to zero. Consequently the range of integration over which the integrand is of $O(1)$ is very much greater than when $N \neq 1$.

As the two integrals $C(1,R)$ and $D(1,R)$ are very similar, they will be analysed together. The discussion surrounding (B39-40) highlights the fact that one requires an estimate of

$$\int_{\frac{t}{\pi R}+\mathrm{i}\beta\sqrt{\frac{2t}{\pi}}}^{\frac{t}{\pi R}+\mathrm{i}\infty} \frac{exp\left[-\mathrm{i}t(1 \pm 1/u) + \mathrm{i}(t/2)\left(log\left(pc_{\pm}(u)\right) + 1/pc_{\pm}(u)\right)\right]}{u^2[(1 \pm 1/u)^2 - 8\pi/t]^{1/4}} du, \quad \text{(B41)}$$

where $pc_{\pm}(u) = pc(t(1 \pm u)/\pi)$ for $C(1,R)$ and $D(1,R)$ respectively, and $\beta$ is constant which must be fixed to an appropriate value. Obviously from the previous discussion it is clearly desirable that $\beta < 1/\sqrt{3}$ to ensure that the value of the real part of the numerator's phase is negative, but it is also necessary to establish a lower bound too. Using (B3, 11) one can write the phase of the numerator in (B41) as the following asymptotic expansion (both real and imaginary parts)

$$-\,\mathrm{i}\left(t-\frac{t}{2}log\left(\frac{t}{2\pi}\right)+\pi+\frac{2\pi^2}{t}\dots\right)+\frac{A_\pm}{\delta}+\frac{\mathrm{i}B_\pm}{\delta^2}+\frac{C_\pm}{\delta^3}+O\left(\frac{t}{\delta^4}\right), \qquad \text{(B42a)}$$

with

$$A_\pm=\pm\left(2\pi+\frac{8\pi^2}{t}\right),\ B_\pm=\frac{t}{2}+\pi\left[3\pm2\left(\frac{t}{\pi R}\right)\right]+\frac{\pi^2}{t}\left[20\pm8\left(\frac{t}{\pi R}\right)\right],\text{ and}$$

$$C_\pm=t\left[\mp\frac{1}{3}-\left(\frac{t}{\pi R}\right)\right]\mp2\pi\left[\left(\frac{t}{\pi R}\right)^2\pm3\left(\frac{t}{\pi R}\right)+2\right]\mp\frac{8\pi^2}{t}\left[\left(\frac{t}{\pi R}\right)^2\pm5\left(\frac{t}{\pi R}\right)+5\right]. \qquad \text{(B42b)}$$

Using (B6, 20) and (B42) integrals (B41) can be transformed into

$$\frac{-\mathrm{i}e^{-\mathrm{i}\left(t-\frac{t}{2}log\left(\frac{t}{2\pi}\right)+\pi+\frac{2\pi^2}{t}\right)}}{(1-8\pi/t)^{1/4}}\int_{\frac{\beta\sqrt{2t}}{\sqrt{\pi}}}^{\infty}\left(\frac{1}{\delta^2}+\frac{\mathrm{i}}{\delta^3}\left(\frac{2t}{\pi R}\pm\frac{1/2}{(1-8\pi/t)}\right)\right)exp\left[\frac{A_\pm}{\delta}+\frac{\mathrm{i}B_\pm}{\delta^2}+\frac{C_\pm}{\delta^3}\right]d\delta. \qquad \text{(B43)}$$

Now over this range of integration the real part of the phase, predominantly comprising the $A_\pm/\delta$ and $C_\pm/\delta^3$ terms, never exceeds a value of $O\left(t^{-1/2}\right)$ as discussed above. By contrast the term $\mathrm{i}B_\pm/\delta^2$ varies much more rapidly and can reach a value of $O(1)$ because $B_\pm\sim t/2$. So it is appropriate to approximate the exponential term in (B43) by $e^{\mathrm{i}B_\pm/\delta^2}\left(1+A_\pm/\delta+C_\pm/\delta^3\right)$, giving

$$\frac{-\mathrm{i}e^{-\mathrm{i}\left(t-\frac{t}{2}log\left(\frac{t}{2\pi}\right)+\pi+\frac{2\pi^2}{t}\right)}}{(1-8\pi/t)^{1/4}}\int_{\frac{\beta\sqrt{2t}}{\sqrt{\pi}}}^{\infty}\frac{e^{\mathrm{i}B_\pm/\delta^2}}{\delta^2}+\frac{e^{\mathrm{i}B_\pm/\delta^2}}{\delta^2}\left(\frac{A_\pm}{\delta}+\frac{\mathrm{i}}{\delta}\left(\frac{2t}{\pi R}\pm\frac{1/2}{(1-8\pi/t)}\right)+\frac{C_\pm}{\delta^3}\right)d\delta. \qquad \text{(B44)}$$

Concentrating on the first term in (B44) and making the substitution $1/\delta=v\sqrt{\pi/2B_\pm}$ one obtains

$$\frac{-\mathrm{i}e^{-\mathrm{i}\left(t-\frac{t}{2}log\left(\frac{t}{2\pi}\right)+\pi+\frac{2\pi^2}{t}\right)}}{(1-8\pi/t)^{1/4}}\sqrt{\frac{\pi}{2B_\pm}}\int_0^{\frac{1}{\beta}\sqrt{2B_\pm/2t}}e^{\mathrm{i}\frac{\pi}{2}v^2}\,dv. \qquad \text{(B45)}$$

The latter integral can be written in terms of the error function to give

$$\frac{-\mathrm{i}e^{-\mathrm{i}\left(t-\frac{t}{2}log\left(\frac{t}{2\pi}\right)+\pi+\frac{2\pi^2}{t}\right)}}{(1-8\pi/t)^{\frac{1}{4}}}\sqrt{\frac{\pi}{2B_\pm}}\frac{e^{\mathrm{i}\pi/4}}{\sqrt{2}}erf\left\{\frac{1}{\beta}\sqrt{\frac{\pi B_\pm}{2t}}e^{-\mathrm{i}\pi/4}\right\}$$

$$\approx\frac{e^{-\mathrm{i}\pi/4}e^{-\mathrm{i}\left(t-\frac{t}{2}log\left(\frac{t}{2\pi}\right)+\pi+\frac{2\pi^2}{t}\right)}}{(1-8\pi/t)^{\frac{1}{4}}}\sqrt{\frac{\pi}{2t\left(1+\frac{2\pi}{t}\left\{3\pm2\left(\frac{t}{\pi R}\right)\right\}\right)}}erf\left\{\frac{1}{2\beta}\sqrt{\pi\left(1+\frac{2\pi}{t}\left\{3\pm2\left(\frac{t}{\pi R}\right)\right\}\right)}e^{-\mathrm{i}\pi/4}\right\}.$$

(B46)

Notice that if the parameter $\beta$ is *small enough* so that the error function is approximately unity, then the sum of the $\pm$ terms given by (B46) is, to first order, just the saddle point contribution (B18) with $N = 1$, although it arises quite differently here. One can determine an appropriate size for $\beta$ from the terms previously neglected in (B44). Integration of these terms gives

$$\frac{-\mathrm{i}e^{-\mathrm{i}\left(t-\frac{t}{2}log\left(\frac{t}{2\pi}\right)+\pi+\frac{2\pi^2}{t}\right)}}{(1-8\pi/t)^{\frac{1}{4}}}\left\{\left(e^{\mathrm{i}B_\pm\pi/2\beta^2 t}-1\right)\left[\frac{C_\pm}{2B_\pm^2}+\frac{A_\pm+\mathrm{i}\left(\frac{2t}{\pi R}\pm\frac{1/2}{(1-8\pi/t)}\right)}{2\mathrm{i}B_\pm}\right]+\frac{C_\pm \pi e^{\mathrm{i}B_\pm\pi/2\beta^2 t}}{4\mathrm{i}B_\pm t\beta^2}\right\}$$

$$\approx \frac{-\mathrm{i}e^{-\mathrm{i}\left(t-\frac{t}{2}log\left(\frac{t}{2\pi}\right)+\pi+\frac{2\pi^2}{t}\right)}}{(1-8\pi/t)^{\frac{1}{4}}}\left(\frac{1}{\mathrm{i}t}\right)\left\{\pm\left(e^{\mathrm{i}B_\pm\pi/2\beta^2 t}-1\right)\left(2\pi-\frac{\mathrm{i}}{6}\right)-\frac{\left((t/\pi R)\pm 1/3\right)\pi}{2\beta^2}e^{\mathrm{i}B_\pm\pi/2\beta^2 t}\right\}$$
$$+\,O\left(\frac{1}{(\beta t)^2}\right). \tag{B47}$$

Now the first term in (B47) is clearly of $O(1/t)$ and will be dominated by the second term of $O(\pi/2t\beta^2)$ when $\beta$ is small. Hence for the magnitude of (B47) to be less than the magnitude of (B46), one requires that

$$\frac{[(t/\pi R)\pm 1/3]\pi}{2t\beta^2}\approx\frac{(1\pm 1/3)\pi}{2t\beta^2}<\sqrt{\frac{\pi}{2t}}\Rightarrow \beta>(1\pm 1/3)^{1/2}\left(\frac{\pi}{2t}\right)^{1/4}. \tag{B48}$$

This lower bound is indicative of the point at which the asymptotic analysis utilised in expressions (B42-44) starts to break down. Consequently $\beta$ should be chosen to lie within the range $\left((\pi/2t)^{1/4}, 1/\sqrt{3}\right)$. However, it is possible to define $\beta_O$ an optimal choice of $\beta$ by considering the second order terms in (B46). Using the asymptotic expansion for the error function $erf(z)\approx 1-e^{-z^2}/z\sqrt{\pi}$, the sum of second order $\pm$ terms in (B46) is given by

$$\frac{2e^{-\mathrm{i}\pi/4}e^{-\mathrm{i}\left(t-\frac{t}{2}log\left(\frac{t}{2\pi}\right)+\pi+\frac{2\pi^2}{t}\right)}}{(1-8\pi/t)^{\frac{1}{4}}}\sqrt{\frac{\pi}{2t}}\left\{-\frac{e^{\mathrm{i}\pi/4\beta^2}2\beta}{\pi e^{-\mathrm{i}\pi/4}}\right\}=-\frac{e^{-\mathrm{i}\left(t-\frac{t}{2}log\left(\frac{t}{2\pi}\right)+\pi+\frac{2\pi^2}{t}\right)}e^{\mathrm{i}\pi/4\beta^2}}{(1-8\pi/t)^{\frac{1}{4}}}\beta\sqrt{\frac{8}{\pi t}}, \tag{B49a}$$

whilst the sum of the two $\pm$ terms in (B47) (setting $B_+\approx B_-\approx t/2$) takes the form

$$\frac{e^{-\mathrm{i}\left(t-\frac{t}{2}log\left(\frac{t}{2\pi}\right)+\pi+\frac{2\pi^2}{t}\right)}}{(1-8\pi/t)^{\frac{1}{4}}}\left(\frac{2(t/\pi R)}{t}\right)\frac{\pi}{2\beta^2}e^{\mathrm{i}\pi/4\beta^2}\approx\frac{e^{-\mathrm{i}\left(t-\frac{t}{2}log\left(\frac{t}{2\pi}\right)+\pi+\frac{2\pi^2}{t}\right)}e^{\mathrm{i}\pi/4\beta^2}}{(1-8\pi/t)^{\frac{1}{4}}}\frac{\pi(t/\pi R)}{t\beta^2}. \tag{B49b}$$

Consequently terms (B49a) and (B49b) will cancel provided $\beta_O=(t\pi/8R^2)^{1/6}=O\left(t^{-1/6}\right)$. This means that, with $\beta=\beta_O$ and the results (B45-48b), the sum of the two integrals defined by (B41) will be given by

$$\frac{\sqrt{2\pi}e^{-\mathrm{i}\left(t-\frac{t}{2}log\left(\frac{t}{2\pi}\right)+\pi+\frac{2\pi^2}{t}\ldots\right)-\mathrm{i}\pi/4}}{\sqrt{t}(1-8\pi/t)^{\frac{1}{4}}}+O\left(\frac{1}{t}\right). \tag{B50}$$

Unfortunately the very small size of the error term in (B50) cannot be translated into an error estimate for $C(1,R)+D(1,R)$, because one must still consider the contribution to the integrals (B41) along the line $u=t/\pi R+\mathrm{i}\delta$ in the region $\delta<\beta_O\sqrt{2t/\pi}$. Although (B39) and (B40) are both negative in this region, so that their respective numerators are exponentially small, their approximate values at $\beta_O=(t\pi/8R^2)^{1/6}$ of $e^{-4/3}$ and $e^{-2/3}$ are still sufficiently large to potentially produce significant additions to (B50). An estimate on the upper bound of these contributions is calculated below.

First the integration range $\delta\in\left[\approx t^{1/4},\beta_O\sqrt{2t/\pi}\right]$ is divided into $M$ binary sections so that $\beta_O\sqrt{2t/\pi}/2^M=t^{1/4}$ and $M=INT\left(log\left(\beta_O t^{1/4}\sqrt{2/\pi}\right)/log(2)\right)$. The choice of $t^{1/4}$ as a lower bound ensures that the numerator phases given by (B39) and (B40) are now an insignificant $e^{-O(1)t^{1/4}}$ in size, whilst it lies just above the cut-off value $\approx 0.89t^{1/4}$ established by (B48), so ensuring that the asymptotic results (B42-44) are still valid. Utilising (B43), initial estimates of the upper bounds on the moduli of the integrals (B41) across the integration range $\delta\in\left[\beta_O\sqrt{2t/\pi}/2^M,\beta_O\sqrt{2t/\pi}\right]$ are given by

$$\left|\frac{-\mathrm{i}e^{-\mathrm{i}\left(t-\frac{t}{2}log\left(\frac{t}{2\pi}\right)+\pi+\frac{2\pi^2}{t}\right)}}{(1-8\pi/t)^{1/4}}\int_{\frac{\beta_O\sqrt{2t}}{2^M\sqrt{\pi}}}^{\frac{\beta_O\sqrt{2t}}{\sqrt{\pi}}}\left(\frac{1}{\delta^2}+\frac{\mathrm{i}}{\delta^3}\left(\frac{2t}{\pi R}\pm\frac{1/2}{(1-8\pi/t)}\right)\right)exp\left[\frac{A_\pm}{\delta}+\frac{\mathrm{i}B_\pm}{\delta^2}+\frac{C_\pm}{\delta^3}\right]d\delta\right|$$

$$\leq\frac{1}{(1-8\pi/t)^{1/4}}\sum_{n=0}^{M-1}\left|\int_{\frac{\beta_O\sqrt{2t}}{2^{n+1}\sqrt{\pi}}}^{\frac{\beta_O\sqrt{2t}}{2^n\sqrt{\pi}}}\left(\frac{1}{\delta^2}+\frac{\mathrm{i}}{\delta^3}\left(\frac{2t}{\pi R}\pm\frac{1/2}{(1-8\pi/t)}\right)\right)exp\left[\frac{A_\pm}{\delta}+\frac{\mathrm{i}B_\pm}{\delta^2}+\frac{C_\pm}{\delta^3}\right]d\delta\right|.\quad\text{(B51)}$$

Now in each sub integral, the real part of the phase is bounded above by its value at its upper integration limit. With $R\approx_p R_{max}\approx t/\pi$ and $\beta_O\approx_p\sqrt{\pi/2}/t^{1/6}$ these phases are approximately $e^{-4(2^{3n})/3}$ and $e^{-2(2^{3n})/3}$ from (B39) and (B40) respectively. Hence (B51) is strictly less than

$$\frac{1}{(1-8\pi/t)^{1/4}}\sum_{n=0}^{M-1}e^{-4\,or\,2(2^{3n})/3}\left|\int_{\frac{t^{1/3}}{2^{n+1}}}^{\frac{t^{1/3}}{2^n}}\frac{exp(\mathrm{i}t/2\delta^2)}{\delta^2}d\delta\right|$$

$$=\frac{1}{(1-8\pi/t)^{1/4}}\sum_{n=0}^{M-1}e^{-4\,or\,2(2^{3n})/3}\left|\sqrt{\frac{\pi}{-2\mathrm{i}t}}\left[erf\left\{\sqrt{\frac{-\mathrm{i}4^{n+1}t^{1/3}}{2}}\right\}-erf\left\{\sqrt{\frac{-\mathrm{i}4^{n}t^{1/3}}{2}}\right\}\right]\right|$$

$$\approx\frac{1}{(1-8\pi/t)^{1/4}}\sum_{n=0}^{M-1}e^{-4\,or\,2(2^{3n})/3}\left|\sqrt{\frac{\pi}{-2\mathrm{i}t}}\left[\frac{\left\{2e^{\mathrm{i}4^nt^{1/3}/2}-e^{\mathrm{i}4^{n+1}t^{1/3}/2}\right\}}{t^{1/6}\sqrt{-\mathrm{i}2\pi4^n}}\right]\right|,$$

$$\leq\frac{1}{(1-8\pi/t)^{1/4}}\frac{3}{2t^{2/3}}\sum_{n=0}^{M-1}\frac{e^{-4\,or\,2(2^{3n})/3}}{4^{n/2}},\quad\text{(B52)}$$

since $\left|2e^{i4^n t^{1/3}/2} - e^{i4^{n+1}t^{1/3}/2}\right| \leq 3$. The sums in (B52) converge extremely rapidly to 0.5158 and 0.2636 to four decimal places respectively. Consequently, combining the results (B5), (B19), (B50) and (B52), one deduces that

$$C(1,R) + D(1,R) = [\text{eqs. (B5)} + \text{(B19)}]_{N=1} + \frac{\sqrt{2\pi}e^{-i\left(t-\frac{t}{2}log\left(\frac{t}{2\pi}\right)+\pi+\frac{2\pi^2}{t}\right)-i\pi/4}}{\sqrt{t}(1-8\pi/t)^{\frac{1}{4}}} + \varepsilon, \quad \text{(B53)}$$

$$\text{where} \qquad |\varepsilon| < \frac{1}{(1-8\pi/t)^{1/4}}\frac{3\times(0.5158+0.2636)}{2t^{2/3}} = \frac{1.1691}{(1-8\pi/t)^{1/4}t^{2/3}}. \qquad \text{(B54)}$$

For $R \approx_p R_{max}$ and $pc_-$ as defined in *Lemma B1.2*, then (B5) and (B19) are terms of $O(t^{-1})$ and $O(t^{-13/12})$ respectively when $N = 1$. Hence the combination of (B53-54) completes the proof of *Lemma B1.6*.

B1.7 *Summary of results for $C(N,R) + D(N,R)$.*

In (24) one actually requires an estimate of the sum of the two integrals $C(N,R)$ and $D(N,R)$. Now examination of equations (B8) and (B22) reveals that they of opposite sign, and so the contributions of the terms of $O(\{t|1-N|\}^{-1})$ in *Lemmas B1.1 & 1.2* cancel exactly. So in general the sum of these two integrals is given by a combination of the results (B5), (B18) and (B19). Specifically if $N = 0$, $N \in [2, N_{2\eta}]$ and $N > N_t$, an approximation for this sum is given by three main terms and their associated errors/corrections

$$C(N,R) + D(N,R) = \frac{-i}{t/2}\left\{\frac{exp[-it\{N(1+\pi R/t) - (log(pc_+) + 1/pc_+)/2\}]}{((1+\pi R/t)^2 - 8\pi/t)^{\frac{1}{4}}\left(a/2\sqrt{pc_+} - 2N\right)} + \right.$$

$$\left.\frac{exp[-it\{N(1-\pi R/t) - (log(pc_-) + 1/pc_-)/2\}]}{((1-\pi R/t)^2 - 8\pi/t)^{1/4}\left(2N - a/2\sqrt{pc_-}\right)}\right\} + \sqrt{\frac{2\pi}{Nt}}e^{-i\pi/4}exp\left(i\frac{t}{2}log\left(\frac{t}{2\pi N^2}\right) - it - i\pi N^2\right) +$$

$$O\left(\frac{1}{t^2\left(a/2\sqrt{pc_+} - 2N\right)^2}\right) + O\left(\frac{(1-\pi R/t)[(1-\pi R/t)^2 - 8\pi/t]^{-1}}{\left(t\left(2N - a/2\sqrt{pc_-}\right)\right)^2}\right) + O_1\left(\frac{2^{3/2}5\sqrt{\pi}(1/N + 2\pi N/t)^2}{16(tN)^{3/2}(1/N - 2\pi N/t)^3}\right).$$

(B55)

The third main term and its associated error are *always excluded* in the cases $N = 0$ and $N > N_t$, irrespective of the choice of $R$ (see *Lemma B1.2*). When $N = 1$, these two terms should be replaced by the estimates (B53-4). The relative sizes of the various error terms are guaranteed to be less than 2% of their respective main terms for the specified $N$ values.

When $N \in (N_{2\eta}, N_t]$, specific modifications to (B55) must be made. The first main term and its error remain unchanged. The second main term (from the integration of $D(N,R)$) can be replaced by the upper bound $0.38 \times 3.55/t^{3/4}$ from *Lemma B1.2*. The third term is unchanged if $N \in (N_{2\eta}, N_\mu]$, vanishes if $N > N_\eta$ and for $N \in (N_\mu, N_\eta]$ the upper bound of $0.64 \times 1.728\sqrt{2\pi/Nt}$ established in *Lemma B1.2* can be substituted. If the radius is reduced so that $R = R_{max} - R_t$ with $R_t \sim O\left(t^{0-1/2}\right)$, this will rule out saddle point contributions for $N$ lying within $t^{-1/6}$ of $N_t$, and mean that the lower limit of integration in (B9) no longer lies close to the singularity of the integrand. Consequently (B55) will provide a good approximation *for all* $N$ in this instance (with the third term excluded for $N = 0$ and $N > Int\left[a/4\sqrt{pc_-}\right]$), as the error terms will always remain (in practice very much) below 2% of the main terms. The terms of (B55) with their modifications, forms the basis of the proof of the *Main Theorem* discussed in Sections 2.4-2.6.

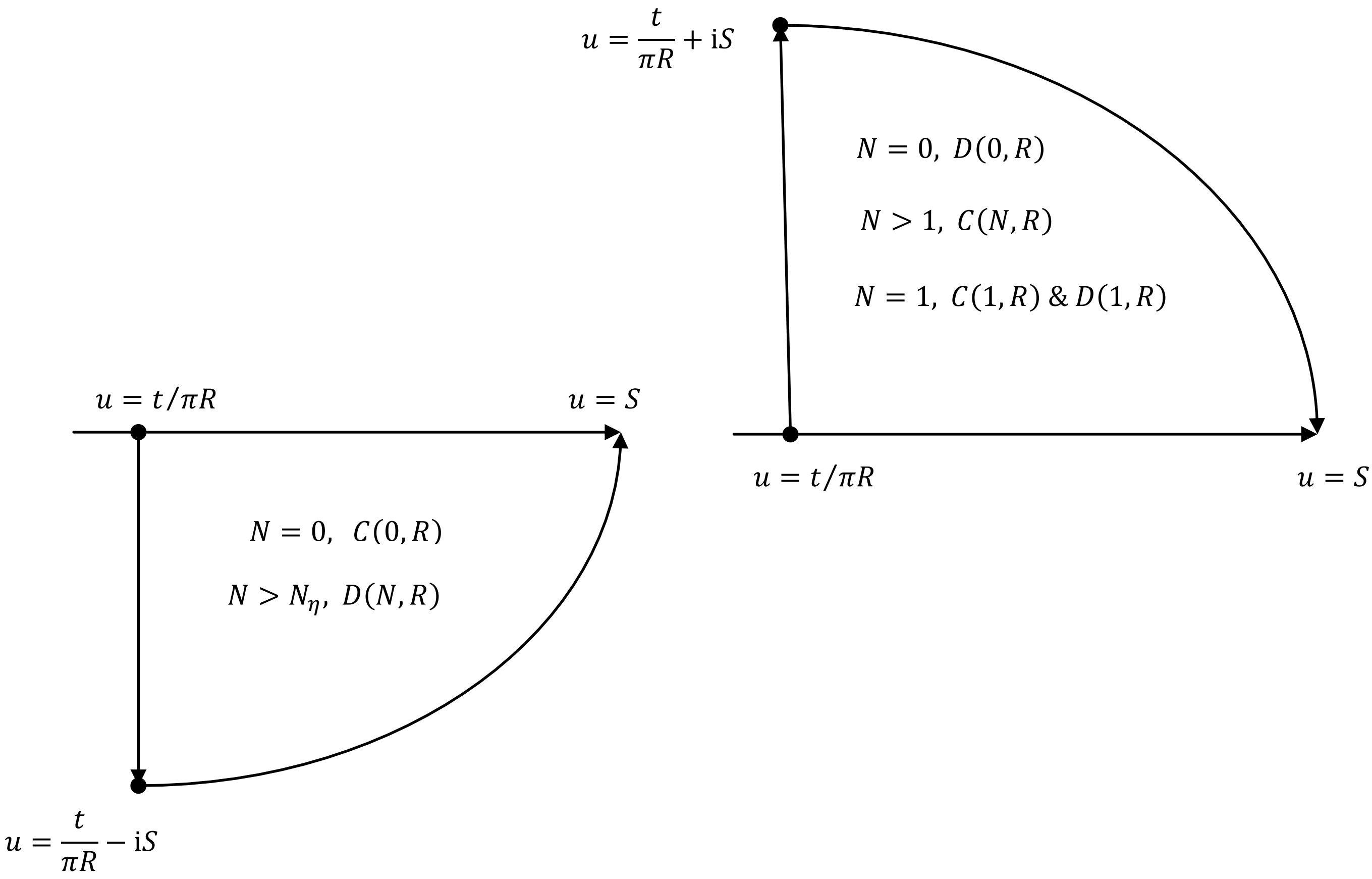

Figure B1. Contour integration paths for estimation of the integrals $C(N,R)$ and $D(N,R)$for the various ranges of $N$ as shown.

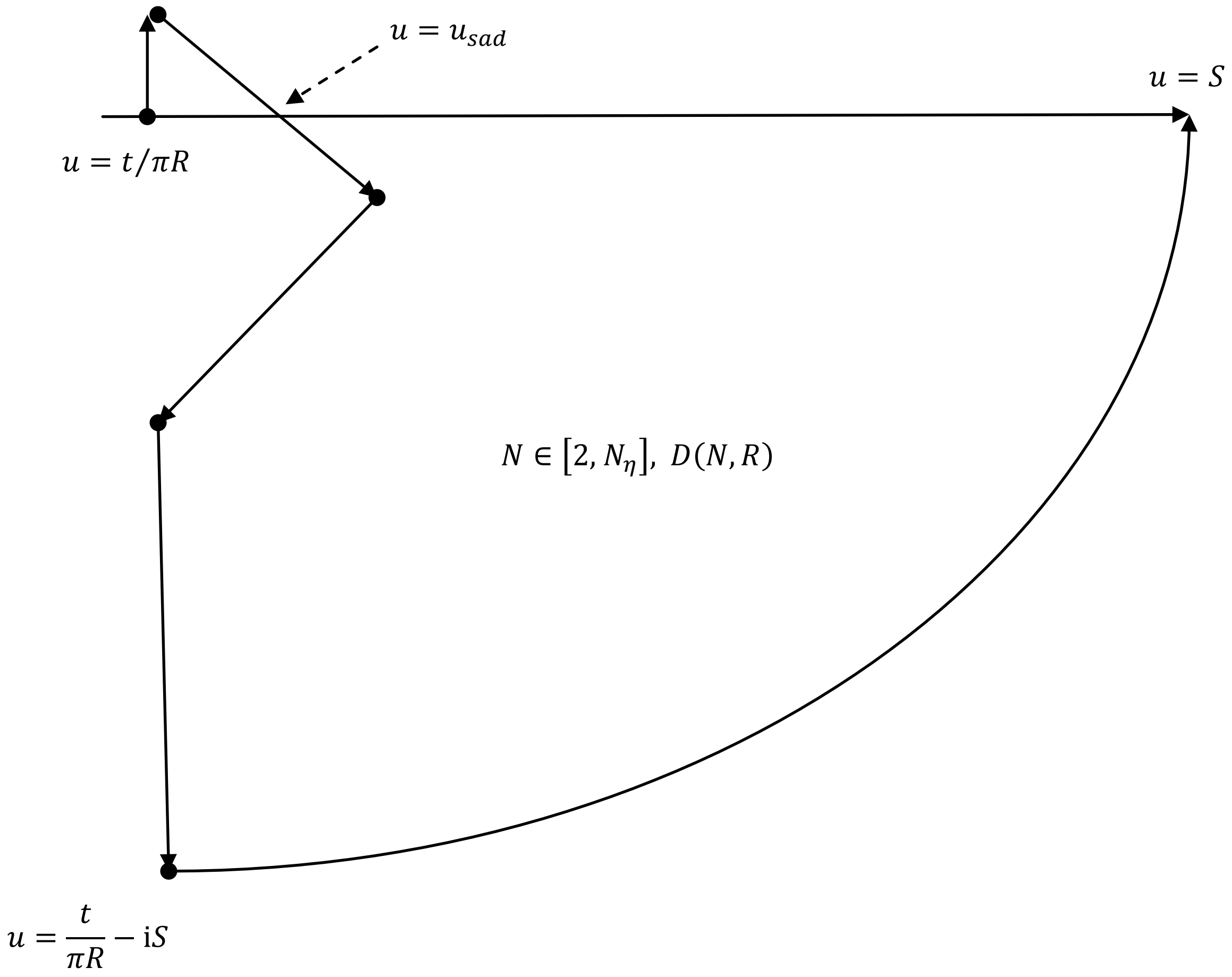

Figure B2. The contour integration path for estimation of the integral $D(N,R)$ when $N \in [2, N_\eta]$, which must pass through the saddle point defined by equation B15a as shown.

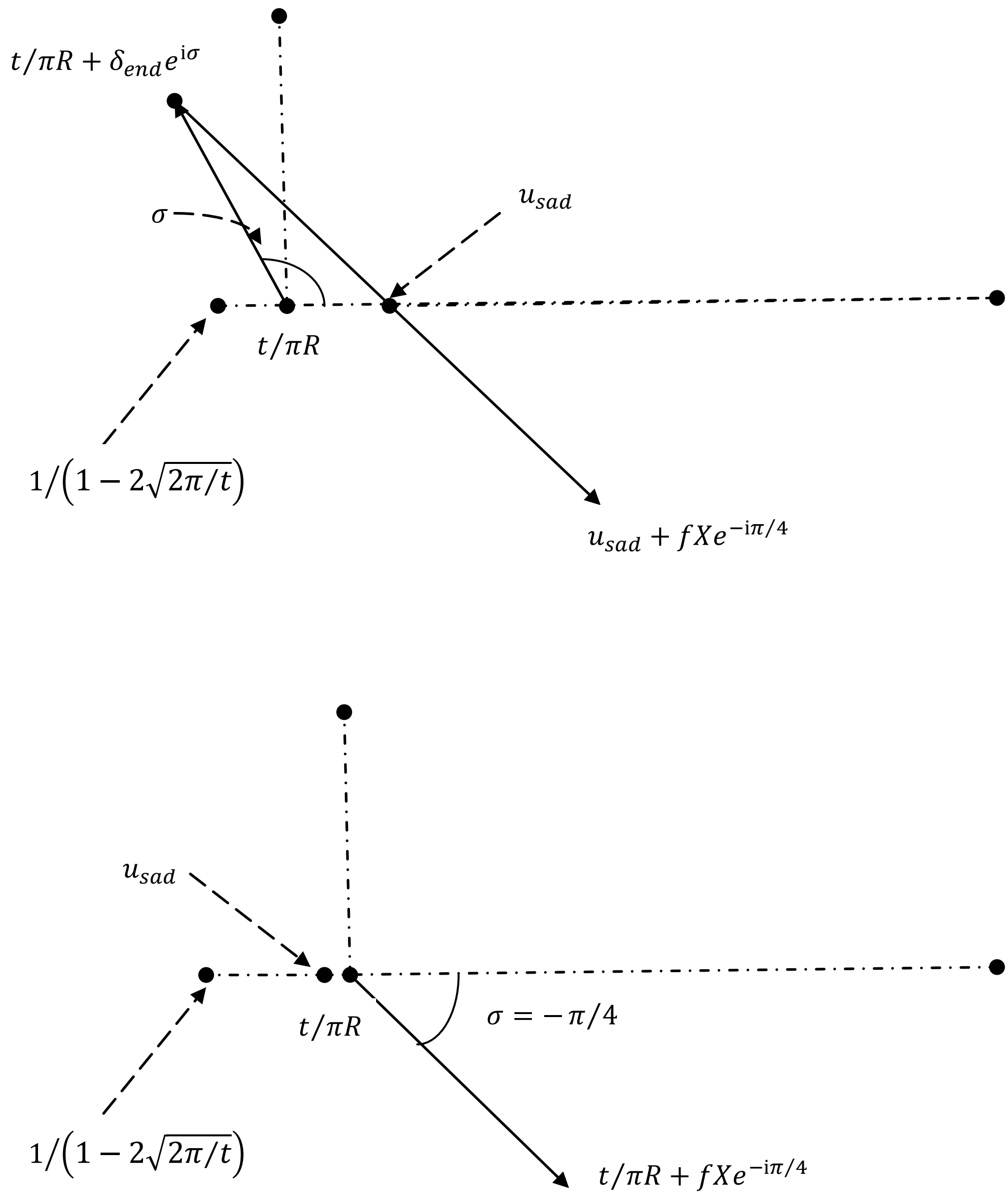

Figure B3. Modifications to the contour shown in Figure B2 in the case when $N \approx N_t$, which means that $u_{sad}$ and $u = t/\pi R$ are lie *extremely* close together (distances here exaggerated for clarity). When $r_N \in (\eta, 1.2)$ (Figure B3a) and $u_{sad} > t/\pi R$ the integration begins along a line an angle $\sigma$ to the horizontal to a point $u = t/\pi R + \delta_{end} e^{i\sigma}$ given by equation B35c. When $r_N \in (0.7, \eta)$ (Figure B3b) and $u_{sad} < t/\pi R$ the integration begins along a line at an angle of $-\pi/4$ to a point $u = t/\pi R + \delta_{end} e^{-i\pi/4}$

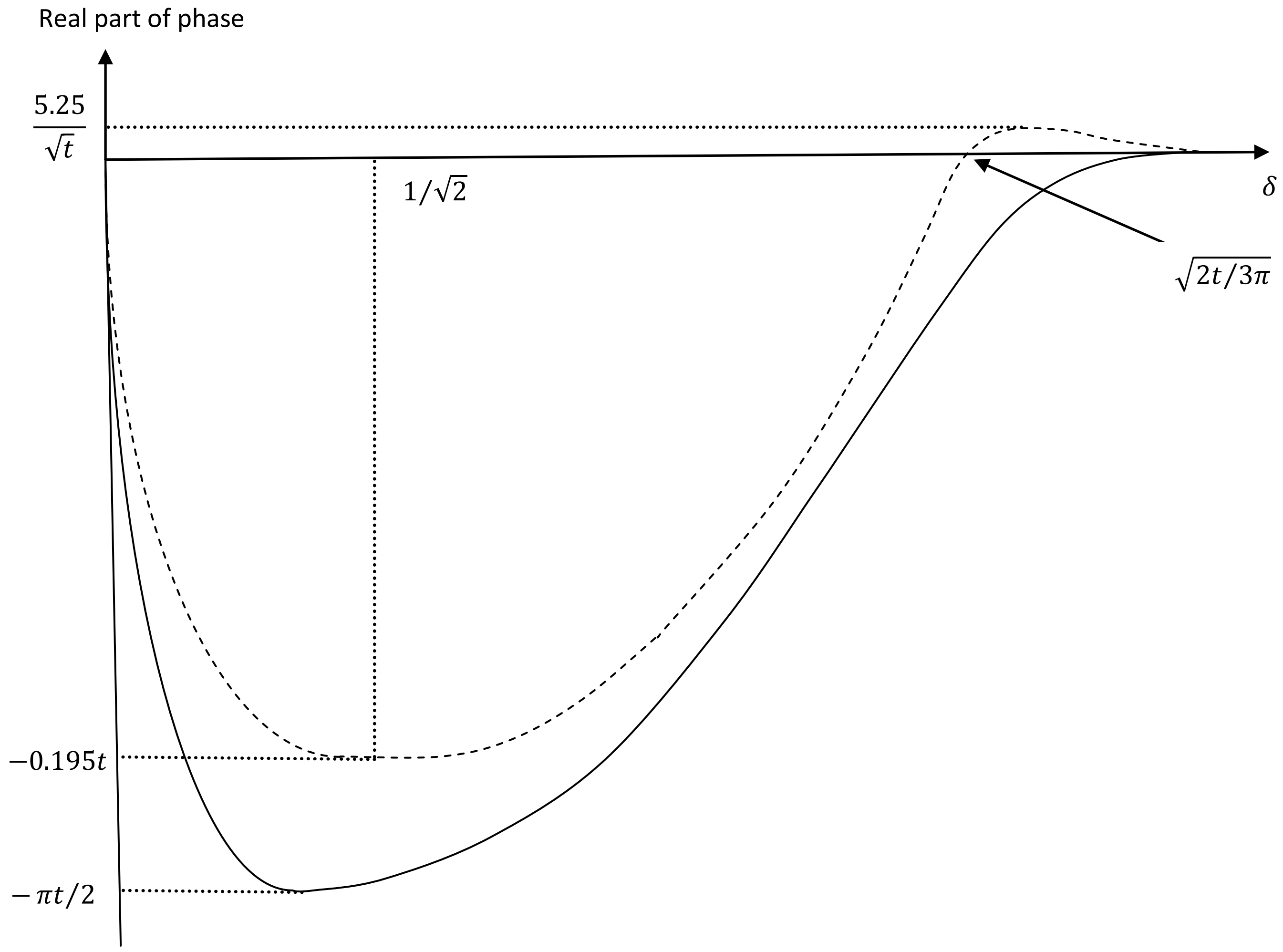

Figure B4. Schematic of the changes to the real part of the phase of the numerator in integrals $C(1,R)$ and $D(1,R)$ observed with $\delta$ when moving along the line $u = t/\pi R + \mathrm{i}\delta$. The solid line —— applies for $D(1,R)$ and the dashed line - - - - - for $C(1,R)$.

# Appendix C. The Imaginary Part of $e^{\mathrm{i}\theta(t)}\times$RSI.

## C1.1 *Introduction*

This appendix is concerned with development of an asymptotic formula for the combined sum of the two largest terms in (A38), developed from the proposed expansion (A33) for the Riemann-Siegel Integral postulated at the end of Section A2. Hence one seeks an estimate of

$$\mathrm{RSI}=-\frac{e^{\frac{\pi t}{2}+\mathrm{i}3\pi/8-\mathrm{i}\theta(t)}}{2\mathrm{i}(1+e^{-2\pi t})}\sum_{\substack{\alpha=1\\ odd}}^{\infty}(-1)^{(\alpha-1)/2}\,e^{-\mathrm{i}\pi\alpha/2}\left[\left(\frac{t\pi}{2}\right)^{-\frac{1}{4}}\sqrt{2\pi}\,e^{[1/32t^2+\varpi]}\Phi\left(\frac{1}{4}-\frac{\mathrm{i}t}{2},\frac{1}{2};\frac{\pi\alpha^2\mathrm{i}}{4}\right)\right.$$
$$\left.-\frac{2\alpha}{\sqrt{2\pi}}\left(\frac{t\pi^5}{2}\right)^{1/4}e^{-[1/32t^2+\varpi]}\Phi\left(\frac{3}{4}-\frac{\mathrm{i}t}{2},\frac{3}{2};\frac{\pi\alpha^2\mathrm{i}}{4}\right)\right],\quad (\mathrm{C1})$$

which is valid for large $t$. In the process it will be established that $e^{\mathrm{i}\theta(t)}\times(\mathrm{C1})$ is purely imaginary, and consequently plays no direct role in calculation of the zeta function, as this depends upon the real part of $e^{\mathrm{i}\theta(t)}$RSI.

## C1.2 *An asymptotic approximation for the confluent hypergeometric function*

In order to make progress one must be able establish an asymptotic approximation for $\Phi(a,c;z)$ valid as $|a|\to\infty$. Fortunately just such an approximation was developed by [35] based on a series of exact expansions in terms of Bessel functions of the first kind. The key results are summarised in [27] Chapter 4.8. The relevant expansion required in this instance is

$$\Phi(a,c;z)=\Gamma(c)e^{z/2}\sum_{m=0}^{\infty}D_m\left(\frac{z}{2}\right)^m(kz)^{-(c+m-1)/2}J_{c+m-1}\big(2\sqrt{kz}\big),\qquad (\mathrm{C2})$$

where $k=c/2-a$. The first four of the coefficients $D_m$ are given by $D_0=1,\ D_1=0, D_2=c/2$ and $D_3=-2k/3$. The remainder can be generated by the following recurrence relation

$$(m+1)D_{m+1}=(m+c-1)D_{m-1}-2kD_{m-2}\qquad m\geq 2.\qquad (\mathrm{C3})$$

Table CI lists the first thirteen such coefficients relevant to the two confluent hypergeometric functions that appear in (C1). As one can see, (C3) generates complicated odd and even polynomials in the parameter $t$, with no obvious pattern to the multiplicative coefficients (one does notice that the leading term of $D_{3n}=(-\mathrm{i}t)^n/(3^n n!)$, $D_{3n-1}=O(t^{n-1})$ and $D_{3n-2}=O(t^{n-2})$, whilst it is elementary to prove that $D_{2n}=real$ and $D_{2n+1}=imaginary$).

As $c=1/2$ and $3/2$, the Bessel functions needed in (C2) are for half integer order. A suitable generating function is given by [16], eq. 8.463

$$J_{-\frac{1}{2}}(z) = \sqrt{\frac{2}{\pi z}}\, cos(z), \qquad J_{n+\frac{1}{2}}(z) = (-1)^n z^{n+\frac{1}{2}} \sqrt{\frac{2}{\pi}} \frac{d^n}{(zdz^n)}\left(\frac{sin(z)}{z}\right), \ \ n \geq 0. \qquad \text{(C4)}$$

Now in (C1), $z = \mathrm{i}\pi\alpha^2/4$, which means that $kz = -\pi\alpha^2 t/8$, $\Rightarrow 2\sqrt{kz} = +\mathrm{i}(\pi\alpha^2 t/2)^{1/2}$ for both of the required confluent hypergeometric functions. Consequently in (C2 & C4) $cos\left(2\sqrt{kz}\right) = cosh(\pi\alpha^2 t/2)^{1/2}$, to be denoted by $CO$, and $sin\left(2\sqrt{kz}\right) = \mathrm{i}sinh(\pi\alpha^2 t/2)^{1/2}$, to be denoted by i$SI$. Now because $z = imaginary$, it means that in (C2) both $D_{2n} \times (z/2)^{2n} = real$ and $D_{2n+1} \times (z/2)^{2n+1} = real$. In addition one has from (C4)

$$J_{n-\frac{1}{2}}(z) = \sqrt{\frac{2}{\pi z}} \begin{cases} even\ function\ of\ z\ for\ n = even \\ odd\ function\ of\ z\ for\ n = odd \end{cases} \Rightarrow J_{n-\frac{1}{2}}(\mathrm{i}x) = \sqrt{\frac{2}{\pi \mathrm{i}x}} \begin{cases} real\ number\ for\ n = even \\ imaginary\ number\ for\ n = odd \end{cases}$$

(C5)

The latter observation implies that both the terms $(kz)^{\frac{1}{4}-m/2} J_{m-1/2}\left(2\sqrt{kz}\right)$ for $c = 1/2$, and $(kz)^{-\frac{1}{4}-m/2} J_{m+1/2}\left(2\sqrt{kz}\right)$ for $c = 3/2$, are *real* for all $m \geq 0$. Consequently the sum in (C2) is simply a sum of real terms. Substituting results $(\mathrm{C2} - \mathrm{C4})$ into (C1) and expanding up to the first four terms ($D_0$ to $D_3$) in detail, one finds after considerable manipulation that

$$\Phi\left(\frac{1}{4} - \frac{\mathrm{i}t}{2}, \frac{1}{2}; \frac{\pi\alpha^2 \mathrm{i}}{4}\right) = e^{\mathrm{i}\pi\alpha^2/8} CO\left[1 - \frac{\sqrt{2}}{96}\frac{(\pi\alpha^2)^{3/2}}{\sqrt{t}} tanh(\pi\alpha^2 t/2)^{1/2} + O\left(\frac{\alpha^{3m}}{t^{m/2}},\ m = 2,3,..\right)\right],$$

$$\Phi\left(\frac{3}{4} - \frac{\mathrm{i}t}{2}, \frac{3}{2}; \frac{\pi\alpha^2 \mathrm{i}}{4}\right) = e^{\mathrm{i}\pi\alpha^2/8} SI\left[\frac{1}{(\pi\alpha^2 t/2)^{1/2}} - \frac{\pi\alpha^2}{48t} coth(\pi\alpha^2 t/2)^{1/2} + O\left(\frac{\alpha^{3m-1}}{t^{(m+1)/2}},\ m = 2,3,..\right)\right].$$

(C6a, b)

Notice that in (C6) both expressions are of the form $e^{\mathrm{i}\pi\alpha^2/8} \times real$. Note too that the $O(\quad)$ terms are just representative of the largest powers of $\alpha$, and there are intermediate terms of $\alpha^{3m-1}$ and $\alpha^{3m-2}$ for each $m$ value (dependant on the period three cycle seen in Table CI).

C1.3 *The asymptotic approximation for the RSI when* $z = 1/2 + \mathrm{i}t$.

Substituting results (C6a, b) into (C1) gives

$$\mathrm{RSI} = -\frac{e^{\frac{\pi t}{2} + \mathrm{i}3\pi/8 - \mathrm{i}\theta(t)}}{2\mathrm{i}(1 + e^{-2\pi t})} \sum_{\substack{\alpha=1 \\ odd}}^{\infty} (-1)^{(\alpha-1)/2}\, e^{-\frac{\mathrm{i}\pi\alpha}{2}} e^{\mathrm{i}\pi\alpha^2/8} \left[\left(\frac{2^3\pi}{t}\right)^{\frac{1}{4}} \left\{e^{[1/32t^2 + \varpi]} CO - e^{-[1/32t^2 + ...]} SI\right\}\right.$$

$$\left. + \left(\frac{2\pi^7}{t^3}\right)^{\frac{1}{4}} \frac{\alpha^3}{48} \left\{e^{-[1/32t^2 + \varpi]} CO - e^{[1/32t^2 + ...]} SI\right\} + O\left(\frac{\alpha^{3m}}{t^{(2m+1)/4}} \times CO\ or\ SI,\ m = 2,3,..\right)\right], \quad \text{(C7)}$$

where again the $O(\quad)$ terms are just representative of the largest powers of $\alpha$.

Two points are of crucial importance. First for (C7) give rise to a *convergent* sum, which is inherent in the postulate that expansion (A33) (which underlies C7) is indeed the correct interpretation of the *convergent* RSI, then it is apparent that the $e^{+\alpha\sqrt{\pi t}/\sqrt{2}}$ factors in both $cosh(\pi\alpha^2 t/2)^{\frac{1}{2}}$ and $sinh(\pi\alpha^2 t/2)^{\frac{1}{2}}$ must *cancel*. Now to $O(1)$, both $e^{[1/32t^2+\varpi]} = e^{-[1/32t^2+\varpi]} = 1$, and this cancellation is then obvious for the first two terms in (C7). That this cancellation should continue indefinitely for all orders of $t$ can be inferred from the fact that the presence of any $e^{+\alpha\sqrt{\pi t}/\sqrt{2}}$ factor in (C7) would result in the sum over $\alpha$ diverging to infinity, completely invalidating expansion (A33) as representative of the RSI. And indeed for all the terms up to of $O\left(t^{-\frac{9}{4}}\right)$ (which is as far as the author has checked) explicit calculation shows that this exact cancellation does indeed continue. (The computation grows ever more lengthy and the $O\left(t^{-\frac{9}{4}}\right)$ calculation requires one to expand the $e^{[1/32t^2+\varpi]}$ and $e^{-[1/32t^2+\varpi]}$ terms as $1 \pm 1/32t^2$ respectively.) It is the author's *assertion* that this cancellation must continue essentially indefinitely. However, how one might go about formally *proving* this is not immediately obvious. Clearly for the cancellation of terms to continue in the way described, there must be some deeper connection between the factors appearing in Table CI and the terms appearing in Stirling's expansion for $Re\big(Log\Gamma(1/4 + \mathrm{i}\,t/2)\big)$, as given by (A35). Another observation is that if one can prove that (C1) is bounded above then this must imply the desired cancellation, for otherwise the sum would diverge to infinity. (These points are revisited in Section A3.8 in connection with the real part of $e^{\mathrm{i}\theta(t)}$RSI.)

However, very strong evidence that (A33) is indeed the correct interpretation of the RSI can be found by comparing the asymptotic formula that results by summing (C7), with some direct numerical computations of the RSI. Writing $\alpha = 2k+1, k = 0,1,2,\ldots$ one can see that the complex variables in (C7) reduce to

$$\frac{e^{\mathrm{i}3\pi/8-\mathrm{i}\theta(t)}}{\mathrm{i}} e^{-\frac{\mathrm{i}\pi\alpha}{2}} e^{\mathrm{i}\pi\alpha^2/8} = \frac{e^{\mathrm{i}3\pi/8-\mathrm{i}\theta(t)}}{\mathrm{i}} e^{\mathrm{i}\pi k(k-1)/2-\mathrm{i}\pi/2+\mathrm{i}\pi/8} = \pm\mathrm{i}e^{-\mathrm{i}\theta(t)}, \qquad \text{(C8)}$$

for successive pairs of $k$ values. Consequently $e^{i\theta(t)} \times$ (C7) is just an imaginary number. One can evaluate (C7) by noting that the cancellation of the $e^{+\alpha\sqrt{\pi t}/\sqrt{2}}$ factors in the terms of $O\left(t^{-1/4}\right)$ and $O\left(t^{-3/4}\right)$ means the sum reduces simply to

$$\mathrm{RSI} = -\frac{e^{\frac{\pi t}{2}+\mathrm{i}3\pi/8-\mathrm{i}\theta(t)}}{2\mathrm{i}(1+e^{-2\pi t})} \sum_{\substack{\alpha=1\\ odd}}^{\infty} (-1)^{(\alpha-1)/2} e^{-\frac{\mathrm{i}\pi\alpha}{2}} e^{\mathrm{i}\pi\alpha^2/8} e^{-\alpha\sqrt{\pi t}/\sqrt{2}} \left[\left(\frac{2^3\pi}{t}\right)^{\frac{1}{4}} + \left(\frac{2\pi^7}{t^3}\right)^{\frac{1}{4}} \frac{\alpha^3}{48} + O\left(\frac{\alpha^{3m}}{t^{(2m+1)/4}}\right)\right],$$

(C9)

which is just a succession of (convergent) geometric series. The result of the first such series can be derived as follows:

$$\sum_{\substack{\alpha=1\\ odd}}^{\infty}(-1)^{(\alpha-1)/2}\, e^{-\frac{\mathrm{i}\pi\alpha}{2}} e^{\mathrm{i}\pi\alpha^2/8} e^{-\alpha\sqrt{\pi t/2}} = -\mathrm{i}e^{\mathrm{i}\pi/8} e^{-\sqrt{\pi t}/\sqrt{2}} \sum_{k=0}^{\infty} e^{\frac{\mathrm{i}k(k+1)\pi}{2}} e^{-2k\sqrt{\pi t/2}},$$

$$=-\mathrm{i}e^{\mathrm{i}\pi/8} e^{-\sqrt{\pi t}/\sqrt{2}}\left[\sum_{k=0 \text{ even}} e^{k\left(\mathrm{i}\pi/2-2\sqrt{\pi t/2}\right)} + \mathrm{i}\sum_{k=1 \text{ odd}} e^{k\left(\mathrm{i}\pi/2-2\sqrt{\pi t/2}\right)}\right],$$

$$=-\mathrm{i}e^{\mathrm{i}\pi/8} e^{-\sqrt{\pi t}/\sqrt{2}}\left[\frac{1}{1-e^{2\left(\mathrm{i}\pi/2-2\sqrt{\pi t/2}\right)}} + \mathrm{i}\,\frac{e^{\left(\mathrm{i}\pi/2-2\sqrt{\pi t/2}\right)}}{1-e^{2\left(\mathrm{i}\pi/2-2\sqrt{\pi t/2}\right)}}\right],$$

$$=\frac{-\mathrm{i}e^{\mathrm{i}\pi/8} e^{-\sqrt{\pi t}/\sqrt{2}}}{1+e^{-2\sqrt{2\pi t}}}\left(1+\mathrm{i}\times \mathrm{i}e^{-\sqrt{2\pi t}}\right) = -\mathrm{i}e^{\mathrm{i}\pi/8}\,\frac{sinh\left(\sqrt{\pi t/2}\right)}{cosh\left(\sqrt{2\pi t}\right)}. \qquad \text{(C10)}$$

By means of a similar but more complicated argument involving the third derivative with respect to $\sqrt{\pi t/2}$, one can sum the second series in (C9) involving $\alpha^3$(but the details are omitted‡). In conclusion one finds that

$$\mathrm{RSI}\left(\frac{1}{2}+\mathrm{i}t\right) = \frac{\mathrm{i}e^{\pi t/2} e^{-\mathrm{i}\theta(t)}}{(1+e^{-2\pi t})}\left[\left(\frac{\pi}{2t}\right)^{1/4}\frac{sinh\left(\sqrt{\pi t/2}\right)}{cosh\left(\sqrt{2\pi t}\right)} - \frac{\pi^{7/4}}{48(2t)^{3/4}}\frac{1}{C_2}\left\{2S_1T_2\left(\frac{24}{(C_2)^2}-7\right)+C_1\left(1+\frac{12((S_2)^2-1)}{(C_2)^2}\right)\right\} + O\left(e^{-\sqrt{\pi t/2}}t^{-5/4}\right)\right],$$

(C11)

where $S_1 = sinh\left(\sqrt{\pi t/2}\right)$, $C_1 = cosh\left(\sqrt{\pi t/2}\right)$, $S_2 = sinh\left(\sqrt{2\pi t}\right)$, $C_2 = cosh\left(\sqrt{2\pi t}\right)$ and $T_2 = S_2/C_2$. Table CII shows a comparison of the values computed from (C11) with numerical calculations of the RSI numerically along the line $0 \searrow 1$. Even for these relatively modest values of $t$ the agreement is strikingly good. In addition, the decline in the relative error is characteristic of the omission of a term of $O(t^{-1})$ smaller than the dominant one, consistent with the analysis. Such agreement can only arise if the postulated expansion (A33) is indeed the correct interpretation of the RSI and indicative (although not a formal proof) that the assertion concerning the cancellation of the $e^{+\alpha\sqrt{\pi t}/\sqrt{2}}$ terms must be true. (Incidentally the fact that the rough order of magnitude of A11 is $e^{\pi t/2}$, shows that mismatch of the asymptotic approximations for the $\mathrm{RSI}_-$ and $\mathrm{RSI}_+$ discussed in section A2.2 occurs in the vicinity of $q = -1/\sqrt{2}$, the point at which the contour crosses the imaginary axis and the arctan function is undefined.) Now in principle one can generate more and more terms in the expansion (C11) to yield approximations for the $Im\{e^{\mathrm{i}\theta(t)} \times \mathrm{RSI}\}$ to arbitrary accuracy. But any such calculations are irrelevant to the zeta function, which depends upon the $Re\{e^{\mathrm{i}\theta(t)} \times \mathrm{RSI}\}$ given by the *third largest* term in equation (A38).

| Coefficient $D_m$ | $a = 1/4 - \mathrm{i}t/2, c = 1/2$ $k = \mathrm{i}t/2$ | $a = 3/4 - \mathrm{i}t/2, c = 3/2$ $k = \mathrm{i}t/2$ |
|---|---|---|
| $D_0$ | $1$ | $1$ |
| $D_1$ | $0$ | $0$ |
| $D_2$ | $1/4$ | $3/4$ |
| $D_3$ | $-\mathrm{i}t/3$ | $-\mathrm{i}t/3$ |
| $D_4$ | $5/32$ | $21/32$ |
| $D_5$ | $-\mathrm{i}17t/60$ | $-\mathrm{i}9t/20$ |
| $D_6$ | $15/128 - t^2/18$ | $77/128 - t^2/18$ |
| $D_7$ | $-\mathrm{i}823t/3360$ | $-\mathrm{i}573t/1120$ |
| $D_8$ | $195/2048 - 29t^2/360$ | $1155/2048 - 13t^2/120$ |
| $D_9$ | $-\mathrm{i}1751t/8064 + \mathrm{i}t^3/162$ | $-\mathrm{i}22177t/40320 + \mathrm{i}t^3/162$ |
| $D_{10}$ | $663/8192 - 9371t^2/100800$ | $4389/8192 - 5177t^2/33600$ |
| $D_{11}$ | $-\mathrm{i}278437t/1419264 + \mathrm{i}41t^3/3240$ | $-\mathrm{i}194741t/337920 + \mathrm{i}17t^3/1080$ |
| $D_{12}$ | $4641/65536 - 120283t^2/1209600 + t^4/1944$ | $33649/65536 - 234049t^2/1209600 + t^4/1944$ |

Table CI. List of the first thirteen polynomials generated by the recursion relation (C3) which appear in the asymptotic expansions for $\Phi(a, c; z)$.

| $t$ | Numerical calculation of $\mathrm{RSI}(1/2 + \mathrm{i}t)$ (to a maximum of 7 significant figures). | Asymptotic approximation, calculated from equation (A11). | Modulus of relative error |
|---|---|---|---|
| 10 | $-6.138923 \times 10^3 - \mathrm{i}8.223933 \times 10^4$ | $-6.058749 \times 10^3 - \mathrm{i}8.115024 \times 10^4$ | $1.32 \times 10^{-2}$ |
| 20 | $8.148139 \times 10^{10} + \mathrm{i}3.291390 \times 10^{10}$ | $8.098092 \times 10^{10} + \mathrm{i}3.271175 \times 10^{10}$ | $6.14 \times 10^{-3}$ |
| 30 | $1.456097 \times 10^{17} - \mathrm{i}3.009586 \times 10^{16}$ | $1.450347 \times 10^{17} - \mathrm{i}2.997701 \times 10^{16}$ | $3.95 \times 10^{-3}$ |
| 40 | $-2.518034 \times 10^{23} - \mathrm{i}1.917761 \times 10^{23}$ | $-2.510735 \times 10^{23} - \mathrm{i}1.912200 \times 10^{23}$ | $2.90 \times 10^{-3}$ |
| 50 | $7.552 \times 10^{29} + \mathrm{i}1.865 \times 10^{29}$ | $7.533717 \times 10^{29} + \mathrm{i}1.860975 \times 10^{29}$ | $2.41 \times 10^{-3}$ |

Table CII. Comparison of a sample of numerical calculations of the $\mathrm{RSI}(1/2 + \mathrm{i}t)$ with the corresponding asymptotic approximation given by equation (C11), for some modest values of $t$.